\documentclass[%
	a4paper,
	twoside,
	10 pt,
    toc=bibliography
]{article}%

    \usepackage[ngerman, english]{babel} 
    \usepackage{amsmath} 
    \usepackage{extarrows} 
    \usepackage{amsthm}    
    \usepackage[arrow, matrix, curve]{xy}   
    \usepackage{amssymb} 
    \usepackage[usenames, dvipsnames]{xcolor}
    \usepackage{minibox} 
    \usepackage{leftidx} 
    \usepackage{MnSymbol} 
    \usepackage{fancyhdr}			
	\usepackage{faktor}
	\usepackage{graphicx} 
	\usepackage{enumitem}
	\usepackage{mathrsfs}
	\usepackage{hyperref}
	\usepackage{soul} 
	\usepackage{makeidx}
	\makeindex
	\usepackage[refpage, noprefix]{nomencl}
	
	\makenomenclature
	\usepackage[style=alphabetic, backend=bibtex]{biblatex}
    \bibliography{literatur.bib}
	\usepackage[utf8]{inputenc}
    \usepackage[expansion=false]{microtype}
	
\title{Functoriality of\\ Moduli Spaces of Global $\mathbb G$-Shtukas}
\author{Paul Breutmann}
\date{}

	\newtheorem{satz}{Satz}[section]
	 \newtheorem{definition}[satz]{Definition}
	 \newtheorem{lemma}[satz]{Lemma}
	 \newtheorem{korollar}[satz]{Corollary}
	 \newtheorem{pr}[satz]{Proposition}
	 \newtheorem{theorem}[satz]{Theorem}
	 \theoremstyle{definition}
	 \newtheorem{bssp}[satz]{Example}
	 \newtheorem{remark}[satz]{Remark}
	 \newtheorem{axiom}{Axiom}
		\usepackage[left=2.7cm, right=2.3cm,top=1cm,bottom=2cm, includeheadfoot]{geometry}
	
		\makeatletter
		\renewcommand{\section}{%
    		\@startsection {section}{1}{\z@}%
         		{-3.5ex plus -1ex minus -.2ex}%
          	{2.3ex plus.2ex}%
                 	{\normalfont\Large\bfseries\centerline}}
	\makeatother
		\makeatletter
\renewcommand\subsubsection{\@startsection{subsubsection}{3}{\z@}{-3.25ex
    plus -1ex minus -.2ex}{1.5ex plus .2ex}{\bfseries\centerline}}
	\makeatother

	 	\newcommand\klg{\leqslant}   
	\newcommand\grg{\geqslant}    
	 
    \newcommand\ds{\displaystyle}
	\newcommand\N{{\mathbb N}} 
	\newcommand\Z{{\mathbb Z}} 

	\newcommand\F{{\mathbb F}}

	\newcommand\OO{{\mathbb O}}

	\newcommand\A{{\mathbb A}}

	\newcommand\p{{\mathfrak p}}
	\newcommand\oo{\mathcal{O}}

	\newcommand\bo{\hfill$\Box$}
    
    \newcommand\sortkey[1]{}
    \newcommand\sem[1]{\lsem #1 \rsem}

    \newcommand{\semr}[1]{(\mkern-4.5mu(#1)\mkern-4.5mu)}

\newcommand\sa{\nabla_n \mathscr H^1(C,\mathbb G)}

\newcommand{\san}[2]{\nabla_n^{#1} \mathscr H^1 (C,\mathbb G)^{#2} }

\newcommand{\uv}{{\underline v}}
\newcommand{\pv}{{v}}
\newcommand{\pw}{{w}}
\newcommand{\gtor}{{\mathcal G}}
\newcommand{\fl}{{\mathcal F\mathit l}}

	 \newlength\breite
	 
    \newcommand\isom{\xrightarrow{\,\smash{\raisebox{-0.60ex}    
               {\ensuremath{\scriptstyle\sim}}}\ }}

	 \newlength\proofeinzug
	\newcommand\prof[1]{
   		\textsl{Proof:} #1 \bo
	} 
	 
	\newcommand\df[2]{
	        \begin{definition}
            [#1]
            #2 
            \end{definition}
            }
	\newcommand\rem[2]{\begin{remark}
            [#1]#2
            \end{remark}
            }
    \newcommand\lem[2]{\begin{lemma}
            [#1]#2
            \end{lemma}
            }
    \newcommand\ko[2]{\begin{korollar}
            [#1]#2
            \end{korollar}
            }
    \newcommand\prop[2]{\begin{pr}
            [#1]#2
            \end{pr}
            }
    \newcommand\theo[2]{\begin{theorem}
            [#1]#2
            \end{theorem}
            }

    \newcommand\ger[1]{}
    \newcommand\gerr[1]{}   
\begin{document}
\pagestyle{plain}

\maketitle
\begin{abstract}
	Moduli spaces of global $\mathbb G$-shtukas play a crucial role in the Langlands program for function fields. We analyze their functoriality properties following a change of the curve and a change of the group scheme $\mathbb G$ under various aspects. In particular, we prove two finiteness results which are of interest in the study of stratifications of these moduli spaces and which potentially allow the formulation of an analog of the Andr\'{e}-Oort conjecture for global $\mathbb G$-shtukas.
\end{abstract}

\tableofcontents

\section{Introduction}

Global $\mathbb G$-shtukas are the function field analogue of abelian varieties. Their moduli spaces play a crucial role in the Langlands program for function fields. This article is concerned about functoriality properties of these moduli spaces in various aspects. Let us look in more detail. 
We choose a smooth, projective, geometrically irreducible curve $C$ over a finite field $\F_q$ with $q$ elements. Let $\mathbb G$ be a smooth, affine group scheme over $C$ and denote by $\sigma$ the $\F_q$-Frobenius of a scheme $S$ over $\F_q$. Then a global $\mathbb G$-shtuka $\underline\gtor=(\gtor, s_1, \dots, s_n, \tau_\gtor)$ over $S$ consists of a $\mathbb G$-torsor $\gtor$ over $C_S:=C\times_{\F_q}S$, $n$ sections $s_i:S\to C$ called legs and an isomorphism $\tau_\gtor:\sigma^\star\gtor|_{C_S\backslash{\cup_i\Gamma_{s_i}}}\to \gtor|_{C_S\backslash{\cup_i\Gamma_{s_i}}}$ outside the union of the graphs $\Gamma_{s_i}$ of the $s_i$. 
The precise definitions of all the notations used in this article are given in the second section. The stack whose $S$-valued points parametrize the global $\mathbb G$-shtukas over $S$ with $n$ legs is denoted by $\nabla_n\mathscr H^1(C,\mathbb G)$. Once we fix $n$ closed points $(\pv_1,\dots, \pv_n)=\uv$ in $C$ to which we refer as characteristic places we can introduce boundedness conditions $Z_\pv$ for all $\pv\in\uv$ and $H$-level structures. Here a bound $Z_\pv$ is roughly a $L^+G_\pv$-invariant closed subscheme of the affine flag variety $\hat \fl_{\mathbb G_\pv}$ (see \ref{boundsdef} for a precise definition) and $H$ is an open, compact subgroup of $\mathbb G(\A^\uv)$, where $\A^\uv$ denotes the ring of the adeles outside $\uv$. Then we denote by $\nabla_n^{\hat Z_\uv, H}\mathscr H^1(C,\mathbb G)$ the stack which parametrize $\mathbb G$-shtukas $\underline \gtor$ over $S$ bounded by $(\hat Z_\pv)_{\pv\in\uv}$ together with an $H$-level structure. At the beginning of the third section we define all the parameters $(C,\mathbb G,\uv,Z_\uv,H)$ as a shtuka datum. This definition results in the natural question of if an appropriate change of this shtuka datum induces a morphism of the corresponding moduli spaces and if so, which properties it has.\\
In subsection 3.1 we define a morphism of shtuka data and clarify what an appropriate change of the shtuka datum should be. Roughly a morphism from $(C,\mathbb G,\uv,\hat Z_\uv, H)$ to $(C',\mathbb G',\underline w, \hat Z'_{\underline w}, H')$ is a pair $(\pi, f)$, where $\pi:C\to C'$ is a finite morphism and $f$ is a morphism of group schemes from the Weil restriction $\pi_\star \mathbb G$ to $\mathbb G'$ such that $\underline w, \hat Z'_{\underline w}$ and $H'$ satisfy certain conditions.\\
In the following subsections we then answer the questions about the functoriality of $\nabla_n^{\hat Z_\uv, H}\mathscr H^1(C,\mathbb G)$. More precisely, we first consider in subsection \ref{changeofcurve} the case that we basically only change the curve $C$, which yields the following main result of that subsection.

\theo{cf. Theorem \ref{basechangegeneral}}{
Let $(C,\mathbb G, \uv, \hat Z_{\uv}, H)$ be a shtuka datum and $\pi:C\to C'$ a finite morphism of smooth, projective, geometrically irreducible curves over $\F_q$ with $\pw_i:=\pi(\pv_i)$ and $\underline\pw:=(\pw_1,\dots,\pw_n)$. Then the morphism $(\pi,id_{\pi_\star\mathbb G}):(C,\mathbb G, \uv, \hat Z_{\uv}, H)\to (C',\pi_\star\mathbb G, \underline\pw, \pi_\star \hat Z_{\underline\pv}, \pi_\star H)$ of shtuka data (see definition \ref{defshtukadata} and remark \ref{remshtukadata}) induces a finite morphism of the moduli stacks
\begin{equation*}
\pi_\star: \nabla_n^{\hat Z_{\uv},H}\mathscr H^1(C,\mathbb G)\to \nabla_n^{\pi_\star \hat Z_{\underline\pv},\pi_\star H}\mathscr H^1(C',\pi_\star\mathbb G)
.
\end{equation*}
}

The construction of this morphism and the proof of the theorem relies on a lemma in subsection \ref{shtukdat} that states an equivalence of categories between $\mathbb G$-torsors over $C$ and $\pi_\star \mathbb G$ over $C'$.\\ 
The next subsection \ref{secchangegroup} addresses the questions about functoriality in the case that we only change the group scheme\ $\mathbb G$. While we construct  a morphism 
\begin{equation*}
\nabla_n^{\hat Z_\uv,H}\mathscr H^1(C,\mathbb G)\to \nabla_n^{\hat Z'_\uv,H'}\mathscr H^1(C,\mathbb G')
\end{equation*}
for all morphisms $(id_C,f)$ of shtuka data, we need to make different assumptions to state different results on the properties of this morphism. Assuming that $f:\mathbb G\to \mathbb G'$ is generically an isomorphism, we get the following first main result of this subsection.

\theo{cf. Theorem \ref{changeparahoric2}}{
Let $\underline w=(w_1,\dots,w_m)$ be a finite set of closed points in $C$ and let $(id_C,f):(C,\mathbb G,\uv,\hat Z_\uv,H)\to (C,\mathbb G',\uv,\hat Z'_\uv,H)$ be a morphism of shtuka data, where $f:\mathbb G\to \mathbb G'$ is an isomorphism over $C\backslash 
\underline w$. Then the morphism 
\begin{equation*}
f_\star:\ \nabla_n^{\hat Z_\uv, H}\mathscr H^1(C,\mathbb G)\to \nabla_n^{\hat Z'_\uv, H} \mathscr H^1(C,\mathbb G'), \quad (\underline\gtor,\gamma H) \mapsto (f_\star \underline\gtor,\gamma H)
\end{equation*}
is schematic and quasi-projective. In the case that $\mathbb G$ is a parahoric Bruhat-Tits group scheme this morphism is projective.
For any morphism $(\underline\gtor', \gamma' H):S\to \nabla_n^{\hat Z'_\uv, H} \mathscr H^1(C,\mathbb G')$ 
\begin{equation*}
\mbox{the fiber product\ \ }
S\times_{\nabla_n^{\hat Z'_\uv, H} \mathscr H^1(C,\mathbb G')}\nabla_n^{\hat Z'_\uv, H} \mathscr H^1(C,\mathbb G)\mbox{\quad is given by a   closed subscheme of}
\end{equation*}
\begin{equation*}
S\times_{\F_q}\left((L^+_{w_1}(\gtor')/ L^+\widetilde{\mathbb G_{w_1}})\times_{\F_q}\dots \times_{\F_q}(L^+_{w_m}(\gtor')/ L^+\widetilde{\mathbb G_{w_m}})\right).
\end{equation*}
If $\hat Z_\uv$ arises as a base change of $\hat Z'_\uv$ for all $\pv\in \uv$, the morphism $f_\star$ is surjective.
}

This result will again be important for the verification of the axioms on the moduli space $\nabla_n^{\hat Z_\uv, H}\mathscr H^1(C,\mathbb G)$ in \cite{Bre19}.
If we assume that $f:\mathbb G\to \mathbb G'$ is a closed immersion instead of a generic isomorphism we get the second main result of this subsection.

\theo{cf. Theorem \ref{theounram}}{
Let $f:\mathbb G\to \mathbb G'$ be a closed immersion of smooth, affine group schemes over $C$. Then the induced morphism $f_\star: \nabla_n\mathscr H^1(C, \mathbb G)\to \nabla_n\mathscr H^1(C,\mathbb G')$ is unramified and schematic.
}

Assuming additionally that $\mathbb G$ is a parahoric Bruhat-Tits group scheme, we prove as well:
\theo{cf. Theorem \ref{theoclosedfinite}}{
Let $\mathbb G$ be a parahoric Bruhat-Tits group scheme and $f:\mathbb G\to \mathbb G'$ be a closed immersion of smooth, affine group schemes and $\uv=(\pv_1,\dots,\pv_n)$ be a set of closed points in $C$. Then the induced morphism
\begin{equation*}
f_\star: \nabla_n\mathscr H^1(C, \mathbb G)^{\uv}\to \nabla_n\mathscr H^1(C,\mathbb G')^{\uv}\mbox{\quad is proper and in particular finite.}
\end{equation*}

}
Apart from the interest of these morphisms in the study of $\nabla_n^{\hat Z_\uv, H}\mathscr H^1(C,\mathbb G)$ in general, there are two other motivations. The first one is that the finiteness results in theorem \ref{basechangegeneral} and theorem \ref{theoclosedfinite} potentially enable us to formulate 
an analog of the Andr\'e-Oort conjecture for moduli spaces of global $\mathbb G$-shtukas, as explained in more detail in remark \ref{remAO}. The second motivation arises from the study of Newton and Kottwitz-Rapoport stratifications of $\nabla_n^{\hat Z_\uv, H}\mathscr H^1(C,\mathbb G)$ in \cite{Bre19}, where the above results are again needed.\\
We remark that the way to choose a bound $\hat Z_v$ as we do in this article is quite general. In \cite{Bre19} we will restrict our interest to moduli spaces of global $\mathbb G$-shtukas that are bounded by a cocharacter of a maximal torus of $\mathbb G$, which can also be formulated using Grassmanians and seems more natural when working with local models.\\

\textbf{Acknowledgements.} I would like to thank my advisor Urs Hartl, for all his helpful discussions. I am grateful to Timo Richarz and Johannes Ansch\"utz for their comments and remarks. During the work of this project, the author was supported by the SFB 878 "Groups, Geometry \& Actions" of the German Science Foundation (DFG), the CNRS and the ERC Advanced Grant 742608 "GeoLocLang".

\section{Preliminaries}\label{preli}
Before we start with the functoriality of $\nabla_n\mathscr H^1(C,\mathbb G)$, we introduce the basic objects and notations that we use in this article. 
Most of the notations introduced in this section, can also be found in \cite{AH13} and \cite{AH14}. 
Let $q$ be a power of some prime number $p$. 
\index{Preliminaries}
We start with a smooth, projective, geometrically irreducible curve $C$
\nomenclature{\sortkey{C}$C$}{smooth projective irreducible curve over $\F_q$}
 over the field $\F_q$ with $q$ elements.
We denote by $Q:=\F_q(C)$ 
\nomenclature{\sortkey{Q}$Q$}{function field of $C$}
its function field. For a closed point $v\in C$ we denote by $A_v$ the completion of the stalk $\oo_{C,v}$ and by $Q_v$ the fraction field of $A_v$. Furthermore we choose a uniformizer $z_v$ in $A_v$, denote the residue field of $v$ by $\F_v$ and set $deg\ v=[\F_v:\F_q]$.\\
Let \nomenclature{\sortkey{G}$\mathbb G$}{smooth, affine group scheme over $C$} 
$\mathbb G$ be a smooth, affine group scheme over $C$ and $\mathbf G:=\mathbb G\times_C Q$ its generic fiber. 
\nomenclature{\sortkey{G}$\mathbf G$}{generic fiber of $\mathbb G$}
We write $\mathbb G_{\pv}:=\mathbb G\times_CSpec\ A_\pv$ and $\mathbf G_\pv:=\mathbb G\times_CSpec\ Q_\pv=\mathbb G_{\pv}\times_{A_\pv}Spec\ Q_\pv$ for the appropriate base changes.\\
For an $\F_q$-scheme $S$ we denote by $\sigma_S:S\to S$ the absolute $\F_q$-Frobenius, which acts as the $q$-power map on the structure sheaf. Further we define $\sigma$ as the endomorphism $id_C\times\sigma_S$ of $C_S:=C\times_{\F_q}S$. For a morphism $s:S\to C$ we denote as usual by $\Gamma_s:S\to C_S$ the graph of $s$, which is a closed immersion.\\
Let $\mathcal X$ be a site with a final object $x$ and $G$ a sheaf of groups on $\mathcal X$. Then a (right) $G$-torsor is a sheaf $\mathcal G$ on $\mathcal X$ with a right action of $G$ on $\mathcal G$ such that $G\times\mathcal G\simeq \mathcal G\times \mathcal G,\ \ (g,h)\mapsto (hg,h)$ is an isomorphism and $\mathcal G(U)\neq\emptyset$ for some covering $U\to X$. When we speak about a torsor, we always mean a right torsor and if nothing else is mentioned we mean a sheaf on the big \'{e}tale site of a scheme. For any scheme $S$ we write $S_{\acute Et}$ for the big \'{e}tale site of this scheme. We denote by $\mathscr H^1(C,\mathbb G)$ the stack fibered over $(\F_q)_{\acute Et}$, whose fiber category $\mathscr H^1(C,\mathbb G)(S)$ is given by the category of $\mathbb G$-torsors over $C_S$.

\paragraph{Global $\mathbb G$-Shtukas:}\mbox{}

\smallskip
Let $S$ be an $\F_q$-scheme. A global $\mathbb G$-shtuka over $S$ is a tuple
$\underline{\mathcal G}=(\mathcal G, s_1,\dots, s_n, \tau_{\mathcal G})$, where
\begin{itemize}[label=$-$]
 \item $\mathcal G$ is a $\mathbb G$-torsor over $C_S$, 
 \item $s_1,\dots, s_n$ are morphisms $S\to C$ and
 \item $\tau_{\mathcal G}:\sigma^\star\mathcal G|_{C_S\backslash(\underset{i}{\bigcup}\Gamma_{s_i})} \to \mathcal G|_{C_S\backslash(\underset{i}{\bigcup}\Gamma_{s_i})} $ is an isomorphism of the $\mathbb G$-torsors $\sigma^\star\mathcal G$ and $\mathcal G$ restricted to $C_S\backslash(\Gamma_{s_1}\bigcup\dots\bigcup\Gamma_{s_n})$.
\end{itemize}

We take the notation $\sa$ from \cite[Definition 2.12]{AH14} for the stack fibered over $(\F_{q})_{\acute Et}$ whose $S$-valued points for a scheme $S$ are given by $\mathbb G$-shtukas $\underline{\mathcal G}$ over $S$. Morphisms from $(\mathcal G, s_1,\dots, s_n,\tau)$ to $(\mathcal G', s'_1,\dots, s'_n,\tau')$ in the fiber category $\sa(S)$ only exist if $s_i=s'_i$ and are given by morphisms $f:\mathcal G\to \mathcal G'$ of $\mathbb G$-torsors over $C_S$ such that $f\circ \tau=\tau'\circ \sigma^\star f$.
Given two $\mathbb G$-shtukas $\underline {\mathcal G}$ and $\underline {\mathcal G}'$ over $S$ with $s_i=s'_i$, we also define a quasi-isogeny from $\underline{\mathcal G}$ to $\underline{\mathcal G'}$ to be an isomorphism $f:\mathcal G|_{C_S\backslash D_S}\to \mathcal G'|_{C_S\backslash D_S}$ of $\mathbb G$-torsors satisfying $f\circ \tau=\tau'\circ \sigma^\star f$, where $D$ is some effective divisor on $C$.
The moduli space $\sa$ is an ind-algebraic stack that is ind-separated and locally of ind-finite type \cite[Theorem 3.14]{AH13}.

\paragraph{Loop Groups: }\mbox{}

\smallskip
Let $\F$ be a finite field and $\mathbb H$ be a smooth, affine group scheme over $\mathbb D:=Spec\ \F\sem{z}$, with generic fiber $\mathbf H:=\mathbb H\times_{\mathbb D}{\dot{\mathbb D}}$ where $\dot{\mathbb D}:=Spec\ \F\semr{z}$. We are mainly interested in the case that $\mathbb D\simeq Spec\ A_\pv$ and $\mathbb H=\mathbb G_\pv$ for some closed point $\pv\in C$.\\
We recall that the sheaf of groups $L^+\mathbb H$ on $\F_{\acute Et}$, whose $R$-valued points for an $\F$-algebra $R$ are given by 
\begin{equation*}
L^+\mathbb H(R):=\mathbb H(R\sem{z}):=\mathbb H(\mathbb D_{R}):=Hom_\mathbb D(\mathbb D_R,\mathbb H)\ \quad\mbox{ with }\mathbb D_R:=Spec\ R\sem{z},
\end{equation*}
is an infinite-dimensional affine group scheme over $\F$. It is called the group of positive loops associated with $\mathbb H$.
The group of loops associated with $\mathbf H$ is the sheaf of groups $L\mathbf H$ on $\F_{\acute Et}$, whose $R$-valued points are defined by 
\begin{equation*}
L\mathbf H(R):=\mathbf H(R\semr{z}):=\mathbf H(\dot{\mathbb D}_R):=Hom_{\dot{\mathbb D}}(\dot{\mathbb D}_R,\mathbf H)\ ,
\end{equation*}
where we write $R\semr{z}=R\sem{z}\left[\frac{1}{z}\right]$ and $\dot{\mathbb D}_R=Spec\ R\semr{z}$. The loop group $L\mathbf H$ is an ind-scheme of ind-finite type over $\F$.

\paragraph{Torsors for Loop Groups: }\mbox{}\label{torsorloop}

\smallskip
We write $\mathscr H^1(\F, L^+\mathbb H)$ for the stack fibered over $(\F)_{\acute Et}$ whose fiber category $\mathscr H^1(\F, L^+\mathbb H)(S)$ is the category of $L^+\mathbb H$-torsors over $S$. In the same way $\mathscr H^1(\F, L\mathbf         H)$ denotes the stack fibered over $(\F)_{\acute Et}$ whose fiber category $\mathscr H^1(\F, L\mathbf H)(S)$ is the category of $L\mathbf H$-torsors over $S$. There is a natural $1$-morphism 
\begin{equation}\label{associatltor}
L:\mathscr H^1(\F, L^+\mathbb H)\to \mathscr H^1(\F, L\mathbf H),\quad \mathcal L^+\mapsto \mathcal L
\end{equation}
induced by the inclusion of sheaves $L^+\mathbb H\subset L\mathbf H$.\\
We now consider also the $z$-adic completions of $\mathbb D$ and $\mathbb H$ and denote them by $\hat{\mathbb D}:=Spf\ \F\sem{z}$ and $\hat{\mathbb H}:=\mathbb H\times_{\mathbb D} \hat{\mathbb D}$. Later when we pass from global $\mathbb G$-shtukas to local $\mathbb G_\pv$-shtukas we often need to know that $L^+\mathbb H$-torsors are equivalent to formal $\hat{\mathbb H}$-torsors. So we recall that for an $\F$-scheme $S$ a $z$-adic formal scheme $\mathcal H$ over $\hat{\mathbb D}_S:=\hat{\mathbb D}\hat{\times}_\F S$ together with an action $\hat{\mathbb H}\hat{\times}_{\hat{\mathbb D}}\mathcal H\to \mathcal H$ of $\hat{\mathbb H}$ is called a formal $\hat{\mathbb H}$-torsor if there is an \'{e}tale covering $S'\to S$ and an $\hat{\mathbb H}$-equivariant isomorphism
$\mathcal H\hat{\times}_{\hat{\mathbb D}_S}\hat{\mathbb D}_{S'}\longrightarrow \hat{\mathbb H}\times_{\hat{\mathbb D}}\hat{\mathbb D}_{S'}$, where $\hat{\mathbb H}$ is acting on itself by right multiplication.\\
We denote by $\mathscr H^1(\hat{\mathbb D}, \hat{\mathbb H})$ the category fibered in groupoids over $(\F)_{\acute Et}$ whose fiber category $\mathscr H^1(\hat{\mathbb D}, \hat{\mathbb H})(S)$ is the groupoid of formal $\hat{\mathbb H}$-torsors over $S$. We remark that Arasteh Rad and Hartl proved in \cite[Proposition 2.4]{AH14} that there is a natural isomorphism of stacks $\mathscr H^1(\hat{\mathbb D}, \hat{\mathbb H})\longrightarrow \mathscr H^1(\F,L^+\mathbb H)$.
It sends a formal $\hat{\mathbb H}$-torsor $\mathcal H$ to the sheaf
\begin{equation*}
\mathcal L^+:S_{\acute Et}\to \mathbf{Sets},\quad T\mapsto \mathrm{Hom}_{\hat{\mathbb D}_S}(\hat{\mathbb D}_T,\mathcal H)
\end{equation*}
which becomes a $L^+\mathbb H$-torsor under the action of $L^+\mathbb H(T)=\mathrm{Hom}_{\hat{\mathbb D}}(\hat{\mathbb D}_T,\hat{\mathbb H})$.

\paragraph{Local $\mathbb H$-Shtukas: }\label{localshtukas}\mbox{}

\smallskip
Let $S$ be a $\F$-scheme and $\hat\sigma$ its absolut $\F$-Frobenius. If $\F$ equals $\F_q$ or $\F_\pv$ we will write $\sigma$ and $\sigma_\pv$ respectively, instead of $\hat \sigma$. 
\nomenclature{\sortkey{o}$\sigma,\ \sigma_\pv,\ \hat\sigma$}{absolut $\F_q$ (resp. $\F_\pv$, $\F$) Frobenius of a scheme $S$ over these fields}
With the previous notations a local $\mathbb H$-shtuka over $S$ is a pair $(\mathcal L^+,\tau_{\mathcal L})$ where
\begin{itemize}[label=$-$]
\item $\mathcal L^+$ is a $L^+\mathbb H$-torsor over $S$ and
\item $\tau_{\mathcal L}:\hat\sigma^\star\mathcal L\longrightarrow\mathcal L$ is an isomorphism of the associated loop group torsors from \eqref{associatltor} in \ref{torsorloop}.
\end{itemize}
A morphism from $\underline{\mathcal L}=(\mathcal L^+,\tau_{\mathcal L})$ to $\underline{\mathcal {L}'}=({{\mathcal L}^+}',\tau'_{\mathcal L'})$ of two local $\mathbb H$-shtukas over $S$ is a morphism $f:\mathcal L^+\to \mathcal {L^+}'$ of $L^+\mathbb H$-torsors over $S$ satisfying $\tau_{\mathcal L}\circ f=\tau_{\mathcal L'}\circ \hat\sigma^\star f$.
A quasi-isogeny from $\underline{\mathcal L}=(\mathcal L^+,\tau_{\mathcal L})$ to $\underline{\mathcal {L}'}=(\mathcal {L^+}',\tau'_{\mathcal L'})$ is an isomorphism $f:\mathcal L\to \mathcal L'$ of the associated $L\mathbf H$-torsors satisfying $\tau_{\mathcal L}\circ f=\tau_{\mathcal L'}\circ \hat\sigma^\star f$.\\
A local $\mathbb H$-shtuka $(\mathcal L^+,\tau_{\mathcal L})$ is called \'{e}tale if $\tau_{\mathcal L}:\hat\sigma^\star\mathcal L\to \mathcal L$ comes already from an isomorphism  $\tau_{\mathcal L}:\hat\sigma^\star\mathcal L^+\to \mathcal L^+$ of the $L^+\mathbb H$-torsors.
We denote the category of local $\mathbb H$-shtukas over $S$ by $Sht_{\mathbb H}(S)$ and the category of \'{e}tale local $\mathbb H$-shtukas over $S$ by $\acute{E}tSht_{\mathbb H}(S)$.\\
 We recall the Corollary \cite[Corollary 2.9]{AH14} that states that if $\mathbb H$ has a connected special fiber, then any \'{e}tale local shtuka over an separably closed field $k$ is already isomorphic to $(L^+\mathbb H_k,\ 1\cdot\hat\sigma^\star)$.\\
Let $\F\sem{\zeta}$ be the power series ring over $\F$ in a variable $\zeta$. We denote by $\mathcal Nilp_{\F\sem{\zeta}}$ the category of schemes over $Spec\ \F\sem{\zeta}$ on which $\zeta$ is locally nilpotent in the structure sheaf. Therefore $\mathcal Nilp_{\F\sem{\zeta}}$ is the full subcategory of formal schemes over $Spf\ \F\sem{\zeta}$ consisting of ordinary schemes. We will define boundedness conditions only for local shtukas over a scheme $S$ in $\mathcal Nilp_{\F\sem{\zeta}}$.\\

\paragraph{The Affine Flag Variety:}\mbox{}

\smallskip
Let $\mathbb H$ be a smooth, affine group scheme over $Spec\ \F\sem{z}$ as before. Then the affine flag variety $\mathcal F\mathit l_{\mathbb H}$ is defined as the quotient sheaf $L\mathbf H/L^+\mathbb H$ on $\F_{\acute Et}$, that is the sheaf associated to the pre-sheaf 
\begin{equation*}
T\mapsto L\mathbf H(T)/L^+\mathbb H(T).
\end{equation*}
By \cite[Theorem 1.4]{PR08} $\mathcal F\mathit l_{\mathbb H}$ is represented by an ind-scheme which is ind-quasi-projective and in particular ind-separated and of ind-finite type over $\F$. By \cite[Theorem A]{2Ric16} $\fl_{\mathbb H}$ is ind-projective if and only if $\mathbb H$ is a Bruhat-Tits group scheme over $\F\sem{z}$ in the sense of \cite[Definition 5.2.6]{BT84}. We also remark that $L^+\mathbb H$ acts from the left on $\mathcal F\mathit l_{\mathbb H}$.

\paragraph{Bounds in $\hat{\mathcal F\mathit l}_{\mathbb H}$}\label{boundsdef}\mbox{}

\smallskip
We fix an algebraic closure $\F\semr{\zeta}^{alg}$ of $\F\semr{\zeta}$. For a finite extension $R$ of discrete valuation rings $\F\sem{\zeta}\subset R\subset\F\semr{\zeta}^{alg}$ with residue field $\kappa_R$ we denote similar as before with $\mathcal Nilp_R$ the category of $R$-schemes on which $\zeta$ is locally nilpotent. Furthermore we set $\hat{\mathcal F\mathit l}_{\mathbb H,R}:=\fl_\mathbb H \hat\times_\F Spf\ R$ as well as $\hat\fl_{\mathbb H}:=\fl_{\mathbb H,\F\sem{\zeta}}$.\\
Now let $R$ and $R'$ be two such finite extensions of discrete valuation rings and let $\hat Z_R\subset\hat\fl_{\mathbb H,R}$ and $\hat{Z'}_{R'}\subset\hat\fl_{\mathbb H,R'}$
be two closed ind-subschemes. We call $\hat Z_R$ and $\hat{Z'}_{R'}$ equivalent if there is a finite extension $\widetilde R$ of discrete valuation rings as above containing $R$ and $R'$ such that $\hat Z_R\times_{Spf\ R}Spf\ \widetilde R=\hat{Z'}_{R'}\times_{Spf\ R'}Spf\ \widetilde R$ as closed ind-subschemes of $\hat\fl_{\mathbb H,\widetilde R}$.\\
Now a bound is defined (compare \cite[Definition 4.8]{AH14} and \cite[Definition 4.5]{AH13}) as an equivalence class $\hat Z=[\hat Z_R]$ of closed ind-subschemes $\hat Z_R\subset\hat\fl_{\mathbb H,R}$ satisfying
\begin{itemize}[label=$-$]
\item firstly, that all subschemes $\hat Z_R$ are stable under the left action of $L^+\mathbb H$ on $\hat\fl_{\mathbb H,R}$, and
\item secondly, that all the special fibers $Z_R:=\hat Z_R\times_{Spf\ R}Spec\ \kappa_R$ are quasi-compact and connected subschemes of $\fl_{\mathbb H}\hat\times_{\F}\kappa_R$.
\end{itemize}
We remark that in \cite{AH13} and \cite{AH14} the definition of a bound does not require the special fibers to be connected. We make this assumption because it does not change the theory and simplifies the formulation of certain statements. In fact, if $\hat Z$ is the disjoint union of two bounds $\hat Z_1\coprod \hat Z_2$ then the moduli space $\nabla_n^{\hat Z}\mathscr H^1(C,\mathbb G)$ that will be defined in paragraph \ref{BogGS} is the disjoint union of the moduli spaces $\nabla_n^{\hat Z_1}\mathscr H^1(C,\mathbb G)$ and $\nabla_n^{\hat Z_2}\mathscr H^1(C,\mathbb G)$. 

\paragraph{The Reflex Ring: }\label{parreflex}\mbox{}

\smallskip
For an equivalence class $\hat Z=[\hat Z_R]$ as above we set $G_{\hat Z}:=\{\gamma\in Aut_{\F\sem{\zeta}}(\F\semr{\zeta}^{alg})\ |\ \gamma(\hat Z)=\hat Z \}$. 
The ring $R_{\hat Z}$ is defined as the intersection of the fixed field of $G_{\hat Z}$ in $\F\semr{\zeta}^{alg}$ with all the finite extensions $R$ over which a representative $\hat Z_R$ of $\hat Z$ exists. In the case that $\hat Z$ is a bound, we call $R_{\hat Z}$ the reflex ring of $\hat Z$.\\
It is not always clear if there exists a representative of $\hat Z$ over $R_{\hat Z}$. We write $\kappa_{\hat Z}$ and $\kappa_R$ for the residue fields of $R_{\hat Z}$ and $R$ respectively. Then the special fiber $Z_R:=\hat Z_R\times_{Spf\ R}\kappa_R$ arises from a unique closed subscheme $Z\subset \fl_{\mathbb H}\times_{\F}\kappa_{\hat Z}$. This follows from Galois descent for closed ind-subschemes of $\fl_{\mathbb H}$, which is effective. The subscheme $Z$ is called the special fiber of $\hat Z$.

\paragraph{Boundedness of Local $\mathbb H$-Shtukas: }\mbox{}

\smallskip
Let $\hat Z$ be a bound with reflex ring $R_{\hat Z}$. Furthermore let $\mathcal L^+$ and $\mathcal {L^+}'$ be two $L^+\mathbb H$-torsors over a scheme $S$ in $\mathcal Nilp_{R_{\hat Z}}$ and $\delta:\mathcal L\to \mathcal L'$ an isomorphism of the associated $L\mathbf H$-torsors. We choose a covering $S'\to S$ such that there are trivializations $\alpha:\mathcal L^+\isom L^+\mathbb H_{S'}$ and $\alpha':\mathcal {L^+}'\isom L^+\mathbb H_{S'}$. Then the automorphism $\alpha'\circ\delta\circ\alpha^{-1}:L\mathbf H_{S'}\isom L\mathbf H_{S'}$ defines a morphism $S'\to L\mathbf H$. 

For any finite extension $R$ of $R_{\hat Z}$ we have an induced morphism 
\begin{equation}\label{critbound}
S'\times_{Spf\ R_{\hat Z}}Spf\ R\longrightarrow
L\mathbf H\hat\times_{\F}Spf\ R\longrightarrow
\hat\fl_{\mathbb H,R}.
\end{equation}
Then $\delta$ is said to be bounded by $\hat Z$ if for all trivializations $\alpha$ and $\alpha'$ and all finite extensions $R$ of $R_{\hat Z}$ with a representative $\hat Z_R$, this morphism \eqref{critbound} factors through $\hat Z_R$.
By \cite[Remark 4.9]{AH14} $\delta$ is bounded if and only if this condition is satisfied for one trivialization and for one such extension $R$.
By definition a local $\mathbb H$-shtuka $(\mathcal L^+,\tau_{\mathcal L})$ is bounded by $\hat Z$ if $\tau_{\mathcal L}$ is bounded by $\hat Z$.

\paragraph{A Version of the Theorem of Beauville-Laszlo: }\label{beauvillelaszlo}\mbox{}

\smallskip
Let $\pv\in C$ be a closed point and set $C^{\pv}:=C\backslash\{\pv\}$ 
 as well as $C^{\pv}_S:=C^{\pv}\times_{\F_q}S$. We define $\mathscr H_e^1 (C^\pv, \mathbb G)$ as the category fibered in groupoids over $(\F_q)_{\acute Et}$ whose fiber category $\mathscr H_e^1(C^\pv, \mathbb G)(S)$ consists of those $\mathbb G$-torsors over $C^\pv_S$ that can be extended to a $\mathbb G$-torsor over $C_S$. By restricting a $\mathbb G$-torsor $\gtor$ over $C_S$ to $C_S^\pv$ we get a morphism $\mathscr H^1(C,\mathbb G)\to \mathscr H_e^1(C^\pv,\mathbb G)$.
 We further introduce the notation $\widetilde{\mathbb G_{\pv}}=Res_{A_\pv/ \F_q\sem{z_\pv}}\mathbb G_\pv$ and $\widetilde{\mathbf G_\pv}:=\widetilde{\mathbb G_\pv}\times_{\F_q\sem{z_\pv}}\F_q\semr{z_\pv}$. For $\gtor\in \mathscr H^1(C,\mathbb G)$ the base change $\gtor_\pv:=\gtor\times_{C_S}(Spf\ A_\pv\times_{\F_q}S)$ defines a formal $\mathbb G_\pv$-torsor over $Spf\ A_\pv\hat\times_{\F_q}S$. Its Weil restriction $Res_{A_{\pv}/ \F_q\sem{z_\pv}}\gtor_{\pv}$ defines a formal $\widetilde{\mathbb G_\pv}$-torsor over $Spf\ \F_q\sem{z_\pv}\hat\times_{\F_q}S$. Using the category equivalence in \ref{torsorloop} it corresponds to an object in $\mathscr H^1(\F_q,L^+\widetilde{\mathbb G_\pv})(S)$ that we denote by $L^+_\pv(\mathcal G)$ which defines a functor 
 \begin{equation*}
 L^+_\pv:\mathscr H^1(C,\mathbb G)\longrightarrow \mathscr H^1(\F_q, L^+\widetilde{\mathbb G_\pv}),\quad \mathcal G\mapsto L^+_\pv(\mathcal G).
 \end{equation*}
 
 Furthermore we have the functor
 \begin{equation*}
 L_\pv: \mathscr H^1_e(C^\pv,\mathbb G)\to \mathscr H^1(\F_q,L\widetilde{\mathbf G_\pv}),\quad \gtor\big|_{C^\pv_S}\mapsto L(L^+_\pv(\gtor))=L_\pv(\gtor)
 \end{equation*}
which is independent of the extension $\gtor$ of $\gtor\big|_{C^\pv_S}$. Now a version of the theorem of Beauville-Laszlo, that is proven in \cite[Lemma 5.1]{AH14}, states that the following diagram is cartesian.
\begin{equation*}
\xymatrix{
{\mathscr H^1(C,\mathbb G)} \ar[r] \ar_{L^+_\pv}[d]& {\mathscr H^1_e(C^\pv,\mathbb G)} \ar[d]_{L_\pv} \\
\mathscr H^1(\F_q, L^+\widetilde{\mathbb G}_{\pv}) \ar[r] & \mathscr H^1(\F_q, L\widetilde{\mathbf G}_{\pv})
}
\end{equation*}

\paragraph{The Global-Local Functor: }\mbox{}\label{parglobloc}

\smallskip
Now we fix $n$ closed points $\uv=\{\pv_1,\dots,\pv_n\}$ of $C$. Then define 
$A_\uv$ as the completion of the local ring $\mathcal O_{C^n, \uv}$ at the closed point $\uv$. We set
\begin{equation*}
\nabla_n\mathscr H^1(C,\mathbb G)^\uv := \nabla_n\mathscr H^1(C,\mathbb G)\times_{C^n}Spf\ A_\uv.
\end{equation*}

An $S$-valued point of $\nabla_n\mathscr H^1(C,\mathbb G)^\uv $ is therefore given by a global $\mathbb G$-shtuka $\underline\gtor=(\gtor,s_1,\dots, s_n,\tau_\gtor)$ such that $s_i:S\to C$ factors through $Spf\ A_{\pv_i}$. We now want to associate with such a global $\mathbb G$-shtuka a local $\mathbb G_{\pv_i}$-shtuka for all $\pv_i\in \uv$. 
We write $\mathbb D_{\pv_i}:=Spec\ A_{\pv_i}$ and $\hat{\mathbb D}_{\pv_i}:=Spf\ A_{{\pv_i}}$ as well as $\hat{\mathbb D}_{{\pv_i},S}:={\mathbb D}_{\pv_i}{\hat\times}_{\F_{\pv_i}}S$. Then we have:
\begin{equation*}
\hat {\mathbb D}_{{\pv_i}}{\hat\times}_{\F_q}S=\coprod_{l\in \Z/deg\ {\pv_i}} V(\mathfrak a_{{\pv_i},l})=\coprod_{l\in\Z/deg\ {\pv_i}}\hat{\mathbb D}_{{\pv_i},S}
\end{equation*}
 where $\mathfrak a_{\pv_i,l}:=\langle a\otimes 1-1\otimes s^\star(a)^{q^l}\ |\ a\in \F_{\pv_i} \rangle$ and $V(\mathfrak a_{{\pv_i},l})$ is the closed subscheme given by this ideal. We remark that $\sigma$ cyclically permutes these components and that the $\F_{\pv_i}$-Frobenius $\sigma^{deg\ {\pv_i}}$ leaves all these components $V(\mathfrak a_{{\pv_i},l})$ stable.
 For $\underline\gtor\in \nabla_n\mathscr H^1(C,\mathbb G)^{\uv}$ the base change 
 \begin{equation*}
 \gtor_{\pv_i}:=\gtor\hat\times_{C_S}(Spf\ A_{\pv_i}\hat\times_{\F_q} S)=\coprod_{l\in\Z/ deg\ \pv_i}\gtor\hat\times_{C_S}V(\mathfrak a_{v_i,l})
 \end{equation*}
defines a formal $\hat{\mathbb G}_{\pv_i}$-torsor over $\ds\coprod_{l\in\Z/deg\ \pv_i}\hat{\mathbb D}_{\pv_i,S}$ which is an object in $\ds\mathscr H^1(\hat{\mathbb D}_{\pv_i},\hat{\mathbb G}_{\pv_i})(\coprod_{l\in \Z/deg\ {\pv_i}} S)$.
Each component $\gtor\times_{C_S}V(\mathfrak a_{\pv_i,l})$ defines a formal $\hat{\mathbb G}_{\pv_i}$-torsor. Similar to the notation in  \ref{torsorloop} we denote by $\mathcal L_{i,0}^+$ the $L^+\mathbb G_{\pv_i}$-torsor associated by \cite[Proposition 2.4]{AH14} with the formal $\hat{\mathbb G}_{\pv_i}$-torsor $\gtor\times_{C_S}V(\mathfrak a_{\pv_i,0})$. Then $(\mathcal L^+_{i,0},\tau^{deg\ \pv_i})$ is a local $\mathbb G_{\pv_i}$-shtuka, where 
$\tau^{deg\ \pv_i}:(\sigma^{deg\ \pv_i})^\star\mathcal L_{i,0}\isom \mathcal L_{i,0}$ is the isomorphism of $L\mathbf G_{\pv_i}$-torsors induced by $\tau_{\mathcal G}$
 (compare also \cite[Lemma 5.1]{AH14}). More precisely, it can be written as $\tau^{deg\ \pv_i}=\tau\circ\sigma^\star\tau\circ \dots\circ(\sigma^{deg\ \pv_i-2})^\star\tau\circ(\sigma^{deg\ \pv_i-1})^\star\tau$.\\
 This now defines the global-local functor:
 \begin{align*}
 \Gamma_{\pv_i}:\nabla_n\mathscr H^1(C,\mathbb G)^{\uv}(S)&\longrightarrow Sht_{\mathbb G_{\pv_i}}
 (S)\\
 \underline \gtor=(\gtor, s_1,\dots,s_n,\tau) & 
 \longmapsto (\mathcal L^+_{i,0}\ ,\ \tau^{deg\ \pv_i})=\Gamma_{\pv_i}(\underline \gtor).
 \end{align*}
We remark that this functor transforms by \cite[Definition 5.4]{AH14} quasi-isogenies into quasi-isogenies.
If $v\notin\underline{v}$ the component $V(\mathfrak{a}_{v,0})$ exists only if $S$ is an $\mathbb{F}_v$-scheme.

If we do not restrict $\mathcal G\times_{C_S}(Spf\ A_{\pv_i}\times_{\F_q}S)$ to the component $V(\mathfrak a_{v_i,0})$ but consider its Weil restriction $L^+_{\pv_i}(\gtor)$ we get in a similar way a local $\widetilde{\mathbb G_{\pv_i}}$-shtuka $(L_{\pv_i}^+(\mathcal G),\tau_{\pv_i})$ where $\tau_{\pv_i}:=L_{\pv_i}(\tau_\gtor):\sigma^\star L_{\pv_i}(\gtor)\to L_{\pv_i}(\gtor)$.
We denote this local $\widetilde{\mathbb G_{\pv_i}}$-shtuka by $L^+_{\pv_i}(\underline \gtor)=(L_{\pv_i}^+(\mathcal G),\tau_{\pv_i})$. We remark that $L^+_{\pv}(\gtor)$ does not only exist for $\pv\in\uv$ but also for other places $\pv \in C$. In the case that $\pv\notin\uv$,   the local shtuka $L^+_\pv(\underline\gtor)$ is \'{e}tale.

\paragraph{Boundedness of Global $\mathbb G$-Shtukas: }\mbox{}\label{BogGS}

\smallskip
Recall that we fixed $n$ closed points $\uv=(\pv_1,\dots,\pv_n)$ in $C$. If the group scheme $\mathbb G$ is fixed we write for each of these points $\fl_{\pv_i}:=\fl_{\mathbb G_{\pv_i}}$ for the corresponding affine flag variety over $\F_{\pv_i}$ and $\hat\fl_{\pv_i,R}=\hat\fl_{\mathbb G_{\pv_i},R}=\fl_{\pv_i}\hat\times_{\F_{\pv_i}}Spf\ R$ for a finite extension $A_{\pv_i}\subset R$.
In each of these affine flag varietes $\hat\fl_{\pv_i}=\hat\fl_{\pv_i,A_{\pv_i}}$ we choose a bound $\hat Z_{\pv_i}=[\hat Z_{\pv_i,R}]$ with reflex ring $R_{\hat Z_{\pv_i}}$ and we write $\hat Z_{\uv}$ for the tuple $(\hat Z_{\pv_1},\dots, \hat Z_{\pv_n})$. Choosing a uniformizer $\pi_{\pv_i}$ in $R_{\hat Z_{\pv_i}}$ and defining $\F_R$ as the compositum of all the residue fields $R_{\hat Z_{\pv_i}}/(\pi_{\pv_i})$, we set $R_{\hat Z_\uv}:=\F_R\sem{\pi_{\pv_1}\dots, \pi_{\pv_n}}$.
In particular the morphism $Spf\ R_{\hat Z_\uv}\to C^n$ factors through $Spf\ A_{\uv}$. This means that every point $\underline\gtor=(\gtor, s_1,\dots,s_n, \tau_\gtor)$ in $\nabla_n\mathscr H^1(C,\mathbb G)\times_{C^n}Spf\ R_{\hat Z_\uv}(S)$ also defines an $S$-valued point in $\nabla _n\mathscr H^1(C,\mathbb G)^{\uv}$ so that we write $\Gamma_{\pv_i}(\underline\gtor)$ for its associated local $\mathbb G_{\pv_i}$-shtuka over $S$. The fact that $S\in Nilp_{R_{\hat Z_{\pv_i}}}$ allows us to ask if $\Gamma_{\pv_i}(\underline\gtor)$ is bounded by $\hat Z_{\pv_i}$.\\
We define $\nabla_n^{\hat Z_{\uv}}\mathscr H^1(C,\mathbb G)$ to be the stack consisting of these bounded global $\mathbb G$-shtukas. That means the fiber category $\nabla_n^{\hat Z_\uv}\mathscr H^1(C,\mathbb G)(S)$ is the full subcategory of $\nabla_n\mathscr H^1(C,\mathbb G)\times_{C^n}Spf\ R_{\hat Z_\uv}(S)$ that consists of those global $\mathbb G$-shtukas $\underline \gtor$ over $S$ that are bounded by $\hat Z_\uv$. By \cite[Remark 7.2]{AH13} the moduli space $\nabla_n^{\hat Z_{\uv}}\mathscr H^1(C,\mathbb G)$ is a closed ind-substack of $\nabla_n\mathscr H^1(C,\mathbb G)\times_{C^n}Spf\ R_{\hat Z_\uv}$. Moreover we 
denote by $\nabla_n^{\hat Z_\uv}\mathcal H^1(C,\mathbb G)_{{\F_R}}:=\nabla_n^{\hat Z_\uv}\mathcal H^1(C,\mathbb G){\times_{Spf R_{\hat Z_\uv}}\F_R}$  the special fiber of $\nabla_n^{\hat Z_\uv}\mathcal H^1(C,\mathbb G)$.

\paragraph{D-Level Structures: }\label{dlevel}\mbox{}

\smallskip
Let $D$ be a proper closed subscheme of $C$ and let $D_S:=D\times_{\F_q}S$ for some $\F_q$-scheme $S$ and $\gtor$ a $\mathbb G$-torsor on $C_S$. By \cite[Definition 3.1]{AH13} a $D$-level structure on $\gtor$ is a trivialization $\Psi:\gtor\times_{C_S}D_S\to \mathbb G\times_C D_S$ and $\mathscr H^1_D(C,\mathbb G)$ denotes the stack fibered over $(\F_q)_{\acute Et}$ whose fiber category $\mathscr H^1_D(C,\mathbb G)(S)$ consists of pairs $(\gtor, \Psi)$ where $\gtor\in \mathscr H^1(C,\mathbb G)(S)$ and $\Psi$ is a $D$-level structure. A morphism from $(\gtor,\Psi)$ to $(\gtor',\Psi')$ in this fiber category is given by an isomorphism $f:\gtor\to \gtor'$ of $\mathbb G$-torsors such that $\Psi=\Psi'\circ(f\times id_{D_S})$. The moduli stack of global $\mathbb G$-shtukas with $D$-level structure is denoted by $\nabla_n\mathscr H^1_D(C,\mathbb G)$. Its fiber category over $S$ is given by tuples $(\underline \gtor, \Psi)=(\gtor,s_1,\dots, s_n,\tau_\gtor,\Psi)$ where $\underline\gtor\in \nabla_n\mathscr H^1(C,\mathbb G)\times_{C^n}(C\backslash D)^n(S)$ (i.e. $s_i:S\to C$ factors through $C\backslash D$) and $\Psi$ is a $S$-level structure on $\gtor$ satisfying $\Psi\circ (\tau\times id_{D_S})=\sigma^\star(\Psi)$. A morphism from $(\underline\gtor,\Psi)$ to $(\underline \gtor', \Psi')$ in this fiber category is a morphism $f\in \nabla_n\mathscr H^1(C,\mathbb G)(S)$ (in particular an isomorphism $f:\gtor\to \gtor'$ of $\mathbb G$-torsors) satisfying $\Psi=\Psi'\circ(f\times id_{D_S})$.\\
If $D=\emptyset$ we have $\nabla_n\mathscr H_D^1(C,\mathbb G)=\nabla_n\mathscr H^1(C,\mathbb G)$. If $\pv_1,\dots,\pv_n\notin D$ and $\hat Z_\uv$ is a bound as before, we use the intuitive notations $\nabla_n\mathscr H_D^1(C,\mathbb G)^\uv$ for the base change $\nabla_n\mathscr H_D^1(C,\mathbb G)\times_{C^n}Spf\ A_\uv$ and $\nabla_n^{\hat Z_{\uv}}\mathscr H_D^1(C,\mathbb G)$ for the stack of $\mathbb G$-shtukas $\underline \gtor$ in $\nabla_n^{\hat Z_\uv}\mathscr H^1(C,\mathbb G)$ with a $D$-level structure.

\paragraph{Local Shtukas and Local $\mathrm{GL}_r$-Shtukas: }\mbox{}

\smallskip
The category of local $\mathrm {GL}_r$-shtukas over an $F_q$-scheme $S$  \ 
can be defined more explicitly. We briefly describe this here since it is useful for the definition of the Tate functors.\\
We denote by $\oo_S\sem{z}$ the sheaf of $\oo_S$-algebras on $S_{\acute Et}$ which associates with every $S$-scheme $Y$ the ring $\oo_S\sem{z}(Y):=\Gamma(Y,\oo_Y)\sem{z}$. Now every sheaf of $\oo_S\sem{z}$-modules that is fqqc-locally free of rank $r$ is by \cite[Prop 2.3]{HV11} already Zariski locally free of rank $r$. We call these locally free sheafes of $\oo_S\sem{z}$-modules of rank $r$. For a commutative ring $R$ we set $R\semr{z}:=R\sem{z}\left[\frac{1}{z}\right]$. This leads to the intuitive notation $\oo_S\semr{z}$ for the sheaf on $S_{\acute Et}$ associated to the pre-sheaf $Y\mapsto \Gamma(Y,\oo_Y)\semr{z}$. The absolut $\F$ Frobenius was denoted $\hat\sigma$ and we use the same notation for the endomorphism of $\oo_S\sem{z}$ and $\oo_S\semr{z}$ that acts as $\hat\sigma$ on sections of $\oo_S$ and as the identity on $z$. For a sheaf $M$ of $\oo_S\sem{z}$-modules we can consider the pullback $\hat\sigma^\star M:=M\otimes_{\oo_S\sem{z},\hat\sigma}\oo_S\sem{z}$. Now by \cite[Definition 4.1]{HV11} a local shtuka of rank $r$ over $S$ is a pair $(M, \tau_M)$ consisting of a locally free sheaf $M$ of $\oo_S\sem{z}$-modules of rank $r$ and an isomorphism 
$$
\tau_M:\hat\sigma^\star M \otimes_{\oo_S\sem{z}}\oo_S\semr{z}\to M\otimes_{\oo_S\sem{z}}\oo_S\semr{z}. 
$$
The local shtuka $(M,\tau_M)$ is called \'{e}tale if $\tau_M$ arises from an isomorphism $\hat\sigma^\star M\to M$ of $\oo_S\sem{z}$-modules. A morphism from $(M,\tau_M)$ to $(M',\tau_{M'})$ between two local shtukas over $S$ is a morphism $f:M\to M'$ of $\oo_S\sem{z}$-modules satisfying $\tau_{M'}\circ\hat\sigma^\star f=f\circ \tau_M$. 
A quasi-isogeny from $(M,\tau_M)$ to $(M',\tau_{M'})$ is a morphism $f:M\otimes_{\oo_S\sem{z}}\oo_S\semr{z}\to M'\otimes_{\oo_S\sem{z}}\oo_S\semr{z}$ of $\oo_S\semr{z}$-modules satisfying $\tau_{M'}\circ\hat\sigma^\star f=f\circ \tau_M$.
We denote the category of local shtukas over $S$ by $Sht_{\F_q}(S)$ and the category of \'{e}tale local shtukas over $S$ by $\acute 
EtSht_{\F_q}(S)$.
\\
Now there is a category equivalence between local $\mathrm {GL}_r$-shtukas as defined in \ref{localshtukas} and the category of local shtukas of rank $r$ over $S$ with isomorphisms as the only morphisms. It is naturally induced by the category equivalence \cite[Lemma 4.2]{HV11} of $\mathscr H^1(\F,L^+\mathrm{GL}_r)(S)$ and the category of locally free sheaves of $\oo_S\sem{z}$-modules of rank $r$ with isomorphisms as morphisms.\\

\paragraph{Tate Functors on Local $\mathbb H$-Shtukas: }\mbox{}

\smallskip
Now let 
$S$ be a connected $\F_q$-scheme with geometric base point $\overline s\in S$ and algebraic fundamental group $\pi_1(S,\overline s)$. We denote by $\mathfrak {FMod}_{\F\sem{z}[\pi_1(S,\overline s)]}$ (resp. $\mathfrak {FMod}_{\F\semr{z}[\pi_1(S,\overline s)]}$ ) the category of finite and free $\F_q\sem{z}$-modules (resp. $\F\semr{z}$ vector spaces) equipped with a continuous action of $\pi_1(S,\overline s)$. Then the dual Tate functor $\check T$ on \'{e}tale local shtukas is defined as 
\begin{equation*}
\check T_-:\acute EtSht_{\F_q}(S)\to \mathfrak {FMod}_{\F\sem{z}[\pi_1(S,\overline s)]}\qquad \underline M:=(M,\tau_M)\mapsto \check T_{\underline M}:=(M\otimes_{\oo_S\sem{z}}\kappa(\overline s)\sem{z})^{\tau_M},
\end{equation*}
where the superscript $\tau_M$ denotes the $\tau_M$ invariants. 
The rational dual Tate functor  is defined by 
\begin{equation*}
\check V_-:\acute EtSht_{\F_q}(S)\to \mathfrak {FMod}_{\F\semr{z}[\pi_1(S,\overline s)]}\qquad \underline M:=(M,\tau_M)\mapsto \check V_{\underline M}:=\check T_{\underline M}\otimes_{\F\sem{z}}\F\semr{z}.
\end{equation*}

We also need Tate functors for local $\mathbb H$-shtukas. To define these we denote by $Rep_{\F\sem{z}}\mathbb H$ the category of representations $\rho:\mathbb H\to \mathrm{GL}(V)$, where $V$ is a finite-free $\F\sem{z}$-module and $\rho$ a morphism of algebraic groups over $\F\sem{z}$. 
Any such $\rho$ naturally induces, as described in \cite[section 3, above Definition 3.5]{AH14}, a functor $\rho_\star:\acute EtSht_{\mathbb H}(S)\to \acute EtSht_{\F}(S)$ that is compatible with quasi-isogenies.\\
Let $Funct^\otimes(Rep_{\F\sem{z}}\mathbb H,\ \mathfrak{FMod}_{\F\sem{z}[\pi_1(S,\bar s)]} )$ and $Funct^{\otimes}(Rep_{\F\semr{z}}\mathbb H,\ \mathfrak{FMod}_{\F\semr{z}[\pi_1(S,\bar s)]})$ be the categories of the appropriate tensor functors whose morphisms are isomorphisms of functors. Now the dual Tate functor $\check{\mathcal T}_-$ and the rational dual Tate functor $\check{\mathcal V}_-$ are defined by 
\begin{align*}
\check{\mathcal T}_-:\ & \acute EtSht_\mathbb H(S)\longrightarrow Funct^\otimes(Rep_{\F\sem{z}}\mathbb H,\ \mathfrak{FMod}_{\F\sem{z}[\pi_1(S,\bar s)]} )& \underline{\mathcal L}\mapsto (\check{\mathcal T}_{\underline{\mathcal L}}:\rho\mapsto \check T_{\rho_\star\underline{\mathcal L}})
\\
\check{\mathcal V}_-:\ &\acute EtSht_\mathbb H(S)\longrightarrow Funct^{\otimes}(Rep_{\F\semr{z}}\mathbb H,\ \mathfrak{FMod}_{\F\semr{z}[\pi_1(S,\bar s)]})&
\quad \underline{\mathcal L}\mapsto (\check{\mathcal V}_{\underline{\mathcal L}}:\rho\mapsto \check V_{\rho_\star\underline{\mathcal L}}).
\end{align*}

\paragraph{Tate Functors on Global $\mathbb G$-Shtukas: }\mbox{}

\smallskip
Now we assume that the tuple $\underline v=(\pv_1,\dots, \pv_n)$ is given by $n$ pairwise different places on $C$ and set $\widetilde C=C\backslash\{\pv_1,\dots,\pv_n\}$. We denote by $\OO^\uv:=\prod_{\pv\in (C\backslash \underline \pv)} A_\pv$ the integral adeles of $C$ outside $\uv$ and by $\mathbb A^\uv:=\OO^{\uv}\otimes_{\oo_{\widetilde C}}Q=\prod'_{\pv\in (C\backslash \uv)}Q_\pv$ the adeles of $C$ outside $\uv$.
Let $Rep_{\OO^\uv}\mathbb G$ be the category of representations $\rho:\mathbb G\times_CSpec\ \OO^\uv\to \mathrm{GL}_{\OO^\uv}(V)$
where $V$ is a finite free $\OO^\uv$-module and $\rho$ a morphism of group schemes over $\OO^\uv$.
Let $S$ be a connected scheme over $Spf\ A_\uv$ with a fixed geometric base point $\bar s$. We denote by $\mathfrak {Mod}_{\OO^\uv[\pi_1(S,\bar s)]}$ (resp. $\mathfrak {Mod}_{\mathbb A^\uv[\pi_1(S,\bar s)]}$) the category of $\OO^\uv$-modules (resp. $\mathbb A^\uv$-modules) with a continuous $\pi_1(S,\bar s)$ action. For a finite subscheme $D\subset C$ we set $D_{\bar s}=D\times_{\F_q}{\bar s}$ as well as $\gtor|_{D_{\bar s}}=\gtor\times_{C_S}D_{\bar s}$. Then the dual Tate functor $\check{\mathcal T}_-$ and the rational dual Tate functor $\check{\mathcal V}_-$ on global $\mathbb G$-shtukas are defined by
\begin{align*}
\check{\mathcal T}_-:\  \nabla_n\mathscr H^1(C,\mathbb G)^\uv(S)
&\longrightarrow Funct^\otimes(Rep_{\OO^\uv}\mathbb G,\ \mathfrak{Mod}_{\OO^\uv[\pi_1(S,\bar s)]} ) \\ 
\underline{\gtor}&\mapsto \left(\check{\mathcal T}_{\underline{\gtor}}:\rho\mapsto \varprojlim_{D\subset \widetilde C} (\rho_\star\underline{\gtor}|_{D_{\bar s}})^{\tau_\gtor}\right)
\\
\check{\mathcal V}_-:\  \nabla_n\mathscr H^1(C,\mathbb G)^\uv(S)
&\longrightarrow Funct^\otimes(Rep_{\A^\uv}\mathbb G,\ \mathfrak{Mod}_{\A^\uv[\pi_1(S,\bar s)]} ) \\
\underline{\gtor}&\mapsto \left(\check{\mathcal V}_{\underline{\gtor}}:\rho\mapsto \varprojlim_{D\subset \widetilde C} (\rho_\star\underline{\gtor}|_{D_{\bar s}})^{\tau_\gtor}\otimes_{\OO^\uv}\mathbb A^\uv\right).
\end{align*}

We remark that the functor $\check{\mathcal V}$ transforms by 
\cite[section 6]{AH13} quasi-isogenies into isomorphisms. Besides it is useful to know that there is a natural isomorphism $\varprojlim_{D\subset C}(\rho_\star\gtor|_{D_{\bar s}})^{\tau_\gtor}\simeq\prod_{\pv\in C\backslash \uv}\check{\mathcal T}_{L^+_\pv(\underline\gtor)}(\rho_\pv)$ writing $\rho=(\rho_\pv)_{\pv\in \widetilde C}$ with $\rho_\pv:=\rho\times id_{A_\pv}$. Here $L^+_\pv(\underline \gtor)$ is the \'{e}tale local $\widetilde{\mathbb G}_\pv$-shtuka and $\check{\mathcal T}_{L^+_\pv(\underline\gtor)}(\rho_\pv):=\check{\mathcal T}_{L^+_\pv(\underline\gtor)}(\widetilde{\rho_\pv})$ where $\widetilde {\rho_\pv}$ is the representation of $\widetilde {\mathbb G}_\pv$ induced from $\rho_\pv$ by Weil restriction (see \cite[remark 5.6]{AH14}).

\paragraph{H-Level Structures: }\label{hlevel}\mbox{}

\smallskip
Let $H$ be an open, compact subgroup of $\mathbb G(\mathbb A^\uv)$. In this paragraph we define $H$-level structures which are a generalization of the previous $D$-level structures. We denote by 
\begin{equation*}
\omega_{\OO^\uv}^\circ: Rep_{\OO^\uv}\mathbb G\longrightarrow \mathfrak {Mod}_{\OO^\uv}\ ,\qquad 
\omega_{\mathbb A^\uv}^\circ: Rep_{\mathbb A^\uv}\mathbb G\longrightarrow \mathfrak {Mod}_{\mathbb A^\uv}
\end{equation*}
the forgetful functors and by $Isom^\otimes(\omega_{\OO^\uv}^\circ,\check {\mathcal T}_{\underline\gtor})$ and $Isom^\otimes(\omega_{\mathbb A^\uv}^\circ,\check {\mathcal V}_{\underline\gtor})$ the sets of isomorphisms of tensor functors 
which are defined for every global  $\mathbb G$-shtuka $\underline\gtor$ over $S$, where $S$ is as before a scheme over $Spf\ A_\uv$ with geometric base point $\bar s\in S$. By the definition of the Tate functor $\pi_1(S,\bar s)$ acts on $\check{\mathcal T}_{\underline \gtor}$ and 
$\mathbb G(\OO^\uv)$ (resp. $\mathbb G(\mathbb A^\uv)$) acts on $\omega^\circ_{\OO^\uv}$ (resp. $\omega^\circ_{\mathbb A^\uv}$)
since we have $\mathbb G(\OO^\uv)=Aut^{\otimes}(\omega^\circ_{\OO^\uv})$ by the generalized tannakian formalism \cite[corollary 5.20]{Wed04}.
This induces an action of $\mathbb G(\OO^\uv)\times\pi_1(S,{\bar s})$ on $Isom^\otimes(\omega_{\OO^\uv}^\circ,\check{\mathcal T}_{\underline \gtor})$ and of $\mathbb G(\mathbb A^\uv)\times\pi_1(S,{\bar s})$ on $Isom^\otimes(\omega_{\mathbb A^\uv}^\circ,\check{\mathcal V}_{\underline \gtor})$.
Now by \cite[Definition 6.3]{AH14} a rational $H$-level structure $\bar\gamma$ on a global $\mathbb G$-shtuka $\underline \gtor$ in $\nabla_n\mathscr H^1(C,\mathbb G)^{\uv}(S)$ is defined as a $\pi_1(S,\bar s)$-invariant $H$-orbit $\bar\gamma=\gamma H$ in $Isom^\otimes( \omega_{\A^\uv}^\circ,\check{\mathcal V}_{\underline\gtor} )$. We denote by $\nabla_n^H\mathscr H^1(C,\mathbb G)^\uv$ the category fibered in groupoids over $(\F_q)_{\acute Et}$ with the following fiber categories.
An object in $\nabla_n^H\mathscr H^1(C,\mathbb G)^\uv(S)$ is a tuple $(\underline \gtor,\gamma)$, where $\underline \gtor\in \nabla_n\mathscr H^1(C,\mathbb G)^\uv(S)$ and $\bar\gamma$ is a $H$-level structure on $\underline\gtor$. A morphism from $(\underline \gtor,\bar\gamma)$ to $(\underline\gtor',\bar\gamma')$ over $S$ is a quasi-isogeny $f:\underline\gtor\to\underline\gtor'$ that is an isomorphism at the characteristic places $\pv_i$ and that satisfies $\check{\mathcal V}_{f}\circ \gamma H=\gamma'H$. \\(So $f:\underline\gtor|_{C_S\backslash T_S}\isom \underline\gtor'|_{C_S\backslash T_S}$ for a finite subscheme $T\subset C$ with $\pv_1,\dots,\pv_n\notin T$.)\\
Now let $D$ be  a finite subscheme $D\subset C$ with $\pv_1,\dots,\pv_n\notin D$ and let $H_D$ be the open and compact subgroup $ker(\mathbb G(\OO^\uv)\to \mathbb G(\oo_D))$ of $\mathbb G(\A^\uv)$. Then we remark hat by \cite[Theorem 6.4]{AH13} there is a canonical isomorphism of stacks 
\begin{equation}\label{leveliso}
\nabla_n\mathscr H_D^1(C,\mathbb G)^\uv\isom \nabla^{H_D}_n\mathscr H^1(C,\mathbb G)^\uv
\end{equation}

Furthermore we note that for the conjugated group $gHg^{-1}$ with $g\in \mathbf G(\A^\uv)$ there is by \cite[Remark 6.6]{AH13} a natural isomorphism $\nabla_n^{H}\mathscr H^1(C,\mathbb G)^\uv\isom \nabla_n^{gHg^{-1}}\mathscr H^1(C,\mathbb G)^\uv$ sending $(\underline\gtor,\gamma H)$ to $(\underline\gtor, \gamma  g^{-1}(gHg^{-1}))$. \\
In addition we remark that by \cite[section 6]{AH13} for a  open, compact subgroup $\widetilde H\subset \mathbb G(\A^\uv)$ contained in $H$ we have a natural finite \'etale morphism
\begin{equation}\label{levelproj}
\nabla_n^{\widetilde H}\mathscr H^1(C,\mathbb G)^\uv\isom \nabla_n^{H}\mathscr H^1(C,\mathbb G)^\uv,\ \ (\underline\gtor,\gamma \widetilde H)\mapsto (\underline\gtor, \gamma  H)
\end{equation}

If we have additionally given a bound $\hat Z_{\pv}$ at all places $v\in \uv$ we denote by $\nabla_n^{\hat Z_{\uv},H}\mathscr H^1(C,\mathbb G)$ the closed substack of $\nabla^H_n\mathscr H^1(C,\mathbb G)^\uv\times_{Spf\ A_\uv}Spf\ R_{\hat Z_\uv}$ that consists of those points $(\underline \gtor,\bar\gamma)$ such that $\underline \gtor$ is bounded by $\hat Z_\uv$.

\paragraph{Parahoric Bruhat-Tits Group Schemes}\mbox{}\label{BTgs}

\smallskip
At some points, mainly in theorem \ref{theoclosedfinite}, we will assume $\mathbb G$ to be a parahoric Bruhat-Tits group scheme, where we call a smooth, affine group scheme $\mathbb G$ over $C$ a parahoric Bruhat-Tits group scheme, if
\begin{itemize}[label=$-$]
\item all fibers are connected,
\item the generic fiber $\mathbf G$ is a connected reductive group over $Q$ and
\item for all $\pv\in |C|$ the group $\mathbb G(A_\pv)\subset \mathbb G(Q_\pv)=\mathbf G(Q_\pv)$ is a parahoric subgroup in the sense of \cite[Definition 5.2.6]{BT84}.
\end{itemize}

For each parahoric subgroup in $\mathbf G_\pv(Q_\pv)$ there is a unique smooth, affine group scheme $\mathbb H$ over $A_\pv$ with connected special fiber, with generic fiber equal to $\mathbf G_\pv$ and with $\mathbb H(A_\pv)$ equal to this parahoric subgroup. Since this group scheme is exactly given by $\mathbb G_\pv$ our definition of parahoric  Bruhat-Tits group scheme coincides with the one in \cite[Definition 3.11]{AH13}.

Now Bruhat-Tits group schemes can be constructed as follows. We start with a reductive group scheme $\mathbf G$ over the function field $Q$, which has a reductive model $G$ over an open subscheme $C\setminus\{w_1,\dots, w_m\}$ of $C$. For each of the pairwise different closed points $w\in \underline w:=\{w_1,\dots, w_m\}$ we choose furthermore a parahoric subgroup $H_w\in \mathbb G(Q_\pw)$. Then $H_\pw$ corresponds as explained above to a smooth, affine group scheme $\mathbb H_\pw$ over $A_\pw$ with generic fiber $\mathbf G\times_{Q}Q_\pw$. Consequently $\big(\coprod_{\pw\in\underline w} \mathbb H_\pw \big) \coprod G$ is a group scheme over $\big(\coprod_{\pw\in\underline w} Spec A_\pw \big) \coprod C\backslash\{\pw_1,\dots,\pw_n\}$. 
Using the theorem of Beauville-Laszlo \cite[Theorem in section 3]{BL95} the identification $\mathbb H_\pw\times_{A_\pw}Q_\pw=G\times_{C\setminus\{\underline w\}}Q_\pw$ allows us to glue this group scheme to a group scheme $\mathbb G$ over $C$.
This group scheme is by \cite[Proposition 2.7.1]{Gro65} smooth and by \cite[Proposition 17.7.1]{Gro67} affine over $C$. Therefore $\mathbb G$ is by construction a parahoric Bruhat-Tits group scheme satisfying $\mathbb G_\pv=\mathbb H_\pv$ and $\mathbb G\times_C Q=\mathbf Q$.

Further we remark that if $\pi:\widetilde C\to C$ is a generically \'{e}tale covering of $C$ and $\mathbb G$ is a parahoric Bruhat-Tits group scheme over $\widetilde C$ then by \cite[Example (3) page 2]{Hei10} the Weil restriction $\pi_\star\mathbb G$ (see also lemma \ref{pismooth}) of $\mathbb G$ along $\pi$ is again a parahoric Bruhat-Tits group scheme.
In addition we remark that parahoric Bruhat-Tits group schemes give an interesting class of smooth, affine group schemes over $C$ since moduli spaces of global $\mathbb G$-shtukas for such parahoric Bruhat-Tits group schemes $\mathbb G$ are used by Lafforgue to establish in \cite{Laf12} and \cite{Laf14} the Langlands-parametrization over the function field $Q$.

\section{Functoriality of $\nabla_n^{\hat Z_{v},H}\mathscr H^1(C,\mathbb G)$}

In this section we establish and analyze morphisms between moduli spaces of global $\mathbb G$-shtukas, which are functorial in changing the curve $C$ and the group scheme $\mathbb G$. As mentioned in the introduction, apart from the general interest of these morphisms in the study of $\nabla_n^{\hat Z_\uv, H}\mathscr H^1(C,\mathbb G)$, there are two other motivations. The first motivation concerns a potential formulation of an Andr\'e-Oort conjecture for moduli spaces of global $\mathbb G$-shtukas using the finiteness results in theorem \ref{basechangegeneral} and theorem \ref{theoclosedfinite}. This potential formulation is exlained in more detail in remark \ref{remAO}.
The second motivation arises from the study of stratifications of $\nabla_n^{\hat Z_\uv, H}\mathscr H^1(C,\mathbb G)$ in \cite{Bre19}, where the results of this third section are needed again.\\
The third section is divided into three subsections. In the first subsection we define a shtuka datum and morphisms of these. A shtuka datum contains all the necessary parameters to define a moduli space of $\mathbb G$-shtukas. Then a morphism is defined in such a way that it exactly satisfies the properties to induce a morphism of the corresponding moduli spaces. The fact that $\nabla_n^{\hat Z_\uv, H}\mathscr H^1(C,\mathbb G)$ is indeed functorial in the shtuka datum is then shown in the following two subsections. The second subsection discusses the case where we only change the curve $C$ in the shtuka datum. The induced morphism is constructed and is proven in theorem \ref{basechangegeneral}  to be finite. 
\\
In the third subsection a change of the group scheme $\mathbb G$ by $f:\mathbb G\to \mathbb G'$ is analyzed. Before making any assumptions on $f$ the morphism $\nabla_n^{\hat Z_\uv, H}\mathscr H^1(C,\mathbb G)\to \nabla_n^{\hat Z'_\uv, H'}\mathscr H^1(C,\mathbb G')$ is constructed in general. Then, assuming that $f$ is generically an isomorphism, we prove in theorem \ref{changeparahoric2} a projectivity and a surjectivity result. 
Afterwards, we consider closed immersions of group schemes. In this situation of a closed immersion we prove $\nabla_n^{\hat Z_\uv, H}\mathscr H^1(C,\mathbb G)\to \nabla_n^{\hat Z'_\uv, H'}\mathscr H^1(C,\mathbb G')$ to be unramified (theorem \ref{theounram}) and even finite if $\mathbb G$ is a parahoric Bruhat-Tits group scheme (theorem \ref{theoclosedfinite}). 

\subsection{The Shtuka Datum}\label{shtukdat}
In this subsection, we define the category of Shtuka-data. While we can easily define the objects, we need some further 
\index{Shtuka-datum}
explanations to define the morphisms.

\df{}{\label{datumdef}
A Shtuka-datum is a tuple $(C,\mathbb G, \underline{\pv}, Z_{\underline\pv}, H)$ where 
\begin{itemize}[label=$-$]
\item $C$ is a smooth, projective, geometrically irreducible curve over $\F_q$,
\item $\mathbb G$ is a smooth, affine group scheme over $C$, 
\item $\uv=(\pv_1,\dots,\pv_n)$ is a tuple of $n$ closed points in $C$ (not necessarily disjoint),
\item $\hat Z_{\uv}$ is a bound in the sense of \ref{boundsdef},
\item $H$ is an open, compact subgroup of $\mathbb G(\mathbb A^{\uv})$.
\end{itemize}

}

Before we can define morphisms, we need the following lemmas. Let $\pi:X\to Y$ be a morphism of schemes. We recall that for any functor
$F:\ (\mathbf{Sch}/X)^{op}\longrightarrow\mathbf{Set}$
the push forward $\pi_{\star}F:(\mathbf{Sch}/Y)^{op}\longrightarrow\mathbf{Set}$ with respect to $\pi$ is defined by $(T\to Y)\longmapsto F(T\times_{Y}X)$.
In the case that $F$ is a scheme (i.e. representable) and $\pi_\star F$ is also representable, we call $\pi_\star F$ the Weil restriction $\mathfrak{R}_{X/Y}(F)$ of $F$.
The basic properties and some conditions for the existence of Weil restrictions are discussed and developed in \cite[Paragraph 7.6]{BLR90} and \cite{CGP10}. 
We have the following lemma, where we call a morphism of schemes finite locally free, if it is finite, flat and of finite presentation.

\lem{}{\label{pismooth}
Let $\pi:X\to Y$ be a surjective, finite locally free morphism of schemes, let $\mathbb G$ be a smooth, affine group scheme over $X$ and let $\gtor$ be a $\mathbb G$-torsor on the big \'{e}tale site of $X$, then 
\begin{enumerate}
\item $\pi_\star \mathbb G$ is a smooth, affine group scheme over $Y$
\item $\pi_\star \gtor$ is a $\pi_\star\mathbb G$-torsor on the big \'{e}tale site of $Y$.
\end{enumerate}
}

\prof{
Since $\pi$ is finite and faithfully flat we can apply theorem 4 in \cite[Paragraph 7.6]{BLR90} to see that the Weil restriction $\pi_\star \mathbb G$ exists indeed as a scheme. 
Let $U,V,W\in\ Fun\left((\mathbf{Sch}/X)^{op},\mathbf{Set}\right)$ be arbitrary with natural transformations $f_1:U\to W$ and $f_2:V\to W$, then for $S\in (\mathbf{Sch}/Y)$ we have 
\begin{eqnarray*}
\pi_\star(U\times_W V)(S)&=&Hom_{X}(S\times_Y X, U\times_W V)\\
&=&\left\{(f,g)\ |f\in Hom(S\times_Y X, U),\ g\in Hom(S\times_Y X, V),\ f_1\circ f=f_2\circ g\right\}\\
&=&(\pi_\star U\times_{\pi_\star W} \pi_\star V)(S).
\end{eqnarray*}

This shows that $\pi_\star$ commutes with fiber products and it follows that $\pi_\star \mathbb G$ becomes a group scheme over $Y$. 

Let $U\subset Y$ be an affine open. Then $\pi_\star(X\times_YU)=U$ 
and the compatibility with the fiber product implies $\pi_\star(\mathbb G\times_X X\times_Y U)=\pi_\star \mathbb G \times_Y U$. Now $\pi_\star\mathbb G\times_Y U$ is affine because $\mathbb G$ is affine over $X$ and $\pi$ is finite. Since the Weil restriction of an affine scheme is by construction affine we conclude that $\pi_\star\mathbb G$ is affine over $Y$. Furthermore we know by \cite[Chapter 7.6, Proposition 5]{BLR90} that $\pi_\star \mathbb G$ is again of finite type and smooth over $Y$, which proves the first part.

Now let $\gtor$ be a $\mathbb G$-torsor over $X$. Since $\mathbb G$ is smooth and affine, $\gtor$ is represented by a smooth, affine scheme over $X$, by faithfully flat descent,  \cite[Proposition 2.7.1]{Gro65} and \cite[Proposition 17.7.1]{Gro67}.
\cite[paragraph 7.6, Theorem 4]{BLR90} and \cite[paragraph 7.6, Proposition 5]{BLR90} tell us again, that $\pi_\star \gtor$ is a smooth scheme over $Y$. Using once more the compatibility of the fiber product with the Weil restriction, the action of $\mathbb G$ on $\gtor$ induces an action of $\pi_\star \mathbb G$ on $\pi_\star \gtor$ and additionally the isomorphism $\mathbb G\times_X \gtor\simeq \gtor\times_X \gtor$ yields an isomorphism $\pi_\star \mathbb G\times_Y\pi_\star \gtor\simeq\pi_\star \gtor\times_Y\pi_\star \gtor$. It remains to show that $\pi_\star \gtor$ has \'{e}tale locally on $Y$ a section to $\pi_\star \gtor$. Since $\pi_\star\gtor\to Y$ is smooth and surjective
 this is content of proposition \cite[paragraph 2.2, Prop. 14]{BLR90}. 
}\\

Now morphisms between $\mathbb G$-torsors are sent by $\pi_\star$ to morphisms of $\pi_\star \mathbb G$-torsors and in fact we have the following lemma.

\lem{}{\label{changeiso}
Let $\pi:X\to Y$ be a surjective, finite locally free morphism and $\mathbb G$ a smooth, affine group scheme over $X$. Then the functor 
\begin{alignat*}{1}
\qquad \qquad\pi_\star:\Big\{ \mathbb G\mbox{-torsors on }X\Big\}&\longrightarrow  \Big\{\pi_\star \mathbb G\mbox{-torsors on }Y\Big\}\ ,\\
\gtor&\longmapsto  \pi_\star \gtor 
\end{alignat*}
induced by lemma \ref{pismooth}, is an equivalence of categories.
The inverse functor sends some $\pi_\star\mathbb G$-torsor $\widetilde \gtor$ to $\pi^\star\widetilde \gtor\times^{\pi^\star\pi_\star \mathbb G}\mathbb G$.
}

\prof{
First we prove that $\pi_\star$ is fully faithful. So let $\gtor,\widetilde \gtor$ be two $\mathbb G$-torsors over $X$ and $f':\pi_\star \gtor\to \pi_\star \widetilde \gtor$ be a morphism of $\pi_\star \mathbb G$-torsors. We choose an \'{e}tale covering $U'\to Y$ with $\pi_\star \gtor(U')\neq \emptyset\neq \pi_\star \widetilde \gtor(U')$. Now this implies automatically that $U:=U'\times_Y X\to X$ is an \'{e}tale covering with $\gtor(U)\neq\emptyset\neq \widetilde \gtor(U)$. We choose two sections $u\in \gtor(U)=\pi_\star \gtor(U')$ and $\widetilde u\in \widetilde \gtor(U)=\pi_\star\widetilde \gtor(U')$, which determine trivializations
\begin{alignat*}{3}
\alpha':&\ \pi_\star \gtor\times_Y U'\to \pi_\star \mathbb G\times_Y U'&\qquad\qquad   \widetilde\alpha'&:\pi_\star \widetilde\gtor\times_Y U'\to \pi_\star \mathbb G\times_Y U' \\
\alpha:&\ \gtor\times_X U\to\mathbb G\times_X U&\qquad\qquad   \widetilde\alpha&: \widetilde \gtor\times_X U\to \mathbb G\times_X U
\end{alignat*}
with $\alpha'^{-1}(1)=\alpha^{-1}(1)=u$ and $\widetilde\alpha'^{-1}(1)=\widetilde\alpha^{-1}(1)=\widetilde u$. Now we consider the following diagrams, where $U_2:=U\times_XU$, $U_2':=U'\times_YU'$ with projections $p_1,\ p_2,\ p_1'$ and $p_2'$ and $h:=\widetilde\alpha'\circ(f'\times id_{U'})\circ \alpha'^{-1}$. Note that since $h$ is $\pi_\star\mathbb G$-equivariant $h$ is determined by $h:=h(1)\in \pi_\star\mathbb G(U')=\mathbb G(U)$.
This same $h$ defines then a morphism of $h:\mathbb G\times_X U\to \mathbb G\times_X U$ of $\mathbb G$-torsors on $U$ and we set $f\times id_U:=\widetilde \alpha^{-1}\circ h\circ \alpha:\gtor\times_XU\to\widetilde\gtor\times_X U$.

\begin{equation*}
\xymatrix{
& \pi_\star \mathbb G( U') \ar[rr]^(0,4){h} & & \pi_\star\mathbb G(U')  \\
\pi_\star \gtor( U')   \ar[ru]^{\alpha'} \ar[rr]_(0,3){f'\times id_{U'}}& & \pi_\star \widetilde \gtor( U') \ar[ru]^{\widetilde\alpha'} & \\
& \pi_\star \mathbb G( U_2') \ar@<0.7ex>[rr]^(0,3){p_1'^{\star}h} \ar@<-0.7ex>[rr]_(0,3){p_2'^{\star}h} \ar@{<-}@<0.7ex>[uu]^(0,7){p_1'} \ar@{<-}@<-0.7ex>[uu]_(0,7){p_2'} & & \pi_\star\mathbb G( U_2') \ar@{<-}@<0.7ex>[uu]^{p_1'} \ar@{<-}@<-0.7ex>[uu]_{p_2'} \\
\pi_\star \gtor( U_2') \ar@{<-}@<0.7ex>[uu]^{p_1'} \ar@{<-}@<-0.7ex>[uu]_{p_2'}  \ar@<0.7ex>[ru]^{p_1'^\star\alpha'} \ar@<-0.7ex>[ru]_{p_2'^\star \alpha'} \ar[rr]_{f'\times id_{U_2'}}& & \pi_\star \widetilde \gtor( U_2')  \ar@{<-}@<0.7ex>[uu]^(0,7){p_1'} \ar@{<-}@<-0.7ex>[uu]_(0,7){p_2'} \ar@<0.7ex>[ru]^{p_1'^\star\widetilde\alpha'} \ar@<-0.7ex>[ru]_{p_2'^\star \widetilde\alpha'}& \\
}
\xymatrix{
& \mathbb G( U) \ar[rr]^(0,4){h} & & \mathbb G(U)  \\
\gtor( U )  \ar[ru]^{\alpha} \ar[rr]_(0,3){f\times id_U}& &  \widetilde \gtor( U )\ar[ru]^{\widetilde\alpha} & \\
&  \mathbb G( U_2 ) \ar@<0.7ex>[rr]^(0,3){p_1^{\star}h} \ar@<-0.7ex>[rr]_(0,3){p_2^{\star}h} \ar@{<-}@<0.7ex>[uu]^(0,7){p_1} \ar@{<-}@<-0.7ex>[uu]_(0,7){p_2} & & \mathbb G( U_2) \ar@{<-}@<0.7ex>[uu]^{p_1} \ar@{<-}@<-0.7ex>[uu]_{p_2} \\
 \gtor( U_2 )\ar@{<-}@<0.7ex>[uu]^{p_1} \ar@{<-}@<-0.7ex>[uu]_{p_2}  \ar@<0.7ex>[ru]^{p_1^\star\alpha} \ar@<-0.7ex>[ru]_{p_2^\star \alpha} \ar[rr]_{f\times id_{U_2}}& &  \widetilde \gtor( U_2 ) \ar@{<-}@<0.7ex>[uu]^(0,7){p_1} \ar@{<-}@<-0.7ex>[uu]_(0,7){p_2} \ar@<0.7ex>[ru]^{p_1^\star\widetilde\alpha} \ar@<-0.7ex>[ru]_{p_2^\star \widetilde\alpha}& 
}
\end{equation*}
 
Now by definition of the Weil restriction we have equal sets at the corresponding vertices of the two cubes and the maps $p_1'$, $p_2'$ coincide with $p_1$ and $p_2$.
Furthermore by definition of $\alpha',h,\widetilde\alpha'$ these maps coincide with the maps $\alpha, h,\widetilde \alpha$ in the right hand cube. 
The morphisms $p_1'^\star\alpha':\pi_\star \gtor\times_Y U_2'\to \pi_\star \mathbb G\times_Y U_2'$ and $p_1^\star\alpha$ are uniquely determined by the preimage of $1\in\pi_\star \mathbb G(U_2')=\mathbb G(U_2)$. But this preimage is in both cases given as $p_1^\star(u)=p_1(\alpha^{-1}(1))=p'_1(\alpha'^{-1}(1))\in \pi_\star \gtor(U_2')=\gtor(U_2)$. Hence the maps $p_1'^\star\alpha'$ and $p_1^\star\alpha$ coincide and equally $p_1'^\star h,\ p_1'^\star\widetilde\alpha',\ p_2'^\star\alpha',\ p_2'^\star h,\ p_2'^\star\widetilde\alpha'$ coincide with $p_1^\star h,\ p_1^\star\widetilde\alpha,\ p_2^\star\alpha,\ p_2^\star h$ and $p_2^\star\widetilde\alpha$ respectively.

We further denote $g':=p_2'^\star\alpha'\circ p_1'^\star\alpha'^{-1}(1)\in \pi_\star\mathbb G(U_2')$ and $\widetilde g':=p_2'^\star\widetilde\alpha'\circ p_1'^\star\widetilde\alpha'^{-1}(1)\in \pi_\star\mathbb G(U_2')$. So that we have $g'=g:=p_2^\star\alpha\circ p_1^\star\alpha^{-1}(1)\in \mathbb G(U_2)$ and $\widetilde g'=\widetilde g:=p_2^\star\widetilde\alpha\circ p_1^\star\widetilde\alpha^{-1}(1)\in \mathbb G(U_2)$. 
With these notations we get the following bijections:
\begin{equation*}
\xymatrix{
Hom_{\mathbb G}(\gtor,\widetilde \gtor) \ar[rrr]_{1:1}^{f \mapsto h:=(\widetilde\alpha\circ(f\times id_U)\circ \alpha^{-1})(1)} \ar@{-->}[d] &\qquad\mbox{\qquad  \qquad } \qquad && \{h\in \mathbb G(U)\ |\ \widetilde g\circ p_1^\star h= p_2^\star h\circ g\} \ar@{=}[d]\\
Hom_{\pi_\star\mathbb G}(\pi_\star \gtor,\pi_\star\widetilde \gtor) \ar[rrr]_{1:1}^{f' \mapsto h':=(\widetilde\alpha'\circ(f'\times id_{U'})\circ \alpha'^{-1})(1)} &&& \{h'\in \mathbb \pi_\star G(U')\ |\ \widetilde g'\circ p_1'^\star h= p_2'^\star h\circ g'\}
}.
\end{equation*}
Here the horizontal bijections are due to faithfully flat descent \cite[paragraph 6.1, Theorem 6]{BLR90} and the fact that the condition $\widetilde g\circ p_1^\star h=p_2^\star h\circ g$ is equivalent by definition of $g,\widetilde g$ to $p_1^\star(\widetilde\alpha^{-1}\circ h\circ\alpha)=p_2^{\star}(\widetilde \alpha^{-1}\circ h\circ \alpha)$ and for $\widetilde g'\circ p_1'^\star h'=p_2'^\star h'\circ g'$ respectively. The equality on the right follows from the identifications in the above cubes. To prove the fully faithfulness it remains to show that the bijective dashed arrow is given by $\pi_\star$. By definition of $\pi_\star f$ the following diagram commutes 
\begin{equation*}
\xymatrix{
Hom(U',\pi_\star \gtor) \ar[r]^{\pi_\star f} \ar@{=}[d] & Hom(U',\pi_\star \widetilde \gtor) \ar@{=}[d]\\
Hom(U'\times_Y X, \gtor) \ar[r]^f & Hom(U'\times_Y X,\pi_\star \widetilde \gtor)
}
\end{equation*}
which shows that $f$ and $\pi_\star f$ map to the same $h$ on the right hand side.\\
It remains to show that $\pi_\star$ is essentially surjective. So let $\widetilde \gtor$ be a $\pi_\star {\mathbb G}$-torsor over $Y$ and choose again an \'{e}tale covering $U'\to Y$ and a trivialization $\alpha': \widetilde \gtor\times_Y U'\isom \pi_\star\mathbb G\times_Y U'$.
Let $U_2':=U'\times_Y U'$ and $g':=p_2'^\star\alpha'\circ p_1'^\star\alpha'^{-1}(1)\in\pi_\star\mathbb G(U_2')$, where $p'^\star_2\alpha'\circ p_1'^\star\alpha'^{-1}:\pi_\star\mathbb G\times_Y U_2'\isom \pi_\star\mathbb G\times_Y U_2'$. 
So the descent datum of $\widetilde \gtor$ is isomorphic to $(\pi_\star\mathbb G\times_Y U', g')$.
Now $U:=U'\times_YX\to X$ is an ´\'{e}tale covering and we set $\gtor_U:=\mathbb G\times_XU$ as well as $U_2:=U\times_X U=U_2'\times_YX$ with projections $p_1,p_2$. Let $g\in\mathbb G(U_2)$ be equal to $g'$ using $\mathbb G(U_2)=\pi_\star \mathbb G(U_2')$. Then $(\gtor_U, g)$ is a descent datum that comes by \cite[beginning of paragraph 6.5 and paragraph 6.1, Theorem 6]{BLR90} from a $\mathbb G$-torsor $\gtor$ on $X$. Now it is clear that $(\pi_\star \gtor_U, \pi_\star g)=(\pi_\star\mathbb G\times_YU',g')$. Therefore we have $\pi_\star \gtor\simeq\widetilde \gtor$ which proves that $\pi_\star$ is essentially surjective.
We only need to prove that for every $\mathbb G$-torsor $\gtor$ on $X$ the torsor $\pi^\star\pi_\star \gtor\times^{\pi^\star\pi_\star\mathbb G} \mathbb G$ is isomorphic to $\gtor$. 
With the same notation as above $\gtor$ is given by the $(\mathbb G\times_X U, g)$ and $\pi_\star \gtor$ is given by the descent datum $(\pi_\star\mathbb G\times_Y U', g')$. Restricting the latter torsor to $X$ we get the descent datum $(\pi^\star(\pi_\star\mathbb G \times_YU'), g\times id_X)=(\pi_\star\mathbb G \times_YU'\times_YX, g\times id_X)$. Using the adjunction $Hom_Y(\pi_\star\mathbb G,\pi_\star\mathbb G)=Hom_X(\pi^\star\pi_\star\mathbb G,\mathbb G)$ we denote by $\varphi:\pi^\star\pi_\star\mathbb G\to\mathbb G$ the morphism  corresponding to $id_{\pi_\star\mathbb G}$. Now applying the functor $\times^{\pi^\star\pi_\star\mathbb G, \varphi}\mathbb G$ gives us the descent data $\left(\pi^\star(\pi_\star\mathbb G\times_YU')\times^{\pi^\star\pi_\star\mathbb G}\mathbb G, (g'\times id_X, \mathbf 1_{\mathbb G})\right)$. Since $\varphi$ maps $(g'\times id_X)$ to $g$ this descent data is isomorphic to $(\mathbb G\times_X U, g)$, which proves the lemma. 
}

\mbox{}\\
Now let $\pi:C\to C'$ be a finite morphism from $C$ to some other smooth, projective, geometrically irreducible curve $C'$. This morphism is then automatically faithfully flat \cite[chapter II Prop. 6.8 and chapter III Prop. 9.7]{Har77}. Let further $(C,\mathbb G, \underline{\pv}, Z_{\underline\pv}, H)$ be a shtuka datum as in definition \ref{datumdef}. Since $\mathbb G$ is a smooth, affine group scheme over $C$, this allows us to apply lemma \ref{pismooth} and \ref{changeiso} in this situation.

\rem{}{\label{isoofh} We denote by $\mathscr H^1(C,\mathbb G)$ the category fibered in groupoids, whose $S$-valued points for some $\F_q$-scheme $S$ are given by isomorphy classes of $\mathbb G$-torsors over $C_S$. By lemma \ref{changeiso} $\pi_\star$ induces an isomorphism $\mathscr H^1(C,\mathbb G)\longrightarrow \mathscr H^1(C',\pi_\star\mathbb G)$.}

Let $w_i=\pi(v_i)$ and $\underline w=(w_1,\dots, w_n)$. Our next goal is to define the bound $\pi_\star Z_{\underline v}=(Z_{w_i})_i$ at the points $w_i$. We need the following lemma and general remark, where $w$ is a closed point in $C'$, $A'_w$ is the completion of the local ring $\oo_{C',w}$ and $\pi_{w}=\pi\times id_{Spec\ A'_w}:C\times_{C'}Spec\ A_w=\coprod_{v|w}Spec\ A_v\to Spec\ A'_w$.

\rem{}{\label{remaweil}
In the next lemma, we need the following general fact about Weil restrictions. Let $X,Y,S$ be schemes over some base scheme $Z$, $\pi:X\to Y$ a $Z$-morphism, $M$ an $X$-scheme and $Y_S=Y\times_ZS$, $X_S=X\times_ZS$ and $M_S=M\times_ZS$ the appropriate base changes. 
Then we have $(\pi\times id_S)_\star(M\times_{Z} S)=\pi_\star M\times_{Z}S$. This is easily seen by the equation for $T\in\mathbf{Sch}/Y_S$:
\begin{eqnarray*}
(\pi_\star M\times_{Z} S)(T) &=& (\pi_\star M\times_{Y} Y_S)(T) = 
 Hom_{Y}(T,\pi_\star M)
= Hom_X(T\times_{Y} X,M)\\ &=& Hom_{X_S}(T\times_{Y}X,M_S)=(\pi\times id_S)_\star(M\times_{Z} S)(T).
\end{eqnarray*}
}

\lem{}{We have $(\pi_{\pw})_\star(\coprod_{v\in\pi^{-1}(\pw)}{\mathbb G}_\pv)
=({\pi_\star\mathbb G})_\pw$ as a group scheme over $Spec A'_\pw$.}
\prof{
This follows formally from remark \ref{remaweil} with $M=\mathbb G$, $X=C$, $Y=Z=C'$ and $S=Spec\ A_\pw$ since we have
\begin{equation*}
(\pi_w)_\star(\coprod_{v|w}{\mathbb G}_{v})=(\pi_w)_\star(\mathbb G\times_{C}\coprod_{v|w} Spec\ A_{v})=(\pi_w)_\star(\mathbb G\times_{C'} Spec\ A_{w})
=\pi_\star \mathbb G\times_{C'}Spec\ A_\pw=({\pi_\star\mathbb G})_\pw.
\end{equation*}
}

The notation $\widetilde {\mathbb G_\pv}$ was introduced in \ref{beauvillelaszlo} and denotes the Weil restriction of the group scheme $\mathbb G_\pv$ along $Spec\ A_\pv\to Spec\ \F_q\sem{z_\pv} $. The lemma has the following corollary.

\ko{}{\label{pgruppenprodukt}
We have $\ds\prod_{\pv|\pw} L^+\widetilde{\mathbb G_{\pv}}=L^+ \widetilde{(\pi_\star\mathbb G)_\pw}$ as group schemes over $\F_q$.}
\prof{
Let $R$ be a connected $\F_q$-algebra,  
then we have: 
\begin{align*}
L^+\widetilde{(\pi_\star\mathbb G)_\pw}(R)
&=\widetilde{(\pi_\star\mathbb G)_\pw}(R\sem{z_w})
=Hom_{Spec\ A_\pw}(Spec\ R\sem{z_w}\otimes_{\F_q}\F_\pw, (\pi_\star \mathbb G)_\pw)
\\
&=Hom_{Spec\ A_\pw}(Spec\ R\hat\otimes_{\F_q}A_\pw, (\pi_\star \mathbb G)_\pw)\\
&=Hom_{Spec\ A_\pw}(Spf\ R\hat\otimes_{\F_q}A_\pw, (\pi_w)_\star(\coprod_{v|w}{\mathbb G}_{v}))
\\
&= Hom_{\coprod_{\pv|\pw}Spec\ A_\pv}(Spec\ R\hat\otimes_{\F_q}A_\pw\otimes_{A_\pw}\prod_{\pv|\pw}A_\pv, \coprod_{\pv|\pw}{\mathbb G_\pv})\\
&=\prod_{\pv|\pw}Hom_{Spec\ A_\pv}(Spec\ R\sem{z_\pv} \otimes_{\F_q}\F_\pv, {\mathbb G_\pv})
=\prod_{\pv|\pw}\widetilde{\mathbb G_\pv}(R\sem{z_\pv})
\\
&=\prod_{\pv|\pw}L^+\widetilde{\mathbb G_\pv}(R).
\end{align*}
}

We have the following $2$-cartesian diagrams:
\begin{equation*}
\xymatrix{
\underset{\pv|\pw}{\prod}{\mathcal F\it l}_{\widetilde{\mathbb G_\pv}} \ar[r] \ar[d] & \F_q \ar[d] & {\mathcal F\it l}_{\widetilde{(\pi_\star \mathbb G)_\pw}} \ar[d] \ar[r] & \F_q \ar[d] \\
\underset{\pv|\pw}{\prod} \mathscr H^1(\F_q,L^+\widetilde{\mathbb G_\pv}) \ar[r] & \underset{\pv|\pw}{\prod}\mathscr H^1(\F_q,L\widetilde{\mathbb G_\pv}) & \mathscr H^1(\F_q, L^+\widetilde{(\pi_\star\mathbb G)_\pw}) \ar[r] & \mathscr H^1(\F_q, L\widetilde{
(\pi_\star\mathbb G)_\pw})
}
\end{equation*}

By corollary \ref{pgruppenprodukt} the lower stacks in the diagrams are isomorphic, so that we get an isomorphism $\prod_{\pv|\pw} {\mathcal F\it l}_{\widetilde{\mathbb G_\pv}}\isom {\mathcal F \it l}_{\widetilde{(\pi_\star\mathbb G)_\pw}}$ and by the base change with the compositum $\F$ of the finite fields $\F_{\pv}$ for all $\pv|\pw$ we get an isomorphism $\prod_{\pv|\pw}\prod_{l\in\Z/deg\ \pv} {\mathcal F\it l}_{{\mathbb G_\pv}}\times_{\F_\pv}\F\simeq\prod_{l\in \Z/deg\ \pw}{\mathcal F\it l}_{(\pi_\star\mathbb G)_\pw}\times_{\F_\pw}\F$. Since $\sigma^{deg\ \pw}$ invariant components are mapped to $\sigma^{deg\ \pw}$ invariant components, it restricts to an isomorphism
\begin{equation*}
\displaystyle \prod_{\pv|\pw}\ \prod_{\underset{\ deg\ \pw|l}{l\in\Z/deg\ \pv}} {\mathcal F\it l}_{\mathbb G_\pv} \times_{\F_\pv} \F\simeq 
{\mathcal F\it l}_{(\pi_\star\mathbb G)_\pw}\times_{\F_\pw}
\F.
\end{equation*}

Now let $R$ be a DVR with $\F\subset R$ and such that there exists a representative $\hat Z_{v, R}$ of $\hat Z_{ v}$ for all $v\in\underline v$. Consider the ind-closed subscheme 
\begin{equation*}
\prod_{\pv|\pw}\ \prod_{\underset{\ deg\ \pw|l}{l\in\Z/deg\ \pv}} {\hat Z}_{\pv,l}\subset \prod_{\pv|\pw}\ \prod_{\underset{\ deg\ \pw|l}{l\in\Z/deg\ \pv}}{\mathcal F\it l}_{\mathbb G_\pv}\times_{\F_\pv} Spf\ R
\end{equation*}
 where ${\hat Z}_{\pv,l}$ is always the closed stratum $\mathcal S(1)\times_{\F_\pv}Spf\ R$ except for $v\in\underline v$ and $l=0$, where we set ${\hat Z}_{\pv_i,0}={\hat Z}_{\pv_i, R}$. 
Here $\mathcal S(1)$ denotes the closed Schubert cell $1\cdot L^+\mathbb G_v$ in ${\mathcal F\it l}_{\mathbb G_\pv}$.
Via the previous isomorphism this defines an ind-closed subscheme in ${\mathcal F\it l}_{(\pi_\star\mathbb G)_\pw}\times_{\F_\pw}Spf\ R$. This defines a bound ${\hat Z}_w$ in the sense of \ref{boundsdef} in ${\mathcal F\it l}_{(\pi_\star\mathbb G)_\pw}$ and we set \label{defpiz} $\pi_\star {\hat Z}_{\underline v}:={\hat Z}_{\underline w}:=({\hat Z}_{w_i})_{i}$.\\

Next we define $\pi_\star H$. We recall that $H$ was an open, compact subgroup of $\mathbb G(\mathbb A^{\underline\pv})$. 
Since $\underline \pv\subset \pi^{-1}(\underline\pw)\subset |C|$ we have a quotient map of topological rings $\mathbb A^{\uv}\longrightarrow \A^{\pi^{-1}(\underline\pw)}$.

Since this map is open, it induces by \cite[theorem 3.6]{Con12} an open continuous group homomorphism $\mathbb G(\A^{\uv})\longrightarrow \mathbb G(\A^{\pi^{-1}(\underline\pw)})$. We have $\A^{\underline\pw}\times_{C'}C=\A^{\underline\pw}\times_{\eta'}\eta=\A^{\pi^{-1}(\underline\pw)}$ where $\eta$ and $\eta'$ are the generic points of $C$ and $C'$. This gives us with the definition of the Weil restriction $\pi_\star \mathbb G(\A^{\underline\pw})=\mathbb G(\A^{\pi^{-1}(\underline \pw)})$, where both groups carry the same topology by \cite[example 2.4]{Con12}. 
Now the image of $H$ under this morphism gives us an open, compact subgroup in $\pi_\star\mathbb G(\A^{\underline\pw})$ that we denote by $\pi_\star H$.\\

\rem{}{
We have seen in \ref{dlevel} that
there is the possibility to define level structures using finite closed subschemes $D$ of $C$ and  in \ref{hlevel} we have remarked that $D$-level structures of a $\mathbb G$-shtuka correspond bijectively to $H_D$-level structures, where $H_D=ker(\mathbb G(\mathbb O^\uv)\to \mathbb G(\mathcal O_D))$. Now we can also consider the Weil restriction $\pi_\star D$ of $D$. It is a closed finite subset of $C'$ consisting of the points $\{\pw\in |C'|\ |\ \pw\times_{C'}C\subset D\}$. And with a $D$-level structure of some $\mathbb G$-shtuka $\underline{\mathcal G}$ we could associate a $\pi_\star D$-level structure of the corresponding $\pi_\star\mathbb G$-shtuka $\pi_\star\underline{\mathcal G}$ (which will be defined in proposition \ref{finitemorbasechange}). But compared to the associated $\pi_\star H$-level structure that we will define in theorem \ref{basechangegeneral} we would lose some information at the points $D\backslash(\pi_\star D\times_{C'}C)$, which is seen in the following way. Since we have $\pi_\star D\times_{C'}C\subset D$ we have $H_D\subset H_{\pi_\star D\times_{C'}C}$ and hence $\pi_\star H_D\subset \pi_\star H_{\pi_\star D\times_{C'}C}$. Now 
\begin{alignat*}{1}
H_{\pi_\star D}&=ker\big(\pi_\star \mathbb G(\mathbb O^{\underline\pw})\to \pi_\star\mathbb G(\mathcal O_{\pi_\star D})\big)=ker\big(\mathbb G(\mathbb O^{\pi_{-1}{\underline\pw}})\to \mathbb G(\mathcal O_{\pi_\star D\times_{C'}C})\big)\\
&=H_{\pi_\star D\times_{C'}C}\cap \mathbb G(\A^{\pi^{-1}(\underline\pw)})=im\left(H_{\pi_\star D\times_{C'}C}\to\mathbb G(\A^{\pi^{-1}(\underline w)})\right)=:\pi_\star H_{\pi_\star D\times_{C'}C}\supset \pi_\star H_D
\end{alignat*}
shows that $\pi_\star H_D$ is in general a finer level than $\pi_\star D$ (or equivalently $H_{\pi_\star D}$) and the previous equation shows that the information is lost exactly at the points $D\backslash (\pi_\star D\times_{C'} C)$.
}

All these previous explanations concerned the case that we change the curve in the shtuka datum but we can also change the group scheme in this datum. Let $f:\mathbb G\to \mathbb G'$ be any morphism of smooth, affine group schemes over $C$ and $\pv$ a closed point in $C$.  Firstly this induces a morphism $L^+\mathbb G_\pv\to L^+\mathbb G'_\pv$ of the positive loop groups as well as a morphism $L\mathbf G_\pv\to L\mathbf G'_\pv$ of the loop groups. Consequently we also get a morphism $\mathcal F\mathit l_{\mathbb G_\pv}\to \mathcal F\mathit l_{\mathbb G'_\pv}$ of the affine flag varieties.
Secondly such a morphism induces a morphism $f_{\A^{\uv}}:\mathbb G(\A^{\uv})\to \mathbb G'(\A^{\uv})$ of locally compact Hausdorf spaces by \cite[Proposition 2.1]{Con12}.
Now we can define morphisms of shtuka data.

\df{}{\label{defshtukadata}\mbox{}
A morphism between two shtuka data $(C,\mathbb G, \uv, \hat Z_{\uv}, H)$ and $(C',\mathbb G', \underline\pw, \hat Z'_{\underline\pw}, H')$ is a pair $(\pi,f)$ such that:
\begin{itemize}[label=$-$]
 \item $\pi:C\to C'$ is a finite morphism with $\pi(\pv_i)=\pw_i$ 
\item $f:\pi_\star\mathbb G\to \mathbb G'$ is a morphism of smooth, affine group schemes over $C'$
\item The morphism $\ds(\pi_\star Z_{\uv})_R\to \prod_{\pw\in\underline \pw}{\hat\fl}_{(\pi_\star \mathbb G)_\pw,R}\to \prod_{\pw\in\underline \pw}{\hat\fl}_{\mathbb G'_\pw,R}$ factores through $\hat Z'_{\underline \pw,R}$, where $R$ is a DVR such that there exists representatives $(\pi_\star \hat Z_{\uv})_R$ and $\hat Z'_{\underline \pw,R}$ of the corresponding bounds
\item $f_{\A^{\underline\pw}}(\pi_\star H)\subset H'$ 
\end{itemize}

}

With this definition we have reached the goal of this subsection. In the next two subsections we will prove that such a morphism induces a morphism of the corresponding moduli stacks and determine some of its properties. But before we give some remarks.

\rem{}{\label{remshtukadata}\mbox{}

\begin{itemize}[label=$-$]
 \item Let $\pi:C\to C'$ be a finite morphism and $\underline\pw=\pi(\uv)$. With the definition of $\pi_\star H$ and $\pi_\star Z_\uv$ on page \pageref{defpiz} it is clear that $(\pi,id_{\pi_\star\mathbb G}):(C,\mathbb G, \uv, Z_{\uv}, H)\to (C',\pi_\star\mathbb G, \underline\pw, \pi_\star Z_{\underline\pv}, \pi_\star H)$ defines a morphism of shtuka data
\item Every morphism $(\pi,f)$ of shtuka data factorizes as $(id_C,f)\circ(\pi,id_{\pi_\star \mathbb G})$. 
\item If $f:\pi_\star \mathbb G\to\mathbb G'$ is an isomorphism in the generic fiber we have $\pi_\star \mathbb G(\A^{\underline\pw})= \mathbb G'(\A^{\underline\pw})$ so that we can naturally choose $H=H'$.
\item If $f:\pi_\star \mathbb G\to\mathbb G'$ is smooth in the generic fiber, then $f_{\A^{\underline\pw}}:\pi_\star \mathbb G(\A^{\underline\pw})\to \mathbb G'(\A^{\underline\pw})$ is an open map by \cite[Theorem 4.5]{Con12} so that we can naturally choose $H'=f_{\A^{\pw}}(\pi_\star H)$
\item If $f:\pi_\star \mathbb G\to\mathbb G'$ is proper in the generic fiber, then $f_{\A^{\underline\pw}}:\pi_\star \mathbb G(\A^{\underline\pw})\to \mathbb G'(\A^{\underline\pw})$ is a topologically proper map by \cite[Proposition 4.4]{Con12} so that we can naturally choose $H=f_{\A^{\pw}}^{-1}(H')$ 
\end{itemize}

}

\subsection{Changing the Coefficients}\label{changeofcurve}

In this subsection we prove that a morphism of shtuka data  $(\pi, id)$, where we only change the curve, induces a finite morphism of the corresponding moduli stacks. We firstly prove this for the moduli stack $\sa$, where the characteristic sections are not fixed and no boundedness condition or level structures are imposed. For this purpose we need the following lemma:

\lem{}{\label{CSdense}
Let $S$ be an $\F_q$-scheme with $n$ morphisms $s_i:S\to C$ for $i=1,\dots, n$. Then the scheme theoretic image of $\widetilde C_S:=C_S\backslash \bigcup_i \Gamma_{s_i}$ in $C_S$ equals $C_S$.
}
\prof{
Since $D:=\bigcup_i \Gamma_{s_i}$ is an effective Cartier-Divisor on $C_S$ over $S$, we find an affine covering $(U_j)_{j\in J}$ of $C_S$ with $U_j:=Spec\ B_j$ such that $D$ is the vanishing locus of an element $f_j\in B_j$ that can be written as $f_j=\frac{a_j}{b_j}$ with two regular elements $a_j, b_j\in B_j$ (see \cite[after definition 11.24]{GW10}). Now the ring homomorphism $B_j\to \Gamma(U_j\backslash D, \mathcal O)$ is injective, which is seen as follows. An element $x\in B_j$ is send to 0 if and only if $f_j^mx=0$ for some $m\in \N$. The latter condition implies $a_j^mx=0$ and since $a_j$ is a non-zero divisor this means $x=0$ so that $Spec\ B_j\backslash D$ is schematically dense in $Spec\ B_j$ (compare also \cite[Lemma 20.2.9]{Gro67}).
Now gluing all the $U_j$ shows that for every affine open $V\subset C_S$ the ring homomorphism $\Gamma(V,\mathcal O)\to \Gamma(V\backslash D, \mathcal O)$ is injective and we conclude that $\widetilde C_S$ is schematically dense in $C_S$ (see \cite[\mbox{} 20.2.1]{Gro67}).
}\mbox{}\\

\prop{}{\label{finitemorbasechange}
Let $\pi:C\to C'$ be a finite morphism of smooth, projective, geometrically irreducible curves over $\F_q$ and let $\mathbb G$ be a smooth, affine group scheme over $C$. This induces a finite morphism of the moduli stacks 
\begin{equation*}
    \pi_\star :\ \sa \to \nabla_n\mathscr H^1(C',\pi_\star\mathbb G) .
\end{equation*}
 which factors through a closed immersion
 $\nabla_n\mathscr H^1(C,\mathbb G)\hookrightarrow\nabla_n\mathscr H^1(C',\pi_\star \mathbb G)\times_{C'^n}C^n$.
}

\prof{ Let $S$ be an $\F_q$-scheme and $(\gtor,s_1,\dots, s_n, \tau_\gtor)\in \nabla_n\mathscr H^1(C,\mathbb G)(S)$. We describe its image $(\gtor',s_1',\dots,s_n',\tau_{\gtor'})$ to define the morphism. The torsor $\gtor'$ is given by $(\pi_S)_\star \gtor$ and the sections $s_i:S\to C$ are mapped to the composition $s_i':=\pi\circ s_i:S\to C'$. This implies $\pi_S(\bigcup_i \Gamma_{s_i})\subset \bigcup_i \Gamma_{s'_i}\subset C'_S$. Let $\widetilde C'_S=C'_S\backslash (\cup_i\Gamma_{s'_i})$ and $\widetilde{C_S}=C_S\backslash (\cup_i\Gamma_{s_i})$. 
Then $U:=C\times_{C'}\widetilde{C'_S}=C_S\times_{C'_S}\widetilde{C'_S}$ is open in $\widetilde{C_S}$. We denote by $\pi_{U}:=\pi\times_{id_{C'}}id_{\widetilde{C'_S}}:U\to \widetilde{C'_S}$ and we have $(\pi_U)_\star(\mathbb G\times_C U)=\pi_\star\mathbb G\times_{C'}\widetilde{C'_S}$. Now we restrict $\tau_\gtor:\sigma^\star \gtor|_{\widetilde {C_S}}\to \gtor|_{\widetilde {C_S}}$ to $U_S$ and apply lemma \ref{changeiso} to $\pi_U$. The category equivalence gives us the desired morphism $\tau_{\gtor'}:(\pi_U)_\star(\sigma^\star \gtor\times_{{C'_S}}\widetilde{C_S'})=\sigma^\star \gtor'|_{\widetilde{C_S'}}\to (\pi_U)_\star(\gtor\times_{{C'_S}}\widetilde{C'_S})= \gtor'|_{\widetilde{C_S'}}$. This defines a global $\pi_\star\mathbb G$-shtuka $(\gtor',s_1',\dots,s_n', \tau_\gtor')$ over $S$ and therefore the morphism of the moduli stacks. 

We now show that this morphism is representable and finite.
Let $S$ be again an arbitrary scheme over $\F_q$ and $\underline \gtor':S\to \nabla_n\mathscr H^1(C',\pi_\star\mathbb G)$ be given by $\underline \gtor'=(\gtor',\tau', s_1',\dots, s_n')$. Then ${\nabla_n \mathscr H^1(C,\mathbb G)}\times_{\nabla_n\mathscr H^1(C',\pi_\star\mathbb G)}S$ is the category fibered in groupoids over $\mathbf{Sch}/\F_q$ whose fiber category over an $\F_q$-scheme $T$ is given by
\begin{equation*}
\big\{
(\underline \gtor,g:T\to S,\beta)\quad  \big|\ \  \underline \gtor=(\gtor,\tau_\gtor,s_1,\dots,s_n)\in \nabla_n \mathscr H^1(C, \mathbb G)(T)\mbox{ and } \beta:g^\star\underline \gtor'\isom \pi_\star \underline \gtor
 \big\}.
\end{equation*}

Using the $n$ sections $s_1',\dots,s_n':S\to C'$ and the morphism $\pi:C\to C'$ we set $\widetilde S:=S\times_{C'^n}C^n$. 
Since $S \times_{\mathscr H^1 (C',\pi_\star \mathbb G)}  \mathscr H^1(C,\mathbb G)=S$ by remark \ref{isoofh} we know that $Hom_{\F_q}(T,S)$ is in bijection with the tuples $\{ g:T\to S, \gtor:T\to \mathscr H^1(C,\mathbb G), \alpha:g^\star \gtor'\isom \pi_\star \gtor \}$.
Consequently the $\F_q$-morphisms $T\to \widetilde S$ are in bijection with the tuples $(\gtor, s_1,\dots, s_n,g,\alpha)$, where $(\gtor,g,\alpha)\in S(T)$ as before and $s_1,\dots,s_n:T\to C$ are morphisms making the following diagram commutative for all $i=1,\dots,n$
\begin{equation*}
\xymatrix{
T \ar[d]_{s_i} \ar[r]^g & \ar[d]^{s_i'} S\\
C \ar[r]^\pi & C'
\ .}
\end{equation*}

We claim that we get a morphism $\nabla_n\mathscr H^1(C, \mathbb G)\times_{\nabla\mathscr H^1(C',\pi_\star\mathbb G)}S\to \widetilde S$ that is injective on $T$-valued points (hence a monomorphism) and satisfies the valuative criterion for properness. Then this implies by \cite[Proposition 8.11.5]{Gro66} that $\nabla_n\mathscr H^1(C, \mathbb G)\times_{\nabla\mathscr H^1(C',\pi_\star\mathbb G)}S$ is a closed subscheme of $\widetilde S$. So first of all a given object $(\gtor,\tau_\gtor, s_1,\dots, s_n, g, \alpha)$ in $(\nabla_n\mathscr H^1(C, \mathbb G)\times_{\nabla\mathscr H^1(C',\pi_\star\mathbb G)}S)(T)$ is sent to $(\gtor, s_1,\dots,s_n,g,\alpha)$. Since $\alpha:g^\star\underline \gtor'\isom \pi_\star\underline \gtor$ is not only an isomorphism of torsors, but also of $\pi_\star\mathbb G$-shtukas the $n$-sections $g\circ s_i'$ of $g^\star\underline \gtor'$ and $(s_i\circ\pi)$ of $\pi_\star \underline \gtor$ have to coincide, so that $(\gtor,s_1,\dots,s_n,g,\alpha)$ is a well defined object in $\widetilde S(T)$. This induces the morphism $\nabla_n\mathscr H^1(C, \mathbb G)\times_{\nabla\mathscr H^1(C',\pi_\star\mathbb G)}S\to \widetilde S$. Further it was claimed, that this morphism is injective on $T$-valued points. So given two points $(\gtor,\tau_\gtor, s_1,\dots,s_n,g, \alpha)$ and $(\gtor,\widetilde\tau_\gtor,s_1,\dots,s_n,g,\alpha)$ we need to show that this implies $\tau_\gtor=\widetilde \tau_\gtor$. Since $\alpha$ is an isomorphism of $\pi_\star\mathbb G$-shtukas, we have $\alpha^{-1}\circ\pi_\star\tau_\gtor=g^\star\tau'\circ\sigma^\star\alpha^{-1}=\alpha^{-1}\circ\pi_\star\widetilde{\tau_\gtor}:\sigma^\star(\pi_\star \gtor)\big|_{\widetilde {C_T'}}\isom g^\star \gtor'\big|_{\widetilde{C_T'}}$, where we write again $\widetilde{C'_T}:=C_T'\backslash \bigcup_i \Gamma_{s_i'} $. This implies $\pi_\star\tau_\gtor=\pi_\star\widetilde {\tau_\gtor}$ and using lemma \ref{changeiso} applied to $\pi\times id_{\widetilde {C_T'}\times_{C'}C}$ we see $\tau_\gtor\big|_{\widetilde{C_T'}\times_{C'}C}=\widetilde{\tau_\gtor}\big|_{\widetilde{C_T'}\times_{C'}C}$. 
We even need to know that $\tau_\gtor=\tau'_\gtor$. In the following diagram the restriction of $(\tau_\gtor, \tau'_\gtor)$ to $\widetilde{C'_T}\times_{C'}C$ factors by the previous observation through the diagonal $\Delta$. 
\begin{equation*}
\xymatrix{
\sigma^\star\gtor\big|_{{\widetilde C}_T} \ar[rr]^{(\tau_\gtor, \tau'_\gtor)} && \gtor\big|_{{\widetilde C}_T}\times_{{\widetilde C}_T} \gtor\big|_{{\widetilde C}_T} \\
\sigma^\star\gtor\big|_{\widetilde{C_T'}\times_{C'}C} \ar[rr] \ar[u] && \gtor\big|_{{\widetilde C}_T} \ar[u]_{\Delta}\quad .
}
\end{equation*}
Since $\gtor$ is separated over $C_T$ the diagonal is a closed immersion and $(\tau_\gtor, \tau'_\gtor)$ factors already over $\widetilde C_T$ through the diagonal if the scheme theoretic image of $\sigma^\star\gtor\big|_{\widetilde{C_T'}\times_{C'}C}$ in $\sigma^\star\gtor\big|_{{\widetilde C}_T} $ equals $\sigma^\star\gtor\big|_{{\widetilde C}_T}$. Since taking the scheme theoretic image is stable under flat base change by \cite[Th\'{e}oreme 11.10.5]{Gro66}, this is the case if the scheme theoretic image of $\widetilde{C_T'}\times_{C'}C$ in $\widetilde C_T$ equals $\widetilde C_T$. By the same argument this follows if the scheme theoretic image of $\widetilde C'_T$ in $C'_T$ equals $C_T'$.
Now this is content of lemma \ref{CSdense} so that we can conclude as desired $\tau_\gtor=\tau'_\gtor$.

Next we claimed that the morphism satisfies the valuative criterion for properness.
So let 
\begin{equation*}
\xymatrix{
Spec\ K \ar[d]_j \ar[rrrr]^{(\mathcal H,\tau_{\mathcal H}, r_1,\dots, r_n, f,\beta)\quad} &&&& S\times_{\nabla_n \mathscr H^1(C',\pi_\star\mathbb G)}\nabla_n\mathscr H^1(C,\mathbb G) \ar[d] & \\
Spec\ R \ar@{-->}[rrrru] \ar[rrrr]_{(\gtor,s_1,\dots,s_n,g,\alpha)} &&&& \widetilde S \ar[r]_{finite} &  S  
}
\end{equation*}
be a commutative diagram, where $R$ is a complete discrete valuation ring with fraction field $K$, maximal ideal $\mathfrak m$ and algebraically closed residue field ${\kappa_R}=R/\mathfrak m$. Note that $R$ is a ${\kappa_R}$-algebra. We have to prove that there exists a unique dashed arrow making everything commutative.
The commutativity of the square shows $\mathcal H=j^\star \gtor$, $s_i\circ j=r_i:Spec\ K\to C$, $f=g\circ j:Spec\ K\to S$ and $j^\star \alpha=\beta$. To define this dashed arrow we have to extend $(\gtor,s_1,\dots, s_n)$ to a $\mathbb G$-shtuka $(\gtor,\tau_\gtor,s_1,\dots, s_n)$ over $R$ such that $\alpha$ extends to an isomorphism $\alpha:g^\star \underline \gtor'\to \pi_\star \underline \gtor$ and $j^\star\tau_\gtor=\tau_{\mathcal H}$. So we define this isomorphism $\tau_\gtor:\sigma^\star \gtor|_{\widetilde{C_R}}\to \gtor|_{\widetilde{C_R}}$. Since $\mathscr H^1(\widetilde {C'_R},\pi_\star\mathbb G)$ and $\mathscr H^1(\widetilde {C'_R}\times_{C'}C,\mathbb G)$ are isomorphic, $\tau_\gtor|_{\widetilde {C'_R}\times_{C'}C}$ is defined by $\alpha\circ g^\star\tau_{\gtor'}\circ \sigma^\star\alpha^{-1}$. Furthermore we know that $\tau_\gtor$ is defined on the generic fiber $\widetilde {C_K}\subset \widetilde {C_R}$ by $\tau_\gtor|_{\widetilde{C_K}}=j^\star\tau_\gtor=\tau_{\mathcal H}$. So let $p\in \widetilde {C_R}\backslash(\widetilde{C'_R}\times_{C'}C\cup \widetilde{C_K})$, i.e. $p\in \big((\bigcup_i\Gamma_{s_i'}\times_{C'}C)\backslash \bigcup_i\Gamma_{s_i}\big)\bigcap C_{R/\mathfrak m}$. It remains to show, that $\tau_\gtor$ extends to $p$. Since $p$ is closed we choose an open $V\subset \widetilde C_{R/\mathfrak m}$ with $V\bigcap(\bigcup_i\Gamma_{s_i'}\times_{C'}C)=p$ and set $\widetilde V=V\backslash p$. Then we consider the 2-cartesian diagram of stacks fibered over ${\kappa_R}_{\acute Et}$ (compare \cite[Lemma 5.1]{AH14}):
\begin{equation*}
\xymatrix{
\mathscr H^1(V,\mathbb G) \ar[r] \ar[d]_{L^+_p}& \mathscr H^1_e(\widetilde V, \mathbb G) \ar[d]_{L_p}  \\
\mathscr H^1({\kappa_R}, L^+\mathbf G_p) \ar[r] & \mathscr H^1({\kappa_R}, L\mathbb G_p)\ .
}
\end{equation*}
Here $\mathscr H^1_e(\widetilde V, \mathbb G)(X)$ is the full subcategory of $\mathscr H^1(\widetilde V, \mathbb G)(X)$ consisting of those $\mathbb G$-torsors over $\widetilde V_X:=\widetilde V\times_{{\kappa_R}} X$ that can be extended to a $\mathbb G$-torsor over $V_X$.
Now $\sigma^\star \gtor|_{V_R}$ and $\gtor|_{V_R}$ define two $R$-valued points in $\mathscr H^1(V,\mathbb G)$ and $\tau_\gtor|_{\widetilde V_R}$ is an isomorphism in $\mathscr H^1_e(\widetilde V,\mathbb G)(R)$ that is already defined. Since $R$ has algebraically closed residue field, we can choose trivializations $\alpha_1:L_p^+(\sigma^\star \gtor)\to L^+\mathbb G_p$ and $\alpha_2:L_p^+(\gtor)\to L^+\mathbb G_p$ (\cite[Proposition 2.4]{AH14}). Then $\alpha_2^{-1}\circ \tau_\gtor\circ \alpha_1:L\mathbf G_{p,R}\to L\mathbf G_{p,R}$ is an isomorphism in $\mathscr H^1({\kappa_R}, L\mathbf G)(R)$ given by an element $h\in L\mathbf G_{p}(R)$. We know by assumption that the pull back of $h$ to $K$ is given by an element $h_K\in L^+\mathbb G_p(K)$, since $\tau_\gtor$ is generically already an isomorphism over $V$. But since $L^+\mathbb G_p$ is closed in $L\mathbb G_p$ it follows that $h\in L^+\mathbb G_p(R)$. This implies that the isomorphism $\tau_\gtor|_{\widetilde V}$ comes from an isomorphism in $\mathscr H^1(V,\mathbb G)(R)$. So $\tau_\gtor$ extends uniquely to $p$ and the valuative criterion is proved. This proves that $\nabla_n\mathscr H^1(C, \mathbb G)\times_{\nabla\mathscr H^1(C',\pi_\star\mathbb G)}S$ is a closed subscheme of $\widetilde S$ and since $\widetilde S$ is finite over $S$ it proves as well that  $\sa \to \nabla_n\mathscr H^1(C',\pi_\star\mathbb G)$ is a finite morphism.
}
\mbox{}\\

The next goal is to prove that for any shtuka datum $(C,\mathbb G, \underline v, \hat Z_{\underline \pv}, H)$ the morphism $\pi:C\to C'$ induces also a finite morphism $\nabla_n^{\hat Z_{\underline v}, H}\mathscr H^1(C, \mathbb G)\to \nabla_n^{\pi_\star \hat Z_{\underline \pv}, \pi_\star H}\mathscr H^1(C', \pi_\star\mathbb G)$. For this we need the following lemma that concerns the boundedness condition. Given a global $\mathbb G$-shtuka $\underline \gtor$ in $\sa$ over $S$, we recall that we introduced in \ref{parglobloc} the global-local functor $\Gamma_{\pv_i}$ that associates with it a local $\mathbb G_{\pv_i}$-shtuka $\Gamma_{\pv_i}(\underline \gtor)$ over $S$. On the other hand we explained (compare also \cite[Remark\ 5.6]{AH14}) that base change with $Spf\ A_{\pv_i}\times_{\F_q}S\simeq \prod_{l\in\Z/deg\ v_i}V(\mathfrak a_{v_i,l})$ gives a local $\widetilde {\mathbb G_{v_i}}$-shtuka $L_{v_i}^+(\underline \gtor)$ over $S$. Here $\widetilde {\mathbb G_{v_i}}$ denotes the Weil restriction $Res_{A_{v_i}/{\F_q}\sem{z_{v_i}}}\mathbb G_{v_i}$.
Now let $Z_{v_i}$ be a bound in $\hat{\mathcal F\it l}_{\mathbb G_{v_i}}$ and $R$ an DVR over $A_{\pv_i}=\F_{\pv_i}\sem{z_{\pv_i}}$ with a representative $\hat Z_{\pv_i,R}\subset \hat{\mathcal F\it l}_{\mathbb G_{\pv_i},R}$. We have $\hat{\mathcal F\it l}_{\widetilde{\mathbb G_{\pv_i}}}\hat\times_{\F_q\sem{z_{\pv_i}}}Spf\ R=\displaystyle\prod_{\Z/deg\ {v_i}}\hat{\mathcal F\it l}_{{\mathbb G_{v_i}},R}$. Let $\mathcal S(1):= L^+_{\pv_i}\mathbb G_{\pv_i} \cdot 1\subset {\mathcal F\it l}_{\mathbb G_{v_i}}$
be the closed Schubert variety 
and $\mathcal S(1)_{R}=\mathcal S(1)\times_{\F_q} Spf R$
  then $\hat Z_{\pv_i,R}\times \mathcal S(1)_R\times\dots\times \mathcal S(1)_R$ defines a bound in $\hat{\mathcal F\it l}_{\widetilde{\mathbb G_\pv}}$ that we also denote by $Z_{\pv_i}\times \mathcal S(1)\times\dots\times \mathcal S(1)$. We have the following lemma.\\

\lem{}{\label{boundedotherway}\mbox{}
Let $\underline{\gtor}\in \san{}{\uv}(S)$ as before. The local $\mathbb G_{v_i}$-shtuka $\Gamma_{v_i}(\underline \gtor)$ is bounded by $\hat Z_{v_i}$ if and only if the local $\widetilde{\mathbb G_{v_i}}$-shtuka $L^+_\pv(\underline \gtor)$ is bounded by $\hat Z_{v_i}\times \mathcal S(1)\times\dots\times \mathcal S(1)$.
}

\prof{
We choose an \'{e}tale covering $S'$ of $S$ that trivializes $L^+_{v_i}(\gtor)$ as well as $\sigma^\star L^+_{v_i}(\gtor)$. In particular $S'$ trivializes also $\Gamma_{v_i}(\gtor)$ and $\sigma^{d\star}\Gamma_{v_i}(\gtor)$. We fix such trivializations and call them $\widetilde\alpha:L^+_{v_i}(\gtor)_{S'}\longrightarrow L^+\widetilde{\mathbb G_{v_i}}_{S'}$, $\widetilde\alpha':\sigma^\star L^+_{v_i}(\gtor)_{S'}\longrightarrow L^+\widetilde{\mathbb G_{v_i}}_{S'}$,\ \ $\alpha:\Gamma_{v_i}(\gtor)_{S'}\longrightarrow \mathbb G\times_CV(\mathfrak a_{v_i,0}) $ and $\alpha':\sigma^{d\star}\Gamma_{v_i}(\gtor)_{S'}\longrightarrow \mathbb G\times_CV(\mathfrak a_{v_i,0}) $.
Denote by $\tau_j$ the Frobenius morphism $\tau_\gtor$ restricted to $V(\mathfrak a_{{v_i},j})$ for $j=0,\dots, d-1$, where $d=deg\ v_i=[\F_{v_i}:\F_q]$. So the local shtuka $\Gamma_{v_i}(\underline \gtor)$ is given by $(\gtor\times_C V(\mathfrak a_{v_i,0}),\tau_0\circ\sigma^\star\tau_1\circ\dots\circ\sigma^{(d-1)\star}\tau_{d-1})$ and $\alpha\circ\tau_0\circ\sigma^\star\tau_1\circ\dots\circ\sigma^{(d-1)\star}\tau_{d-1}\circ \alpha'^{-1}:L{\mathbf G}_{{v_i},S'}\isom L{\mathbf G}_{{v_i},S'}$ computed in $\mathscr H^1(\F_{v_i}, L{\mathbf G}_{v_i})(S')$ defines a morphism $S'\to L{\mathbf G}_{v_i}$. Now $\Gamma_{v_i}(\underline{\gtor})$ is bounded by $Z_{v_i}$ if and only if the morphism 
$S'\times_{R_{Z_{\pv_i}}}Spf\ R\to L\mathbf G_{v_i}\times_{\F_{v_i}}Spf\ R\to \hat{\mathcal F\it l}_{\mathbb G_{v_i},R}$ factors through $Z_{{v_i},R}$. Since $\tau$ is an isomorphism outside the graphs of $s_i$, $\tau_1, \dots, \tau_{d-1}$ are isomorphisms. Hence $\sigma^\star\tau_1\circ\dots\circ\sigma^{(d-1)\star}\tau_{d-1}\circ \alpha'^{-1}$ comes from some other trivialization $\beta'^{-1}:\mathbb G\times_C V(\mathfrak a_{{v_i}, 0})\to\sigma^\star \Gamma_{v_i}({\gtor})$. This shows that $\tau_0\circ\sigma^\star\tau_1\circ\dots\circ\sigma^{(d-1)\star}\tau_{d-1}$ is bounded by $Z_{v_i}$ if and only if $\tau_0$ is bounded by $Z_{v_i}$. Now $\widetilde\alpha\circ\tau\circ(\widetilde{\alpha'})^{-1}:L\widetilde{{\mathbf G}_{v_i}}_{,S'}\isom L\widetilde{{\mathbf G}_{v_i}}_{,S'}$ computed in $\mathscr H^1(\F_q,L\widetilde{\mathbf G}_{v_i})(S')$ defines in the same way a morphism $S'\to L\widetilde{{\mathbf G}_{v_i}}$ which induces a morphism $S'\times_{R_{Z_{\pv_i}}}Spf\ R\to {\mathcal F\it l}_{\widetilde{\mathbb G}_{v_i}}\times_{\F_q}Spf\ R=\displaystyle\prod_{l\in \Z/deg {v_i}}{\mathcal F\it l}_{{\mathbb G}_{v_i}}\times_{\F_{\pv_i}}Spf\ R$. Note that 
the morphism in the $j$-th component of $\prod \mathcal F\it l_{{\mathbb G}_{v_i}}\times_{\F_{v_i}}Spf\ R$ is exactly defined by $\widetilde\alpha\circ \tau_j\circ (\widetilde\alpha')^{-1}$. Since $\tau_1,\dots, \tau_{d-1}$ are isomorphisms the morphism into the $j$-th component with $j\grg 1$ always factors through $\mathcal S(1)_R$.  This implies that $\tau$ is bounded by $ Z_{v_i}\times \mathcal S(1)\times\dots\times \mathcal S(1)$ if and only if $\tau_0$ is bounded by $Z_{v_i}$.
}
\mbox{}\\

Now we can prove:

\theo{}{\label{basechangegeneral}
Let $(C,\mathbb G, \uv, \hat Z_{\uv}, H)$ be a shtuka datum and $\pi:C\to C'$ a finite morphism of smooth, projective, geometrically irreducible curves over $\F_q$ with $\pw_i:=\pi(\pv_i)$ and $\underline\pw:=(\pw_1,\dots,\pw_n)$. Then the morphism $(\pi,id_{\pi_\star\mathbb G}):(C,\mathbb G, \uv, \hat Z_{\uv}, H)\to (C',\pi_\star\mathbb G, \underline\pw, \pi_\star \hat Z_{\underline\pv}, \pi_\star H)$ of shtuka data (see definition \ref{defshtukadata} and remark \ref{remshtukadata}) induces a finite morphism of the moduli stacks
\begin{equation*}
\pi_\star: \nabla_n^{\hat Z_{\uv},H}\mathscr H^1(C,\mathbb G)\to \nabla_n^{\pi_\star \hat Z_{\underline\pv},\pi_\star H}\mathscr H^1(C',\pi_\star\mathbb G)
.
\end{equation*}
}
\prof{
Let $S$ be an $\F_q$-scheme and $(\gtor,s_1,\dots, s_n, \tau_\gtor, \gamma)\in \nabla_n^{\hat Z_{\uv},H} \mathscr H^1(C,\mathbb G)(S)$. We describe again its image $(\gtor',s_1',\dots,s_n',\tau_{\gtor'}, \gamma')$ in $\nabla_n^{\hat Z_{\underline\pw},\pi_\star H}\mathscr H^1(C',\pi_\star\mathbb G)(S)$ to define the morphism. 
The $\pi_\star \mathbb G$-torsor $\underline \gtor'=(\gtor',s_1', \dots, s_n', \tau_{\gtor'})$ is already defined by the morphism in proposition \ref{finitemorbasechange}, but we have to prove that it lies indeed in $\nabla_n^{\hat Z_{\underline w}}\mathscr H^1(C, \pi_\star \mathbb G)(S)$. We will do this first and then define the $\pi_\star H$-level structure $\gamma'$. Since the section $s_i:S\to C$ is mapped to $s_i':=\pi\circ s_i$ and $s_i$ is required to factor through $Spf\ A_{v_i}$, we easily see with  $\pi(v_i)=w_i$ that $s_i'$ factors through $Spf\ A_{w_i}$. Furthermore this shows that $C_S'\backslash \bigcup_i\Gamma_{s_i'}\supset C'^\star\times_{\F_q}S$ and $C_S\backslash \bigcup_i\Gamma_{s_i}\supset C^\star\times_{\F_q}S$ where we use the notation 
$C^\star=C\backslash \{\pv_1,\dots,\pv_n\}$ and $C'^\star=C'\backslash\{\pw_1,\dots\pw_n\}$. It remains to show that $\tau_{\gtor'}$ is bounded by $Z_{\underline w}$ to see $\underline \gtor'\in \nabla_n^{\hat Z_{\underline w}}\mathscr H^1(C',\pi_\star\mathbb G)(S)$.

Now by assumption $\tau_\gtor$ is bounded by $Z_{v_i}$, which means by definition that the local shtukas $\Gamma_{v_i}(\underline{{\gtor}})$
are bounded by $Z_{v_i}$ for all $i$. 
By lemma \ref{boundedotherway} this is equivalent to the fact that $L^+_{v_i}(\underline{{\gtor}})$
is bounded by $Z_{v_i}\times \mathcal S(1)\times\dots\times \mathcal S(1)$. Now consider the following 2-cartesian diagram (compare \ref{beauvillelaszlo} and \cite[Lemma 5.1]{AH14}) where we set $U:=C'^\star\times_{C'}C$.
\begin{equation}\label{diagh}
\xymatrix{
\mathscr H^1(C,\mathbb G) \ar[r] \ar[d]_{\prod_{\pv}L^+_\pv} & \mathscr H_e^1(U,\mathbb G) \ar[d]_{\prod_{\pv}L_\pv} \\
\prod_{\pv\in\pi^{-1}(\underline\pw)} \mathscr H^1(\F_q, L^+\widetilde{{\mathbb G}_\pv}) \ar[r] & \prod_{\pv\in\pi^{-1}(\underline\pw)} \mathscr H^1(\F_q, L{\widetilde{{\mathbf G}_\pv}})
\ .}
\end{equation}
Here $\mathscr H_e^1 (U, \mathbb G)(S)$ is the full subcategory of $\mathscr H^1(U, \mathbb G)(S)$ consisting of those $\mathbb G$-torsors over $U_S$ that can be extended to a $\mathbb G$-torsor over $C_S$. Now the categories $\mathscr H^1(C,\mathbb G)$ and $\mathscr H^1_e(U, \mathbb G)$ are by lemma \ref{changeiso} equivalent to $\mathscr H^1(C',\mathbb \pi_\star G)$ and $\mathscr H^1_e(C'^\star,\pi_\star{\gtor})$. Furthermore the categories 
\begin{align*}
\prod_{\pv\in\pi^{-1}(\underline\pw)}\mathscr H^1(\F_q, L^+\widetilde{{\mathbb G}_\pv})&=\prod_{\pw\in\underline\pw}\mathscr H^1(\F_q, \prod_{\pv\in\pi^{-1}(\pw)} L^+\widetilde{{\mathbb G}_\pv})\mbox{\quad  and }\\
\prod_{\pv\in\pi^{-1}(\underline\pw)}\mathscr H^1(\F_q, L\widetilde{{\mathbf G}_\pv})&=\prod_{\pw\in\underline\pw}\mathscr H^1(\F_q, \prod_{\pv\in\pi^{-1}(\pw)} L\widetilde{{\mathbf 
G}_\pv})
\end{align*}
are 
by corollary \ref{pgruppenprodukt} equivalent to $\prod_{\pw\in\underline\pw}\mathscr H^1(\F_q, L^+\widetilde{{\pi_\star\mathbb G}_\pw})$ and $\prod_{\pw\in\underline\pw}\mathscr H^1(\F_q, L\widetilde{{\pi_\star\mathbf G}_\pw})$. Therefore the whole diagram \eqref{diagh} is equivalent to the diagram
\begin{equation*}
\xymatrix{
\mathscr H^1(C',\pi_\star\mathbb G) \ar[r] \ar[d]_{\prod_{\pw}L^+_\pw} & \mathscr H_e^1(C'^\star,\pi_\star\mathbb G) \ar[d]_{\prod_{\pw}L_\pw} \\
\prod_{\pw\in\underline\pw} \mathscr H^1(\F_q, L^+\widetilde{{\pi_\star\mathbb G}_\pw}) \ar[r] & \prod_{\pw\in\underline\pw} \mathscr H^1(\F_q, L{\widetilde{{\pi_\star\mathbf G}_\pw}})\ .
}
\end{equation*}

Now we choose some covering $S'$ over $S$ that trivializes $L^+_\pv{\gtor}, \sigma^\star L^+_\pv{\gtor}$ for all $\pv\in\pi^{-1}(\underline\pw)$ and fix some trivializations $\alpha_\pv:L_\pv^+({\gtor})_{S'}\to L^+{\widetilde{{\mathbb G}_\pv}_{S'}}$,  $\alpha'_\pv:\sigma^\star L_\pv^+({\gtor})_{S'}\to L^+{\widetilde{{\mathbb G}_\pv}_{S'}}$. 
\\Then $({\gtor}, \tau)$ defines a tuple $\prod_{\pw\in \underline\pw}\big(\prod_{\pv|\pw}L^+\widetilde{{\mathbb G}_{\pv,S'}},\ \prod_{\pv|\pw}\alpha_\pv\circ L_\pv(\tau|_{U})\circ{\alpha'}_\pv^{-1}\big)$ and the equivalence of the diagrams shows that it corresponds to the tuple $\prod_{\pw\in\underline\pw}\big(L^+\widetilde{{\pi_\star\mathbb G}_{\pw}}_{,S'},\ \alpha_\pw\circ L_\pw(\pi_\star\tau|_{C'^\star})\circ{\alpha'}_\pw^{-1}\big)$ defined in the same way by the shtuka $(\pi_\star{\gtor}, \pi_\star\tau)$. Here $\alpha_\pw,\alpha_\pw'$ are the trivializations corresponding to $\prod_{\pv|\pw}\alpha_\pv$ and $\prod_{\pv|\pw}{\alpha'}_\pv$.
Now choose some finite extension $R\supset \F_q\sem{\zeta}$ such that there are representatives $Z_{\pv,R}$ for all $\pv\in\underline \pv$. Using the $2$-cartesian diagram 
\begin{equation*}
\xymatrix{\displaystyle
\prod_{\pv\in\pi^{-1}(\underline\pw)}{\mathcal F\it l}_{\widetilde{{\mathbb G}_\pv}}
\ar[r] \ar[d] & \F_q
\ar[d] \\ 
\prod_{\pv\in\pi^{-1}(\underline\pw)} \mathscr H^1(\F_q, L^+\widetilde{{\mathbb G}_\pv}) 
\ar[r] & \prod_{\pv\in\pi^{-1}(\underline\pw)} \mathscr H^1(\F_q, L{\widetilde{{\mathbf G}_\pv}})
}
\end{equation*}
the tuple $(\prod_{\pv|\pw}L^+\widetilde{{\mathbb G}_{\pv,S'}}, \ \prod_{\pv|\pw}\alpha_\pv\circ L_\pv(\tau|_{U})\circ{\alpha'}_\pv^{-1})$ defines an $S'\times_{R_{\hat{\underline Z}}}Spf\ R$-valued point in $\prod_{\pv|\pw}{\mathcal F\it l}_{\widetilde{{\mathbb G}_\pv}}\times_{\F_q}Spf\ R=\prod_{\pv|\pw}\prod_{l\in\Z/deg\ \pv}{\mathcal F\it l}_{{{\mathbb G}_\pv}}\times_{\F_\pv}Spf\ R$. 
By lemma \ref{boundedotherway} the boundedness of $\underline{{\gtor}}$ at all the points $v\in (\underline v\cap \pi^{-1}(w))$ by $Z_{v}$ 
is equivalent to the boundedness of $L^+_{v}(\underline{{\gtor}})$ by 
$\prod_{l\in\Z/deg\ \pv}Z_{\pv,l}$ with $Z_{v,0}=Z_{v}$ for all $v\in \underline v$ and $Z_{\pv,l}=\mathcal S(1)_R$ for all $v\notin\underline v$ and $l\neq 0$, which means by definition, that the above $S'\times_{R_{\hat{\underline Z}}}Spf\ R$ valued point factors through $\prod_{\pv|\pw}\prod_{l\in\Z/deg\ \pv}Z_{\pv,l}$. 
The tuple $(L^+\widetilde{\pi_\star\mathbb G}_{\pw, S'}, S',\alpha_\pw\circ L_\pw(\pi_\star\tau|_{C^\star})\circ  \alpha_\pw'^{-1})$ defines in the same way a morphism $S'\times_{R_{\hat{\underline Z}}}Spf\ R\to {\mathcal F\it l}_{\widetilde{(\pi_\star\mathbb G)_\pw}}\times_{\F_q}Spf\ R=\prod_{l\in\Z/deg\ \pw}{\mathcal F\it l}_{{\pi_\star\mathbb G}_\pw}\times_{\F_\pw}Spf\ R$. 
Composing with the isomorphism $\prod_{\pv|\pw}{\mathcal F\it l}_{\widetilde{{\mathbb G}_\pv}}\simeq{\mathcal F\it l}_{\widetilde{{\pi_\star\mathbb G}_\pw}}$, the above morphism factors also through $\prod_{\pv|\pw}\prod_{l\in\Z/deg\ \pv}Z_{\pv,l}$. With lemma \ref{boundedotherway} and the definition of $\pi_\star Z_\uv$ on page \pageref{defpiz} it follows, that $\pi_\star \underline{\mathcal G}$ is bounded by $\pi_\star Z_\uv$.\\ 

Next we have to define the $\pi_\star H$-level structure $\gamma'$. We fix a geometric base point $\overline s\in S$ and we choose for all closed points $\pv\in C\backslash \uv$ a trivialization $L^+_\pv\underline{\mathcal G}_{\overline s}\simeq (L^+\widetilde{{\mathbb G}_\pv}_{,\kappa(\overline s)}, \tau=1)$, which exists by \cite[Corollary 2.9]{AH14}. This provides also trivializations 
\begin{equation*}
L^+_\pw(\pi_\star\mathcal G_{\overline s})=\prod_{\pv|\pw} L^+_\pv\mathcal G_{\overline s}\simeq \prod_{\pv|\pw}L^+{\widetilde{{\mathbb G}_\pv}}_{,\kappa(\overline s)}=L^+{\widetilde {{\pi_\star\mathbb G}_\pw}}_{,\kappa(\overline s)}.
\end{equation*}
We denote by $\underline{\mathcal L_\pv}$ and $\underline{\mathcal L_\pw}$ the shtukas $L^+_\pv\underline{\mathcal G}_{\overline s}$ and $L^+_\pw(\pi_\star\mathcal G_{\overline s})$.
Now these trivializations induce isomorphisms 
\begin{align*}
\beta:\omega^\circ_{\mathbb O^{\uv}}= \prod_{\pv\in C\backslash \uv}\mathcal T_{L^+\widetilde{{\mathbb G}_{\pv}}} \simeq & \prod_{\pv\in C\backslash\uv} \mathcal T_{\underline{\mathcal L_\pv}} =\mathcal T_{\mathcal G}\\
\pi_\star\beta:\omega^\circ_{\mathbb O^{\underline\pw}}= \prod_{\pw\in C'\backslash \underline\pw}\mathcal T_{L^+\widetilde{\pi_\star \mathbb G_\pw }}
\simeq & \prod_{\pw\in C'\backslash \underline\pw}  \mathcal T_{\underline{\mathcal L_\pw}}  =\mathcal T_{\pi_\star \mathcal G}.
\end{align*}
We write $\omega^\circ_{\mathbb A^{\uv}}:=\omega^\circ_{\mathbb O^\uv}\otimes_{\mathbb O^\uv}\mathbb A^\uv$ and ${\omega'}_{\mathbb A^{\underline\pw}}^\circ:=\pw^\circ_{\mathbb O^{\underline\pw}}\otimes_{\mathbb O^{\underline\pw}}\mathbb A^{\underline\pw}$. Now $\beta^{-1}\circ\gamma\in Aut^{\otimes}(\omega^\circ_{\mathbb A^{\uv}})$ is given by an element $g\in \mathbb G(\mathbb A^\uv)$ and the $H$-orbit of $\gamma$ is $\beta\circ gH$. Now we can use the projection $\mathbb G(\mathbb A^{\uv})\twoheadrightarrow\mathbb G(\mathbb A^{\pi^{-1}(\underline\pw)})=\pi_\star\mathbb G(\mathbb A^{\underline\pw})$ to define $\pi_\star g\in\pi_\star\mathbb G(\mathbb A^{\underline\pw})$ as the image of $g$. This corresponds to an element in $Aut({\omega'}_{\mathbb A^{\underline\pw}}^\circ)$. Therefore $\pi_\star\beta\circ \pi_\star g$
defines an element 
$\gamma'$ and consequently an $\pi_\star H$-orbit in $Isom^{\otimes}({\omega'}_{\mathbb A^{\underline\pw}}^\circ, \check {\mathcal V}_{\pi_\star \underline{\mathcal G}})$. This orbit is independent of the representative $\gamma$ since $\pi_\star H$ was defined as the image of $H$ under the above projection. 
Let $\rho\in \pi_1(S,\overline s)$ since $\gamma H\subset Isom^{\otimes}({\omega}^\circ_{\mathbb A^{\uv}}, \check {\mathcal V}_{\underline{\mathcal G}})$ is $\pi_1(S,\overline s)$ invariant, we know that there is $h\in H$ such that $\rho \gamma=\gamma h$. This defines a group homomorphism $\varphi: \pi_1(S,\overline s)\to H$ and we set
$\overline\varphi :\pi_1(S,\overline s)\to H\to \pi_\star H$. Let $\rho \in \pi_\star(S,\overline s)$ and  $\gamma'\in Isom^{\otimes}({\omega'}^\circ_{\mathbb A^{\underline \pw}}, \check {\mathcal V}_{\pi_\star \underline{\mathcal G}})$ be as above, then $\rho$ operates by $\rho \gamma =\gamma'\overline{\varphi}(\rho)$ and in particular $\pi_1(S,\overline s)\gamma'\subset \gamma'\pi_\star H$. This means that the orbit $\gamma' \pi_\star H$ is $\pi_1(S,\overline s)$ invariant and defines a level structure $\gamma'$ of $\underline \gtor'$.\\ 

After constructing this morphism, we now prove that it is representable by a scheme and finite. 
By proposition \ref{finitemorbasechange} it is clear that the morphism $\nabla_n^{\hat \hat Z_{\underline v}}\mathscr H^1(C,\mathbb G)\to \nabla_n^{\hat Z_{\underline w}}\mathscr H^1(C',\pi_\star\mathbb G)$ is finite.
Now we find some finite subscheme $D\subset C$ such that $H_D:=ker(\mathbb G(\OO^\uv)\to \mathbb G(\oo_D))$ is a subgroup of finite index in $H$. 
Then we have by \ref{hlevel} the following diagram: 
\begin{equation*} \xymatrix{
\nabla_n^{\hat Z_{\underline v}, H}\mathscr H^1(C,\mathbb G) \ar[d] 
& \nabla_n^{\hat Z_{\underline v},H_D}\mathscr H^1(C,\mathbb G) \ar[d] \ar[r]^{\sim} \ar[l]_{finite\ \acute{e}tale}
& \nabla_n^{\hat Z_{\underline v}}\mathscr H^1_D(C,\mathbb G) \ar[d] \ar[r]^{finite\ \acute{e}tale}
& \nabla_n^{\hat Z_{\underline v}}\mathscr H^1(C,\mathbb G) \ar[d]^{finite}\\
\nabla_n^{\hat Z_{\underline w},\pi_\star H}\mathscr H^1(C',\pi_\star\mathbb G)
&\nabla_n^{\hat Z_{\underline w},H_{\pi_\star D}}\mathscr H^1(C',\pi_\star\mathbb G) \ar[r]^{\sim} \ar[l]_{finite\ \acute{e}tale}
&\nabla_n^{\hat Z_{\underline w}}\mathscr H^1_{\pi_\star D}(C',\pi_\star\mathbb G) \ar[r]^{finite\ \acute{e}tale}
&\nabla_n^{\hat Z_{\underline w}}\mathscr H^1(C',\pi_\star\mathbb G)}
\end{equation*}

where the horizontal arrows are finite (and even \'{e}tale) by \cite[section\ 6]{AH13}. This implies firstly that the morphism $\nabla_n^{\hat Z_{\underline v},H_D}\mathscr H^1(C,\mathbb G) \to \nabla_n^{\hat Z_{\underline w},H_{\pi_\star D}}\mathscr H^1(C',\pi_\star\mathbb G) $ is finite and consequently that the morphism $\nabla_n^{\hat Z_{\underline v}, H}\mathscr H^1(C,\mathbb G)  \to \nabla_n^{\hat Z_{\underline w},\pi_\star H}\mathscr H^1(C',\pi_\star\mathbb G)$ is finite.
}

\subsection{Changing the Group $\mathbb G$}
\label{secchangegroup}

Now let $f:\mathbb G\to \mathbb G'$ be a morphism of smooth, affine group schemes over $C$. In this subsection we explain how this induces a morphism between the moduli stacks of $\mathbb G$-shtukas and $\mathbb G'$-shtukas. Further we prove some of its properties, depending on $f$. First of all we recall, that given a sheaf $M$ on $C_{\acute Et}$ with an action of $\mathbb G$, we can define the sheaf $M\times^\mathbb G\mathbb G'$ whose $R$-valued points are given by the set $\{(a,b)\ |\ a\in M(R), b\in \mathbb G'(R)\}/\sim$, where $(a,b)\sim (c,d)$ if and only if $(a,b)=(cg, f(g^{-1})d)$ for some $g\in \mathbb G(R)$. Actually this construction works for any sheaf of groups on any site. Now this construction is functorial for $\mathbb G$-equivariant morphisms $\varphi:M_1\to M_2$ and commutes obviously with base change. We also write $f_\star M=M\times^{\mathbb G}\mathbb G'$ and note that if $M$ is a $\mathbb G$-torsor then $f_\star M$ is a $\mathbb G'$-torsor.
With these facts we see that for a given $\mathbb G$-shtuka $(\mathcal G, s_1,\dots, s_n, \tau_{\mathcal G})$ over $S$, the tuple $(f_\star\mathcal G, s_1, \dots, s_n, f_\star\tau_{\mathcal G})$ defines a $\mathbb G'$-shtuka over $S$. Therefore we get a morphism
\begin{equation}\label{changethegroup}
\nabla_n\mathscr H^1(C,\mathbb G) \longrightarrow \nabla_n\mathscr H^1(C,\mathbb G')\qquad (\mathcal G, s_1,\dots, s_n, \tau_{\mathcal G})\mapsto (f_\star\mathcal G, s_1, \dots, s_n, f_\star\tau_{\mathcal G}).
\end{equation}

Now we want to show that this morphism also induces a morphism of these moduli stacks with additional $H$-level structure. So we fix $n$ closed points $\underline v=(v_1,\dots,v_n)$ in $C$ and let $H\subset \mathbb G(\A^{\uv})$ be an open and compact subgroup. Let further $S$ be a connected $\F_q$-scheme with a geometric base point $\overline s\in S$ and $(\underline\gtor,\gamma)=(\gtor, s_1,\dots,s_n,\tau_{\gtor}, \gamma)$ be a $\mathbb G$-shtuka over $S$ with an $H$-level structure $\gamma H$. We already mentioned that by \cite[Proposition 2.1]{Con12} $f:\mathbb G\to \mathbb G'$ induces a continuous homomorphism $f_{\A^\uv}:\mathbb G(\A^\uv)\to \mathbb G'(\A^\uv)$ (see also above definition \ref{defshtukadata}). For an open, compact subgroup $H'\subset \mathbb G'(\A^\uv)$ satisfying $f_{\A^\uv}(H)\subset H'$ we now construct an $H'$-level structure on the shtuka $f_\star\underline \gtor=(f_\star\gtor,s_1,\dots,s_n,f_\star\tau_{\gtor})$. \\
We choose for every $v\in \widetilde C=C\backslash\uv$ a trivialization $\alpha_v:L^+_v(\underline\gtor_{\overline s})\isom (L^+\widetilde{\mathbb G_{v}}_{,\overline s},1\cdot \sigma^\star)$ which exists by \cite[Proposition 2.9]{AH14}. 
Since $f_\star$ commutes with base change this induces trivializations
$f_\star \alpha_v:L^+_v((f_\star\underline \gtor)_{\overline s})\isom (L^+\widetilde{\mathbb G'_{v}}_{,\overline s},1\cdot \sigma^\star)$. We denote by $\omega^\circ_{\mathbb O^\uv}:Rep_{\OO^\uv}\mathbb G\to 
\mathfrak{Mod}_{\OO^\uv[\pi_1(S,\bar s)]}$ and $\omega'^\circ_{\mathbb O^\uv}:Rep_{\OO^\uv}\mathbb G'\to 
\mathfrak{Mod}_{\OO^\uv[\pi_1(S,\bar s)]}$ the forgetful functors and by $\underline{\mathcal L}_v$ and $\underline{\mathcal L'}_v$ the local shtukas 
$L^+_v(\underline\gtor_{\overline s})$ and $L^+_v((f_\star\underline \gtor)_{\overline s})$. Then the previous trivializations provide isomorphisms of tensor functors
\begin{align*}
\beta:\omega^\circ_{\mathbb O^{\uv}}= \prod_{\pv\in C\backslash \uv}\mathcal T_{L^+\widetilde{{\mathbb G}_{\pv}}} \isom & \prod_{\pv\in C\backslash\uv} \mathcal T_{\underline{\mathcal L_\pv}} =\mathcal T_{\underline{\mathcal G}}\\
f_\star\beta:\omega'^\circ_{\mathbb O^{\underline\pv}}= \prod_{\pv\in C\backslash \underline\pv}\mathcal T_{L^+\widetilde{\mathbb G'_\pv }}
\isom & \prod_{\pv\in C\backslash \underline\pv}  \mathcal T_{\underline{\mathcal L'_\pv}}  =\mathcal T_{f_\star\underline{\mathcal G}}.
\end{align*}

It follows that $\beta^{-1}\circ \gamma\in Aut^\otimes(\omega_{\A^\uv}^\circ)$ is given by an element $g\in \mathbb G(\A^\uv)$ and the $H$-orbit of $\gamma$ is given by $\beta\circ gH$. Now we view the image $f(g)$ of $g$ under the map $f_{\A^\uv}:\mathbb G(\A^\uv)\to \mathbb G'(\A^\uv)$ as an automorphism in 
$Aut^\otimes(\omega'^\circ_{\mathbb O^\uv})$
and define $\gamma':=f_\star\gamma := f_\star\beta\circ f(g)$. Since $f_{\A^\uv}(H)\subset H'$ the $H'$ orbit of $f(g)$ is independent of the chosen representative $\gamma$ in the orbit $\gamma H$. 
Since $\pi_1(S,\bar s)$ leaves $\gamma H$ invariant there is for all $\rho\in \pi_1(S,\bar s)$ an $h\in H$ such that $\rho\cdot \gamma=\gamma\cdot h$. This defines a group homomorphism $\varphi:\pi_1(S,\bar s)\to H$ and we set $\varphi ':\pi_1(S,\bar s)\to H\xrightarrow{f_{\A^\uv}|_H}H'$. Now $\rho\in \pi_1(S,\bar s)$ operates on $\gamma'\in Isom(\omega',\check {\mathcal V}_{\underline\gtor})$ by $\rho\cdot\gamma'=\gamma'\cdot \bar \varphi'(\rho)$. In particular $\pi_1(S,\bar s)\gamma'\subset \gamma'H'$ so that $\gamma'H'$ is $\pi_1(S,\bar s)$ invariant and defines a $H'$-level structure on $f_\star\underline \gtor$. A morphism $(\underline \gtor,\overline\gamma)$ to $(\underline{\mathcal F}, \overline \eta)$ induces naturally a morphism $(f_\star\underline \gtor, \overline \gamma')\to (f_\star\underline{\mathcal F}, \overline \eta')$ so that we get a morphism of moduli stacks 
\begin{equation}\label{changethegroup2}
\nabla_n^H\mathscr H^1(C,\mathbb G)\to \nabla_n^{H'}\mathscr H^1(C,\mathbb G')\quad (\underline \gtor, \overline\gamma)\mapsto (f_\star\underline\gtor, \overline\gamma').
\end{equation}

Next we show that this morphism behaves well with respect to boundedness conditions.
We note that for all $\pv\in\uv$ the morphism $f:\mathbb G\to \mathbb G'$ induces a morphism $L^+\mathbb G_{v_i}\to L^+\mathbb G'_{v}$ as well as a morphism $L\mathbb G_{v}\to L\mathbb G'_{v}$ and consequently also a morphism $\mathcal F\mathit l_{\mathbb G_{v}}\to \mathcal F\mathit{l}_{\mathbb G'_{v}}$.

\lem{}{\label{changegroupbound}
Let $\hat Z_{\underline v}$ be a bound in $\prod_{i-1}^n\widehat{\mathcal F\mathit l}_{\mathbb G_{v_i}}$ and $\underline{\mathcal G}$ a $\mathbb G$-shtuka over $S$ bounded by $\hat Z_{\underline v}$.
Let further $\hat Z'_{\underline v}$ be a bound in  $\prod_{i=1}^n\widehat{\mathcal F\mathit l}_{\mathbb G'_{v_i}}$ such that after choosing representatives over some DVR $R$ the morphism $\hat Z_{\underline v, R}\to \prod_{i=1}^n\widehat{\mathcal F\mathit l}_{\mathbb G'_{v_i}}$ factors through $\hat Z'_{\underline v,R}$.
Then $f_\star\underline{\mathcal G}$ is bounded by  $\hat Z'_{\underline v}$.
}
\prof{
We have to prove that for $v\in\uv$ the local shtuka $\Gamma_v(f_\star\underline{\mathcal G})$ is bounded by $\hat Z'_v$. We choose some covering $S'\to S$ with $S'/Spec\ R$ that trivializes $L^+_v\sigma^\star\mathcal G$ and $L^+_v\mathcal G$ at the same time and fix such trivializations, which we denote by
$\alpha:L^+_v\sigma^\star\mathcal G_{S'}\to L^+\mathbb G_{v,S'}$ and $\alpha':L^+_v\mathcal G_{S'}\to L^+\mathbb G_{v,S'}$. Then $f_\star \alpha$ and $f_\star\alpha'$ are trivializations of $L_v^+(f_\star\sigma^\star\mathcal G)_{S'}$ and $L_v^+(f_\star\mathcal G)_{S'}$. Now we have the automorphism $\alpha'\circ \tau^{deg\ v}\circ \alpha^{-1}:L\mathbf                                                                                                                                                                            G_{v,S'}\isom L\mathbf G_{v,S'}$ and we let $1_{S'}:S'\to L\mathbf G_{\pv,S'}$ be the unit morphism. The composition defines an $S'$ valued point $1_{S'}\circ \alpha'\circ \tau^{deg\ v}\circ \alpha^{-1}$ in $L\mathbf G_{v,S'}$. 
The composition of this point with the morphism $L\mathbb G_{v,S'}\to L\mathbb G'_{v,S'}$ induced by $f$ defines an $S'$-valued point in $L\mathbb G'_{v,S'}$. Since the diagram 
\begin{equation*}
\xymatrix{
S' \ar[r]^{1_{S'}\ \ } \ar[rd] & L\mathbf G_{v,S'} \ar[d] \ar[rrr]^{\alpha'\circ \tau^{deg\ v}\circ\alpha^{-1}} &&& L\mathbf G_{v,S'} \ar[d] \\
 & L\mathbf G'_{v,S'} \ar[rrr]^{f_\star(\alpha'\circ\tau^{deg\ v}\circ\alpha^{-1})} &&& L\mathbf G'_{v,S'}
}
\end{equation*}
commutes, this is exactly the $S'$ valued point defined by 
\begin{equation*}
f_\star(\alpha')f_\star(\tau^{deg\ v})f_\star(\alpha^{-1})=f_\star(\alpha'\tau^{deg\ v}\alpha^{-1}).
\end{equation*}
By assumption $\underline{\mathcal G}$ is bounded by $\hat Z_\uv$ and consequently the morphism $1_{S'}\circ \alpha'\circ \tau^{deg\ v}\circ \alpha^{-1}$ factors after projection to $\widehat{\mathcal F\mathit l}_{\mathbb G_v,S'}$ through $\hat Z_{v,R}$ and maps then into $\hat Z'_{v,R}$. This means exactly that $\Gamma_v(f_\star\underline{\mathcal G})$ is bounded by $\hat Z'_v$, so that $f_\star\underline{\mathcal G}$ is bounded by $\hat Z'_\uv$.

}

This lemma and the previous explanations show.

\ko{}{\label{cochangegroup}\mbox{} 
The morphism $(id,f):(C,\mathbb G,\uv,\hat Z_{\uv},H)\to (C,\mathbb G',\uv,\hat Z'_{\uv},H')$ of shtuka data 
\begin{align*}
\mbox{induces a morphism\qquad\quad }
f_\star :\quad \nabla_n^{\hat Z_\uv,H}\mathscr H^1(C,\mathbb G) &\longrightarrow \nabla_n^{\hat Z'_\uv, H'}\mathscr H^1(C,\mathbb G')\\ 
 (\mathcal G, s_1,\dots, s_n, \tau_{\mathcal G},\gamma)&\longmapsto (f_\star\mathcal G, s_1, \dots, s_n, f_\star\tau_{\mathcal G},f_\star\gamma).
\end{align*}}
\prof{Follows directly from lemma \ref{changegroupbound} and the morphism \eqref{changethegroup2} on page \pageref{changethegroup2}.
}

\mbox{}

Now we are interested in some special classes of morphisms $f:\mathbb G\to \mathbb G'$.

\subsubsection{Generic Isomorphisms of $\mathbb G$}

First of all we want to consider morphisms $f:\mathbb G\to \mathbb G'$ which are generically an isomorphism, that means $f\times id_Q:\mathbf G\isom \mathbf G'$. In this case $f:\mathbf G\to \mathbf G'$ is already an isomorphism over some open subscheme in $C$. So we fix such an $f:\mathbb G\to \mathbb G'$ and denote by $U$ the maximal open subscheme in $C$ such that $f\times id_U:\mathbb G_U\to \mathbb G'_U$ is an isomorphism and denote by $\underline w=(w_1,\dots,w_m)$ the finite set of closed points in the complement $C\backslash U$.\\
Before we come to the moduli stacks of the global $\mathbb G$-shtukas, we prove a proposition that describes the morphism $\mathscr H^1(C,\mathbb G)\to \mathscr H^1(C,\mathbb G')$. For this proposition we need the following lemma.

\lem{}{\label{quotientstack}
Let $\mathcal L'_+$ be an $L^+\widetilde{\mathbb G'_w}$ torsor over an $\F_q$-scheme $S$. Then the quotient stack $\left[\mathcal L'_+/L^+\widetilde{\mathbb G_w} \right]$ is represented by a scheme $\mathcal L'_+/L^+\widetilde{\mathbb G_w}$ over $S$ that is \'{e}tale locally on $S$ isomorphic to 
$L^+\widetilde{\mathbb G'_w}/L^+\widetilde{\mathbb G_w}$. In the case that $\mathbb G_w$ is parahoric $\mathcal L'_+/L^+\widetilde{\mathbb G_w} $ is projective.
}

\prof{
Let $\mathcal L':=\mathcal L'_+\times^{L^+\widetilde{\mathbb G'_w}}L\widetilde{\mathbf G'_w}$ be the associated $L\widetilde{\mathbf G'_w}$-torsor of $\mathcal L'_+$. By \cite[Theorem 4.4]{AH14} the quotient stack $\left[\mathcal L'/L^+\widetilde{\mathbb G_w}\right]$ is represented by an ind-quasi-projective ind-scheme $\mathcal L'/L^+\widetilde{\mathbb G_w}$ over $S$. The closed morphism $\mathcal L'_+\to \mathcal L'$ realizes $\left[\mathcal L'_+/L^+\widetilde{\mathbb G_w}\right]$
 as a closed sub-sheaf of $\mathcal L'/L^+\widetilde{\mathbb G_w}$. Since $\mathcal L'_+$ is affine over $S$ the quotient $\left[\mathcal L'_+/L^+\widetilde{\mathbb G_w}\right]$ is given by a closed subscheme in $\mathcal L'/L^+\widetilde{\mathbb G_w}$. It is clear that after passing to a covering $S'\to S$ that trivializes $\mathcal L'_+$, the scheme $\mathcal L'_+/L^+\widetilde{\mathbb G_w}$ becomes isomorphic to $L^+\widetilde{\mathbb G'_w}/L^+\widetilde{\mathbb G_w}\times_{\F_q}S'$. Since $\mathcal L'/L^+\widetilde{\mathbb G_w}$ is by \cite[Theorem 4.4]{AH14} ind-projective if $\mathbb G_w$ is parahoric, we see that the last statement about the projectivity of $\mathcal L'_+/L^+\widetilde{\mathbb G_w}$ follows. 
 }

\mbox{}\\
Now we can prove:

\prop{}{\label{changegroupgeneric}
Let $f:\mathbb G\to \mathbb G'$ be a morphism of smooth, affine group schemes over $C$, which is an isomorphism over $C\backslash 
\underline w$. Then the morphism 
\begin{equation*}
f_\star:\ \mathscr H^1(C,\mathbb G)\to \mathscr H^1(C,\mathbb G'), \quad \gtor\mapsto f_\star \gtor
\end{equation*}
is schematic and quasi-projective. \'{E}tale locally it
is relatively representable by the morphism 
\begin{equation*}
(L^+\widetilde{\mathbb G'_{w_1}}/ L^+\widetilde{\mathbb G_{w_1}})\times_{\F_q}\dots \times_{\F_q}(L^+\widetilde{\mathbb G'_{w_m}}/ L^+\widetilde{\mathbb G_{w_m}})\longrightarrow \F_{q}.
\end{equation*}

That means that for any $\F_q$-morphism $S\to \mathscr H^1(C,\mathbb G')$ there is an \'{e}tale covering $S'\to S$ such that the fiber product $S'\times_{\mathscr H^1(C,\mathbb G')}\mathscr H^1(C,\mathbb G)$ is given by $\ds S'\times_{\F_q}(\prod_{w\in\underline w}L^+\widetilde{\mathbb G'_{w}}/ L^+\widetilde{\mathbb G_{w}})$, where the product is taken over $\F_q$.
In the case that the fibers $\mathbb G_\pw$ for all $\pw\in \underline\pw$ are parahoric group schemes this morphism is projective.
}

\prof{Let $S\to \mathscr H^1(C,\mathbb G')$ be given by a $\mathbb G$-torsor $\gtor'$ over $C_S$. Let $g:T\to S$ be an $S$-scheme. Then a $T$-valued point of the fiber product $S\times_{\mathscr H^1(C,\mathbb G')}\mathscr H^1(C,\mathbb G)$ is given by a tuple $(g,\gtor,\alpha)$ where $\gtor\in \mathscr H^1(C,\mathbb G)$ and $\alpha:f_\star \gtor\isom g^\star\gtor '$. Using the theorem of Beauville-Laszlo from \ref{beauvillelaszlo} we write $\gtor=(\gtor|_{U_T},\prod_{w\in \underline w}\mathcal L_w, \varphi)$ with $U_T:=(C\backslash \underline w)\times_{\F_q} T$,\ $\mathcal L_w\in \mathscr H^1(\F_q,L^+\widetilde{\mathbb G_w})(T)$ and $\varphi=(\varphi_w)_{w\in\underline w}:\prod_{w\in\underline w}L_w(\gtor)\isom \prod_{w\in\underline w}L({\mathcal L_w})$. In the same way we write $\gtor'=(\gtor'|_{U_S},\prod_{w\in \underline w}\mathcal L'_w, \psi)$. In particular $f_\star\gtor$ is given by $(f_\star(\gtor|_{U_T}),f_{w,\star}\mathcal L_w, f_\star \varphi)$ and the isomorphism $\alpha$ is determined by  $\alpha_U:f_\star\gtor\to g^\star \gtor '$ and $\alpha_w:f_{w,\star}\mathcal L_w\to  L_w^+(g^\star\gtor ')$ satisfying
\begin{equation*}
\xymatrix{
L_w(f_\star(\gtor|_{U_T})) \ar[r]^{f_\star\varphi} \ar[d]^{L_w(\alpha_U)} & L(f_{w,\star}(\mathcal L_w)) \ar[d]^{L(\alpha_w)} \\
L_w(g^\star(\gtor'|_{U_T})) \ar[r]^{g^\star\psi} & L(g_w^\star\mathcal L'_w)
\qquad .}
\end{equation*}
Since $f|_U=id$ we have $f_\star(\gtor|_{U_T})=\gtor|_{U_T}$ and the point $(\gtor|_{U_T},\prod_{w\in\underline w}\mathcal L_w, \varphi)$ is equivalent to the point $(g^\star\gtor'|_{U_T},\prod_{w\in\underline w}\mathcal L_w, (\varphi_w\circ L_w(\alpha_U^{-1})))$ by the isomorphism $(\alpha_U^{-1},\prod id_{\mathcal L_w})$.
This shows that the category of tuples $(\gtor,\alpha)$ as above is equivalent to the category of tuples $(\mathcal L_w,\alpha_w)_{w\in \underline \p w}$ where $\mathcal L_w\in \mathscr H^1(C,L^+\widetilde{\mathbb G_w})(T)$ and $\alpha_w:f_\star \mathcal L_w\isom g^\star \mathcal L'_w$. Namely we associate with some arbitrary tuple $(\mathcal L_w,\alpha_w)_{w\in D}$ the tuple $((g^\star\gtor'|_{U_T},\mathcal L_w, \varphi),\beta)$ where $\beta|_U=id$ and $\beta_w=\alpha_w$ and $\varphi$ is uniquely determined by the condition $f_\star\varphi=\psi\circ L(\alpha_w^{-1})$. This is unique because $f\times id_{Q_w}:\mathbf G_w\to \mathbf G'_w$ is an isomorphism. \\
Now we note that the isomorphisms $\alpha_w:f_\star \mathcal L_w\isom g^\star \mathcal L'_w$ are in bijection with the $L^+\widetilde{\mathbb G_w}$ equivariant morphisms $\mathcal L_w\to g^\star \mathcal L'_w$. This shows that the tuples $(\mathcal L_w,\alpha_w)$ parametrize exactly the $T$-valued points of the quotient stack $\left[ g^\star \mathcal L_w/L^+\widetilde{\mathbb G_w}\right]$ over $\F_q$. It follows with the lemma \ref{quotientstack} that the fiber product $S\times_{\mathscr H^1(C,\mathbb G')}\mathscr H^1(C,\mathbb G)$ is given by the scheme $g^\star \mathcal L_{w_1}/L^+\widetilde{\mathbb G_{w_1}}\times\dots\times g^\star \mathcal L_{w_m}/L^+\widetilde{\mathbb G_{w_m}}$. In particular the
morphism $f_\star$ is representable and the remaining statements follow directly from the previous lemma.
}

\mbox{}\\
Now let us turn to the moduli stacks of global $\mathbb G$-shtukas. 
Let us firstly assume that $\underline w\subset \uv$ and that all the closed points $\underline \pw$ are $\F_q$-rational. In particular the group homomorphism $f:\mathbb G\to \mathbb G'$ is an isomorphism outside the fixed characteristic places $\pv_1,\dots,\pv_n$. Then we have the following theorem.

\prop{}{\label{changeparahoric}
Let $f:\mathbb G\to \mathbb G'$ be a morphism of smooth, affine group schemes over $C$, which is an isomorphism over $C\backslash 
\underline w$ with $\underline\pw\subset \uv$ and $w_i\in C(\F_q)$ for all $w_i\in\underline w$. Let $H\subset \mathbb G(\A^\uv)=\mathbb G'(\A^\uv)$ be an open, compact subgroup, let $\hat Z'_{\pv_i}$ be a bound in $\fl_{\mathbb G'_{\pv_i}}$ for all $i$ and let $\hat Z_{\pv_i}$ be the base change of $\hat Z'_{{\pv_i}}$ under the map $\fl_{\mathbb G_{\pv_i}}\to\fl_{\mathbb G'_{\pv_i}}$. Then the morphism 
\begin{equation*}
f_\star:\ \nabla_n^{\hat Z_\uv, H}\mathscr H^1(C,\mathbb G)\to \nabla_n^{\hat Z'_\uv, H} \mathscr H^1(C,\mathbb G'), \quad (\underline\gtor,\gamma H) \mapsto (f_\star \underline\gtor,\gamma H)
\end{equation*}
 is schematic and quasi-projective. \'{E}tale locally it is relatively representable by the morphism
\begin{equation*}
(L^+\widetilde{\mathbb G'_{w_1}}/ L^+\widetilde{\mathbb G_{w_1}})\times_{\F_q}\dots \times_{\F_q}(L^+\widetilde{\mathbb G'_{w_m}}/ L^+\widetilde{\mathbb G_{w_m}})\longrightarrow \F_{q}.
\end{equation*}
That means that for any $\F_q$-scheme $S$ there is an \'{e}tale covering $S'\to S$ such that the fiber product $S'\times_{\nabla_n^{\hat Z'_\uv, H}\mathscr H^1(C,\mathbb G)}\nabla_n^{\hat Z_\uv, H}\mathscr H^1(C',\mathbb G')$ is given by $\ds S'\times_{\F_q}(\prod_{w\in\underline w}L^+\mathbb G'_{w}/ L^+\mathbb G_{w})$, where the product is taken over $\F_q$. In particular $f_\star$ is a surjective morphism.
In the case that $\mathbb G$ is a parahoric Bruhat-Tits group scheme this morphism is projective.
}

\prof{
Since $f$ is an isomorphism outside $\underline \pw$, for two open subgroups $\tilde H
\subset H\subset \mathbb G(\A^{\uv})=\mathbb G'(\A^\uv)$ the diagram 
\begin{equation*}
\xymatrix{
\nabla_n^{\hat Z_\uv, \widetilde H}\mathscr H^1(C, \mathbb G) 
\ar[d]^{f_\star}
\ar[r]
&
\nabla_n^{\hat Z_\uv, H}\mathscr H^1(C, \mathbb G)
\ar[d]^{f_\star}
\\
\nabla_n^{\hat Z'_\uv, \widetilde H}\mathscr H^1(C, \mathbb  G') 
\ar[r]
&
\nabla_n^{\hat Z'_\uv,  H}\mathscr H^1(C, \mathbb  G')  
}
\end{equation*} is cartesian. In particular we can assume $H\subset\mathbb G(\OO^\uv)=\mathbb G'(\OO^\uv)$, because otherwise we can prove the theorem for the compact open subgroup $\widetilde H:=H\cap \mathbb G(\OO^\uv)$. This implies the assertions of the theorem for the group $H$ since the vertical arrows on the left and the right in the previous diagram are relatively represented by the same morphism.
Now for each $S$-valued point $(\underline\gtor,\gamma H)$ in $\nabla_n^{\hat Z_\uv, H}\mathscr H^1(C,\mathbb G)$ we find an isomorphic point $(\underline\gtor',\gamma'H)$ with $\gamma'\in Isom^{\otimes}(\omega_{\OO^\uv}^\circ, {\check {\mathcal T}}_{\underline \gtor'})$.
This is due to the fact, that we can pull back global $\mathbb G$-shtukas along quasi-isogenies of local $\mathbb G_\pv$-shtukas \cite[Theorem 5.2]{AH13} and is explained in the proof of \cite[Theorem 6.4]{AH13}. We get a morphism $\nabla_n^{\hat Z_\uv,H}\mathscr H^1(C,\mathbb G)\to \mathscr H^1(C,\mathbb G)$ sending $(\underline\gtor,\gamma H)=(\underline \gtor',\gamma'H)$ to $\gtor'$. This is the morphism $\nabla_n^{\hat Z_\uv,H}\mathscr H^1(C,\mathbb G)\to \nabla_n^{\hat Z_\uv,\mathbb G(\OO^\uv)}\mathscr H^1(C,\mathbb G)$
from \eqref{levelproj} in \ref{hlevel} composed with the morphism \eqref{leveliso} with $D=\emptyset$ in \ref{hlevel} and the natural morphism $\nabla_n\mathscr H^1(C,\mathbb G)\to \mathscr H^1(C,\mathbb G)$.
  Now using proposition \ref{changegroupgeneric} it suffices to prove that $\nabla_n^{\hat Z_{\uv},H}\mathscr H^1(C,\mathbb G)$ is given by the fiber product 
  \begin{equation*}
  \mathcal M:=\nabla_n^{\hat Z'_{\uv},H}\mathscr H^1(C,\mathbb G')\times_{\mathscr H^1(C,\mathbb G')} \mathscr H^1(C,\mathbb G).
  \end{equation*}
   There is a natural morphism $p:\nabla_n^{\hat Z_\uv, H}\mathscr H^1(C,\mathbb G)\to \mathcal M$ which sends an $S$-valued point $(\underline \gtor, \gamma H)$, where we can assume as before $\gamma\in Isom^{\otimes}(\omega_{\OO^\uv}^\circ, {\check {\mathcal T}}_{\underline \gtor})$, to $((f_\star\underline\gtor,\gamma H), \gtor, id_{f_\star\underline\gtor})$, which is well defined by \ref{cochangegroup}. We need to prove that this morphism induces an equivalence of the fibered categories. First we see that it is fully faithful. Let $(\underline\gtor_1,\gamma_1 H)$ and $(\underline\gtor_2,\gamma_2 H)$ be two $S$-valued points in $\nabla_n^{\hat Z_\uv, H}\mathscr H^1(C,\mathbb G)$, where we assume again $\gamma_1,\gamma_2\in Isom^{\otimes}(\omega_{\OO^\uv}^\circ, {\check {\mathcal T}}_{\underline \gtor})$. Let $g\in Hom((\underline\gtor_1,\gamma_1 H),(\underline\gtor_2,\gamma_2 H))$. Since $\check {\mathcal V}_g\circ\gamma_1=\gamma_2\ \mathrm{mod}\ H$ we see that $\check {\mathcal V}_g= \gamma_2\circ h\circ\gamma_1^{-1}$ for some $h\in H$, which implies that $\check {\mathcal V}_g$ already comes from a tensor isomorphism in $\mathrm{Isom}^\otimes(\check {\mathcal T}_{\gtor_1},\check {\mathcal T}_{\gtor_2})$. By \cite[Proposition 3.6]{AH14} it follows that $g$ is not only a quasi-isogeny but also a morphism of the global $\mathbb G$-shtukas $\underline \gtor_1\to \underline\gtor_2$. Therefore $Hom((\underline\gtor_1,\gamma_1 H),(\underline\gtor_2,\gamma_2 H))$ equals the morphisms of $\mathbb G$-torsors such that $g$ is a morphism of the global $\mathbb G$-shtukas $\underline\gtor_1$ and $\underline \gtor_2$ compatible with the level structure. Since $f_\star$ is an isomorphism outside of $\uv$ the latter condition is equivalent to the statement that $f_\star g$ is a morphism of $\mathbb G'$-shtukas compatible with the level structure. But this says exactly that 
  \begin{equation*}
  Hom((\underline\gtor_1,\gamma_1 H),(\underline\gtor_2,\gamma_2 H))=Hom_{\mathcal M}(((f_\star\underline\gtor_1,\gamma_1 H), \gtor_1, id_{f_\star\underline\gtor_1}),((f_\star\underline\gtor_2,\gamma_2 H), \gtor_2, id_{f_\star\underline\gtor_2})).
  \end{equation*}
  For the essential surjectivity let $((E,s_1,\dots,s_n,\tau_E,\gamma_EH),\gtor,\psi)$ be an $S$-valued point in $\mathcal M$, with $\gamma_E\in Isom^{\otimes}(\omega_{\OO^\uv}^\circ, {\check {\mathcal T}}_{\underline E})$ as before. This is isomorphic to 
  \begin{equation*}
  ((f_\star \gtor, s_1,\dots,s_n,\sigma^\star\psi\circ \tau_E\circ \psi^{-1}, \check{\mathcal T}_\psi\circ\gamma_EH),\gtor,id_{f_\star\gtor})
  \end{equation*}
   by $(\psi^{-1},id_\gtor)$. We need to show that it comes from an element 
  \begin{equation*}
  (\underline \gtor, \gamma_\gtor H)=(\gtor, s_1,\dots,s_n,\tau_\gtor,\gamma_\gtor H)\in \nabla_n^{\hat Z_{\uv}, H}\mathscr H^1(C,\mathbb G).
 \end{equation*}
   Here $s_1,\dots,s_n$ and $\gtor$ are already uniquely defined. 
  Therefore we need to define the isomorphism $\tau_\gtor:\sigma^\star\gtor|_{C_S\backslash\underset{i}{\cup}\Gamma_{s_i}}\to \gtor|_{C_S\backslash\underset{i}{\cup}\Gamma_{s_i}}$. Since all the closed points $\underline \pw$ are $\F_q$-rational $(C\backslash \underline w)_S$ is contained in $C_S\backslash\underset{i}{\cup}\Gamma_{s_i}$. Note that this is not the case if $w_i$ splits, because in this case $w_i\times_{\F_q}S$ has $deg\ w_i$ components isomorphic to S and $\Gamma_{s_i}$ surjects only to one of these. By assumption $f\times id_S:\mathbb G_S\to \mathbb G'_S$ is an isomorphism over $(C\backslash \underline w)_S$. This together with the fact that $\tau_\gtor$ has to satisfy $f_\star\tau_\gtor=\sigma^\star\psi\circ \tau_E\circ \psi^{-1}$ an the inclusion $(C\backslash \underline w)_S\subset C_S\backslash\underset{i}{\cup}\Gamma_{s_i}$ defines therefore a unique $\tau_\gtor:\sigma^\star\gtor|_{C_S\backslash\underset{i}{\cup}\Gamma_{s_i}}\to \gtor|_{C_S\backslash\underset{i}{\cup}\Gamma_{s_i}}$. 
  Now $\underline\gtor$ is a global $\mathbb G$-shtuka with $H$-level structure $\gamma_\gtor:=\check{\mathcal T}_\psi\circ \gamma_E$ that is mapped to $((\underline E,\gamma_E H), \gtor, \psi)\in \mathcal M(S)$. It just remains to prove that $\underline \gtor$ is bounded by $\hat Z_\uv$ to see that $(\underline\gtor, \gamma_\gtor)$ lies indeed in $\nabla_n^{\hat Z_{\uv}, H}\mathscr H^1(C,\mathbb G)$. Let $R$ be an extension of $A_{\pv_i}$ with representatives $\hat Z_{\pv_i,R}$ and $\hat Z'_{\pv_i,R}$ of the bounds $\hat Z_{\pv_i}$ and $\hat Z'_{\pv_i}$. We choose an \'{e}tale covering $S'$ of $S$ and trivializations $\alpha:L^+{\mathbb  G}_{\pv_i,S'}\to \Gamma_{\pv_i}(\gtor_{S'})$ and $\alpha':L^+{\mathbb G}_{\pv_i,S'}\to \Gamma_{\pv_i}(\sigma^\star\gtor_{S'})$. Then $\alpha^{-1}\circ\tau_\gtor\circ \alpha'=(f_\pv)_\star(\alpha^{-1}\circ\tau_\gtor\circ \alpha'):L{\mathbf G}_{\pv_i,S'}\to L{\mathbf G}_{\pv_i,S'}$
defines an $S'$-valued point of $L{\mathbf G}_{\pv_i}$ and hence an induced morphism $S'\to \hat\fl_{\mathbb G_{\pv_i},R}\to \hat\fl_{\mathbb G'_{\pv_i},R}$. By assumption $\underline E$ and hence $f_\star\underline\gtor$ is bounded by $\hat Z'_\uv$. This means that this morphism factors through $\hat Z'_{\pv_i,R}$ and since $\hat Z_{\pv_i}$ arises from base change it factors by the universal property of the fiber product also through $\hat Z_{\pv_i,R}$. This shows that $\underline\gtor$ is bounded by $\hat Z_{\pv_i}$ for all $\pv_i\in \uv$.}\mbox{}\\

If $(id_C,f):(C,\mathbb G,\uv,\hat Z_\uv,H)\to (C,\mathbb G',\uv,\hat Z'_\uv,H)$ is a morphism of shtuka data, where $\hat Z_\uv$ does not arise as a base change of $\hat Z'_\uv$ or if $f:\mathbb G\to \mathbb G'$ is an isomorphism outside  $\underline w$ without any conditions relating $\underline w$ to the characteristic points $\uv$ or their residue field, the morphism 
\begin{equation*}
	f_\star:\ \nabla_n^{\hat Z_\uv, H}\mathscr H^1(C,\mathbb G)\to \nabla_n^{\hat Z'_\uv, H} \mathscr H^1(C,\mathbb G'), \quad (\underline\gtor,\gamma H) \mapsto (f_\star \underline\gtor,\gamma H)
\end{equation*}
 is still representable, but in general not surjective anymore. More precisely, we have the following theorem.

\theo{}{
\label{changeparahoric2}
Let $\underline w=(w_1,\dots,w_m)$ be a finite set of closed points in $C$ and let 
\begin{equation*}
	(id_C,f):(C,\mathbb G,\uv,\hat Z_\uv,H)\to (C,\mathbb G',\uv,\hat Z'_\uv,H)
\end{equation*} be a morphism of shtuka data, where $f:\mathbb G\to \mathbb G'$ is an isomorphism over $C\backslash 
\underline w$. Then the morphism 
\begin{equation*}
f_\star:\ \nabla_n^{\hat Z_\uv, H}\mathscr H^1(C,\mathbb G)\to \nabla_n^{\hat Z'_\uv, H} \mathscr H^1(C,\mathbb G'), \quad (\underline\gtor,\gamma H) \mapsto (f_\star \underline\gtor,\gamma H)
\end{equation*}
is schematic and quasi-projective. In the case that $\mathbb G$ is a parahoric Bruhat-Tits group scheme this morphism is projective.
For any morphism $(\underline\gtor', \gamma' H):S\to \nabla_n^{\hat Z'_\uv, H} \mathscr H^1(C,\mathbb G')$ 
\begin{equation*}
\mbox{the fiber product\ \ }
S\times_{\nabla_n^{\hat Z'_\uv, H} \mathscr H^1(C,\mathbb G')}\nabla_n^{\hat Z'_\uv, H} \mathscr H^1(C,\mathbb G)\mbox{\quad is given by a   closed subscheme of}
\end{equation*}
\begin{equation*}
S\times_{\F_q}\left((L^+_{w_1}(\gtor')/ L^+\widetilde{\mathbb G_{w_1}})\times_{\F_q}\dots \times_{\F_q}(L^+_{w_m}(\gtor')/ L^+\widetilde{\mathbb G_{w_m}})\right).
\end{equation*}
If $\hat Z_\uv$ arises as a base change of $\hat Z'_\uv$ for all $\pv\in \uv$, the morphism $f_\star$ is surjective.
}

\prof{
In the case that $\hat Z_{\pv_i}$ does not arise by base change from $\hat Z_{\pv_i}$ the immersion $\hat Z_{\pv_i}\to \hat\fl_{\mathbb G_{\pv_i}}$ factors through the base change $\hat Z''_{\pv_i}:=\hat Z'_{\pv_i}\times_{\hat\fl_{\mathbb G'_{\pv_i}}} \hat\fl_{\mathbb G_{\pv_i}}$. Since $\nabla_n^{\hat Z_\uv, H}\mathscr H^1(C,\mathbb G)\to \nabla_n^{\hat Z''_\uv,H}\mathscr H^1(C,\mathbb G)$ is a closed substack we may therefore assume from the beginning that $\hat Z_{\pv_i}$ arises by base change from $\hat Z'_{\pv_i}$ for all $\pv_i\in\uv$. Furthermore we can as in the previous theorem assume that $H\subset \mathbb G(\OO^\uv)$.
Let $S\to \nabla_n^{\hat Z'_\uv, H} \mathscr H^1(C,\mathbb G')$ be an $S$-valued point given by $(\underline\gtor',\gamma'H)=(\gtor',s_1',\dots,s_n',\tau_{\gtor'}, \gamma'H)$, where we can assume as before that $\gamma'\in Isom^{\otimes}(\omega_{\OO^\uv}^\circ, {\check {\mathcal T}}_{\underline \gtor'})$. There is a 
natural morphism 
\begin{equation*}
S\times_{\nabla_n^{\hat Z'_\uv, H} \mathscr H^1(C,\mathbb G')}\nabla_n^{\hat Z_\uv, H}\mathscr H^1(C,\mathbb G)\to S\times_{ \mathscr H^1(C,\mathbb G')}\mathscr H^1(C,\mathbb G)
\end{equation*}
 sending an $T$-valued point $(g,\underline \gtor, \gamma H, \psi)$ to $(g,\gtor, \psi)$, where $g:T\to S$ is a morphism of schemes, $(\underline\gtor, \gamma H)$ is a $T$-valued point in $\nabla_n^{\hat Z_\uv, H}\mathscr H^1(C,\mathbb G)$ and $\psi:f_\star(\underline\gtor, \gamma H)\isom g^\star(\underline\gtor',\gamma'H)$ is an isomorphism of global $\mathbb G'$-shtukas. By proposition \ref{changegroupgeneric} it is now enough to show that this is a closed immersion.\\
Given a $T$-valued point $(g,\gtor, \psi)$ in $S\times_{ \mathscr H^1(C,\mathbb G')}\mathscr H^1(C,\mathbb G)$, there can be at most one $T$-valued point $(g,(\underline\gtor, \gamma  H),\psi)$ in $S\times_{\nabla_n^{\hat Z'_\uv, H} \mathscr H^1(C,\mathbb G')}\nabla_n^{\hat Z_\uv, H}\mathscr H^1(C,\mathbb G)$ with $\underline \gtor=(\gtor, s_1,\dots, s_n, \tau_\gtor)$ and $\gamma\in Isom^{\otimes}(\omega_{\OO^\uv}^\circ, {\check {\mathcal T}}_{\underline \gtor})$ mapping to $(g,(\underline\gtor, \gamma  H),\psi)$. This is because $\psi:f_\star(\underline\gtor,\gamma H)\isom g^\star(\underline\gtor',\gamma ' H)$ is an isomorphism of global $\mathbb G$-shtukas. That means namely that $s_1,\dots, s_n$ are determined by $s_1'\circ g,\dots, s_n'\circ g$, that $\gamma H$ equals $\check{\mathcal T}_{\psi^{-1}}\circ g^\star\gamma' H$ and that there is at most one $\tau_\gtor$ since over the open subset $X:=(C\backslash \underline w)_S\bigcap(C_S\backslash\underset{i}{\bigcup}\Gamma_{s_i})\subset C_S$ the isomorphism $\tau_\gtor$ is determined by $f_\star\tau_\gtor=\sigma^\star\psi\circ g^\star\tau_{\gtor'}\circ \psi$. \\
Therefore we have to answer the question if the morphism 
$\tau_\gtor\big|_X: \sigma^\star \gtor\big|_X\to \gtor\big|_X$
can be extended to $C_S\backslash\underset{i}{\bigcup}\Gamma_{s_i}$. Note that if this is possible, then the global $\mathbb G$-shtuka $\underline \gtor$ is automatically bounded by $\hat Z_\uv$ as we have seen at the end of the proof of the previous proposition \ref{changeparahoric}.\\
Let $\F$ be the compositum of all $\F_{v_i}$ with $v_i\in\underline v$ and let $v_i^{(0)}\in C_\F$ be the closed point lying over $v_i$ that equals the image of the characteristic morphism $s_i$. Then the definition $\widetilde {C_\F}:=C_\F\backslash\big(\underset{i}{\bigcup}\ v_i^{(0)}\big)$ satisfies $\widetilde {C_\F}\times_\F S=C_S\backslash\left( \underset{i}{\bigcup}\ \Gamma_{s_i}\right)$ 
Let further 
\begin{equation}\label{eqdefI}
I=\left\{w\in C_\F\ \Big|\ w|w_j\ \mathrm{for\ some\ }w_j\in \underline w,\ w\neq{v_i^{(0)}}\mathrm{\ for\ all\ }v_i\in \uv\ \right\}.
\end{equation}
 In other words that means that $I$ is determined by $\underset{i}{\bigcup}\ \Gamma_{s_i}\subset(C_\F\backslash I)\times_\F S=:C_S^I$
and $((C_\F\backslash I)\times_\F S)\backslash\left( \underset{v_i\in \underline w\cap \uv}{\bigcup}\ \Gamma_{s_i}\right)=(C\backslash \underline w)\times_{F_q}S$. The definition satisfies also the equation $(\widetilde C_\F\backslash I)\times_{\F} S=U$.
Then by the theorem of Beauville-Laszlo from \ref{beauvillelaszlo} we have 
the following cartesian diagram
\begin{equation*}
\xymatrix{
\mathscr H^1(\widetilde{C_\F},\mathbb G_\F) \ar[r] \ar[d]_{\underset{w\in I}{\prod}L^+_w} &  \mathscr H_e^1(\widetilde{C_\F}\backslash I,\mathbb G_\F) \ar[d]^{\underset{w\in I}{\prod} L_w}
\\
\underset{w\in I}{\prod} \mathscr H^1(\F, L^+\widetilde{\mathbb G_w}) \ar[r] & \underset{w\in I}{\prod} \mathscr H^1(\F, L\widetilde{\mathbf G_w})
}
\end{equation*}
which means that $\sigma^\star\gtor|_{C_S\backslash\underset{i}{\cup}\Gamma_{s_i}}$ and $\gtor|_{C_S\backslash\underset{i}{\cup}\Gamma_{s_i}}$ are given by tuples $(\sigma^\star\gtor|_U, \underset{w\in I}{\prod}L^+_w(\sigma^\star\gtor), \underset{w\in I}{\prod}id_{L_w(\sigma^\star\gtor)})$ and $(\gtor|_U, \underset{v\in I}{\prod}L^+_\pv(\gtor), \underset{v\in I}{\prod} id_{L_\pv(\gtor)})$. The morphism $\tau_\gtor\big|_U:\sigma^\star\gtor|_U\to \gtor|_U$ determines for all $w\in I$ an isomorphism $L_\pw(\tau_\gtor):L_\pw(\sigma^\star\gtor)\to L_w(\gtor)$ . The question if $\tau_\gtor$ can be extended to $C_S\backslash\underset{i}{\bigcup}\ \Gamma_{s_i}$ 
is then equivalent to the question if all the isomorphisms $L_w(\tau_\gtor)$ in $\mathscr H^1(\F,L\widetilde{\mathbb G_\pw})$ already come from an isomorphism $L^+_w(\sigma^\star\gtor)\to L_w^+(\gtor)$ in $\mathscr H^1(\F, L^+\widetilde{\mathbb G_w})$. Since $L^+\widetilde{\mathbb G_\pw}\subset \widetilde{L\mathbb G_\pw}$ is a quasi-compact closed subscheme, this is a closed condition on $T$ which shows that 
\begin{equation*}
S\times_{\nabla_n^{\hat Z'_\uv, H} \mathscr H^1(C,\mathbb G')}\nabla_n^{\hat Z_\uv, H}\mathscr H^1(C,\mathbb G)\to S\times_{ \mathscr H^1(C,\mathbb G')}\mathscr H^1(C,\mathbb G)
\end{equation*} is a closed immersion.
It rests to show that under our assumption on $\hat Z_\uv$ the morphism $f_\star$ is surjective. This is not clear yet, since the closed subscheme
\begin{equation*}
S\times_{\nabla_n^{\hat Z'_\uv, H} \mathscr H^1(C,\mathbb G')}\nabla_n^{\hat Z_\uv, H}\mathscr H^1(C,\mathbb G) \hookrightarrow S\times_{\F_q}\left((L^+_{w_1}(\gtor')/ L^+\widetilde{\mathbb G_{w_1}})\times_{\F_q}\dots \times_{\F_q}(L^+_{w_m}(\gtor')/ L^+\widetilde{\mathbb G_{w_m}})\right)
\end{equation*}
does not necessarily surject to $S$. For the proof of the surjectivity we show that for any algebraically closed field $K$ and every global $\mathbb G'$-shtuka $\underline\gtor'=(\gtor',s_1,\dots, s_n, \tau_{\gtor'})$ in $\nabla_n^{\hat Z'_\uv}\mathscr H^1(C,\mathbb G')(K)$, there is a global $\mathbb G$-shtuka $\underline \gtor=(\gtor, s_1,\dots, s_n, \tau_\gtor)$ in $\nabla_n^{\hat Z_\uv}\mathscr H^1(C,\mathbb G)(K)$ with $f_\star\underline\gtor=\underline\gtor'$. By proposition \ref{changegroupgeneric} and the fact that $K$ is algebraically closed the choice of a $\mathbb G$-torsor $\gtor$ over $C_K$ with $f_\star\gtor=\gtor'$ corresponds to an element in $\prod_{w\in \underline w}(L^+\widetilde{\mathbb G'_{w}}/L^+\widetilde{\mathbb G_w})(K)$. Now let $\F'$ be the compositum of the fields $\F_w$ for all $w\in \underline w$. For a closed point $w\in\underline w\subset C$ there are exactly $deg\ w$ different closed points in $C_\F$ lying above $w$. We denote them by $w^{(0)}, \dots, w^{(deg\ w-1)}$, where $w^{(0)}$ is a randomly chosen one and the others arise by applying successively $\sigma$ on the residue field. If $w\in\uv$ we choose $w^{(0)}$ as before to be the image of the characteristic morphism $s_i$. Now once again Beauville and Laszlo help us with the diagram
\begin{equation*}
\xymatrix{
\mathscr H^1(C_{\F'},\mathbb G'_{\F'}) \ar[r] \ar[d]^{\underset{v\in J}{\prod}L^+_v} & \mathscr H_e^1(V,\mathbb G'_{\F'}) \ar[d]^{\underset{v\in J}{\prod}L_v} \\
\underset{v\in J}{\prod} \mathscr H^1(\F', L^+\mathbb G'_v) \ar[r] & \underset{v\in J}{\prod} \mathscr H^1(\F', \mathbf G'_v)
\qquad ,}
\end{equation*}
where $J=\{v\in C_{\F'}\ |\ v|w\ \mbox{ for some }w\in\underline w\}$ and $V:=C_{\F'}\backslash J$. It allows us to identify $\gtor'$ with the tuple $(\gtor'|_{V_K}, \underset{w\in\underline w}{\prod}\underset{i=1}{\overset{deg\ w}{\prod}}L^+\mathbb G'_{w^{(i)}}, (\epsilon_{w}^{(i)})_{w^{(i)}\in J})$ where $\epsilon_w^{(i)}:L_{w^{(i)}}(\gtor')\isom L\mathbf G'_{w^{(i)}}$ already comes from an isomorphism of $L^+\mathbb G'_{w^{(i)}}$-torsors. Consequently $\sigma^\star \gtor'$ is identified with 
$(\sigma^\star\gtor'|_{V_K}, \underset{w\in\underline w}{\prod}\underset{i=1}{\overset{deg\ w}{\prod}}L^+\mathbb G'_{w^{(i)}}, (\sigma^\star\epsilon_{w}^{(i-1)})_{w^{(i)}\in J})$ where $\sigma^\star\epsilon_w^{(i-1)}:L_{w^{(i)}}(\sigma^\star\gtor')\isom L\mathbf G'_{w^{(i)}}$ is coming again from an isomorphism of $L^+\mathbb G'_{w^{(i)}}$-torsors.\\
 Note that the index $i$ is computed in $\Z/deg\ w$ so that $-1=deg\ w-1$. We use again the intuitive notation $\tau'_{w^{(i)}}:=L_{w^{(i)}}(\tau_{\gtor'}):L_{w^{(i)}}(\sigma^\star\gtor')\to L_{w^{(i)}}(\gtor')$ and define for all $w^{(i)}\in J$ the element $c_{w}^{(i)}:=\epsilon_w^{(i)}\circ \tau'_{w^{(i)}}\circ\sigma^\star(\epsilon_w^{(i-1)})^{-1}$ in $L\mathbf G_{w^{(i)}}(K)$.
The fact that $\tau_{\gtor'}$ is an isomorphism over $C_K\backslash \underset{k=1}{\overset{n}{\bigcup}}\Gamma_{s_k}$
 implies that $c_w^{(i)}$ is an element in $L^+\mathbb G'_{w^{(i)}}(K)$ for all $w^{(i)}\in J_0:=J\backslash(\underset{w\in\underline w\cap \uv}{\bigcup}w^{(0)})$. Equivalently we have $c_w^{(i)}\in L^+\mathbb G'_{w^{(i)}}(K)$ for all $w\in \underline w\cap \uv$ and $i=1,\dots,deg\ w-1$ as well as for all $w\in \underline w\backslash(\underline w\cap\uv)$ and all $i=0,\dots, deg\ w-1$. 
 We will now define the tuple $(b_w^{(i)})_{w^{(i)}\in J}\in \underset{w^{(i)}\in J}{\prod}L^+\mathbb G'_{w^{(i)}}/L^+\mathbb G_{w^{(i)}}(K)$ that will determine by proposition \ref{changegroupgeneric} the $\mathbb G$-torsor $\gtor$ over $C_K$ mapped under $f_\star:
 \mathscr H^1(C,\mathbb G)\to \mathscr H^1(C,\mathbb G')$ to $\gtor'$.
If $w\in\underline w\cap \uv$ we define
\begin{equation*}
b_w^{(0)}:=1 \mbox{ and } b_w^{(i)}:=\sigma^\star b_{w}^{(i-1)}\cdot (c_w^{(i)})^{-1}\in L^+\mathbb G'_{w^{(i)}}(K)\mbox{ for all }i=1,\dots, deg\ w-1
\end{equation*}
Now if $w\notin \uv$ we write $\widetilde c_w^{}:=c_w^{(0)}\cdot\sigma^\star (c_w^{(deg\ w-1)}) \cdot \sigma^{2\star}(c_w^{(deg\ w-2)})\cdot\dots\cdot \sigma^{(deg\ w-1)\star}(c_w^{(1)})\in L^+\mathbb G'_{w^{(0)}}(K)$ and we can choose by \cite[Corollary 2.9]{AH14}
an element $b_w^{(0)}\in L^+\mathbb G'_{w^{(0)}}(K)$ satisfying $b_w^{(0)}\cdot \widetilde c_w^{}\cdot \sigma^{deg\ w\star}(b_w^{(0)})^{-1}=1$. Additionally we define 
$b_w^{(i)}:=\sigma^\star b_{w}^{(i-1)}\cdot (c_w^{(i)})^{-1}\in L^+\mathbb G'_{w^{(i)}}(K)\mbox{ for all }i=1,\dots, deg\ w-1$,
so that we have $b_w^{(i)}\cdot c_w^{(i)}\cdot \sigma^\star (b_w^{(i-1)})^{-1}=1$. Further more we see
\begin{align*}
& b_w^{(0)}\cdot c_w^{(0)}\cdot \sigma^\star (b_w^{(deg\ w-1)})^{-1} \\
=\ & b_w^{(0)} \cdot c_w^{(0)}\cdot 
\sigma\left( \sigma^{deg\ w-1}b_w^{(0)}\cdot \sigma^{deg\ w-2}(c_w^{(1)})^{-1}\cdot \sigma^{deg\ w-3}(c_w^{(2)})^{-1}\cdot\dots\cdot c_w^{deg\ w-1} \right)^{-1}\\
=\ & b_w^{(0)}\cdot \widetilde c_w^{}\cdot \sigma ^{deg\ w\star}(b_w^{(0)})^{-1}=1\ .
\end{align*}
Now the $\mathbb G$-torsor over $C_K$, determined by the choice of $(b_w^{(i)})_{w^{(i)}\in J}\in\underset{w\in\underline w}{\prod}\underset{i=0}{\overset{deg\ w-1}{\prod}}L^+\mathbb G'_{w^{(i)}}(K)$ and lying in the pre-image of $\gtor'$ under $f_\star:\mathscr H^1(C,\mathbb G)\to \mathscr H^1(C,\mathbb G')$, is given, as described in proposition \ref{changegroupgeneric}, by 
\begin{equation*}
\gtor=\big(\gtor'|_U, \underset{w\in\underline w}{\prod}\underset{i=0}{\overset{deg\ w-1}{\prod}}L^+\mathbb G_{w^{(i)}}, (b_w^{(i)}\circ \epsilon_w^{(i)})_{w^{(i)}\in J}\big).
\end{equation*}
It lies indeed in the pre-image of $\gtor'$ since $f_\star\gtor=\big(\gtor'|_U, \underset{w\in\underline w}{\prod}\underset{i=0}{\overset{deg\ w-1}{\prod}}L^+\mathbb G'_{w^{(i)}}, (b_w^{(i)}\circ \epsilon_w^{(i)})_{w^{(i)}\in J}\big)$ is isomorphic to $\gtor'$ by $\left(id_{\gtor'}\big|_{U}, \left((b_w^{(i)})^{-1}\right)_{w^{(i)}\in J}\right)$. Now we show that there is $\tau_{\gtor}:\sigma^\star\gtor|_{C_K\backslash\underset{k}{\cup}\Gamma_{s_k}}\to\gtor|_{C_K\backslash\underset{k}{\cup}\Gamma_{s_k}} $ with $f_\star \tau_\gtor=\tau_{\gtor'}$. We set $\tau_\gtor\big|_U:=\tau_{\gtor'}$ and need to convince ourself that it extends to $C_K\backslash\underset{k}{\cup}\Gamma_{s_k}$. This is the case if and only if for all $w^{(i)}\in J_0$ the vertical right hand side morphism $b_w^{(i)}\circ \epsilon_w^{(i)}\circ L_{w^{(i)}}(\tau_\gtor)\circ \sigma^\star(\epsilon_w^{(i-1)})^{-1}\circ \sigma^\star(b_w^{(i-1)})^{-1}\in L\mathbf G_{w^{(i)}}(K)$
in the diagram
\begin{equation*}
\xymatrix{
L_{w^{(i)}}(\sigma^\star\gtor) \ar[d]_{\tau_{w^{(i)}}:=L_{w^{(i)}}(\tau_\gtor)} \ar[rr]^{\sigma^\star\epsilon_w^{(i-1)}} && L\mathbf G_{w^{(i)}} \ar[rr]^{\sigma^\star b_w^{(i-1)}} \ar[d]^{c_w^{(i)}} &&  L\mathbf G_{w^{(i)}} \ar[d] \\
L_{w^{(i)}}(\gtor) \ar[rr]_{\epsilon_w^{(i)}} && L\mathbf G_{w^{(i)}} \ar[rr]_{b_w^{(i)}} && L\mathbf G_{w^{(i)}} 
}
\end{equation*}
is given by an element in $L^+\mathbb G_{w^{(i)}}(K)$.
By construction we have $b_w^{(i)}\circ c_w^{(i)}\circ \sigma^\star b_w^{(i-1)}=1$ for all $w^{(i)}\in J_0$.
This proves $f_\star{\underline \gtor}=\underline\gtor'$ and finally the theorem.
}

\mbox{}\\
\subsubsection{Closed Subgroups of $\mathbb G'$}
Secondly we take a closer look to the case that $f:\mathbb G\to \mathbb G'$ is a closed immersion of group schemes over $C$. 
We start with the following lemma that we mainly need for theorem \ref{theounram}.

\lem{}{\label{diagmono}
Let $f:\mathbb G\to \mathbb G'$ be a closed immersion of smooth, affine group schemes over the curve $C$. Then the diagonal morphism $\Delta:\mathscr H^1(C,\mathbb G)\to\mathscr H^1(C,\mathbb G)\times_{\mathscr H^1(C,\mathbb G')} \mathscr H^1(C,\mathbb G)$ of the induced morphism $f_\star:\mathscr H^1(C,\mathbb G)\to \mathscr H^1(C,\mathbb G')$ is a monomorphism. \\The same is true for the diagonal morphism $\Delta:\nabla_n\mathscr H^1(C,\mathbb G)\to \nabla_n\mathscr H^1(C,\mathbb G)\times_{\nabla_n\mathscr H^1(C,\mathbb G')} \nabla_n\mathscr H^1(C,\mathbb G)$ of the induced morphism $f_\star:\nabla_n \mathscr H^1(C,\mathbb G)\to \nabla_n \mathscr H^1(C,\mathbb G')$.
}
\prof{
For the first diagonal morphism we have to prove that for any $\F_q$-scheme $S$ the functor $\Delta_S:\mathscr H^1(C,\mathbb G)(S)\to \mathscr H^1(C,\mathbb G)\times_{\mathscr H^1(C,\mathbb G')} \mathscr H^1(C,\mathbb G)(S)$ is fully faithful. Let $\gtor\in \mathscr H^1(C,\mathbb G)(S)$, then this functor is cleary always faithful since $\varphi\in \mathrm{Aut}(\gtor)$ is send to $(\varphi,\varphi)\in \mathrm{Aut}(\Delta(\gtor))$, where $\Delta(\gtor)=(\gtor,\gtor,id_{f_\star\gtor})$. Note that it suffices to consider $\varphi\in \mathrm{Aut}(\gtor)$ since all morphisms in $\mathscr H^1(C,\mathbb G)$ are isomorphisms. To show that $\Delta_S$ is full, let $(\varphi,\psi)\in \mathrm{Aut}(\Delta(\gtor))$ which means by definition that 
\begin{equation*}
\xymatrix{
f_\star\gtor \ar[d]_{id_{f_\star\gtor}} \ar[r]^{f_\star\varphi} & f_\star\gtor \ar[d]^{id_{f_\star\gtor}} \\ 
f_\star\gtor  \ar[r]_{f_\star\psi} & f_\star\gtor
}
\end{equation*}
commutes. Therefore we have $f_\star\varphi=f_\star\psi$ and since $f:\mathbb G\to \mathbb G'$ is a closed immersion this implies $\varphi=\psi$ and hence that $\Delta_S$ is full.\\
More precisely, to see this, one chooses a covering $U\to C_S$ that trivializes $\gtor$ so that $\varphi$ and $\psi$ correspond to morphisms $\varphi, \psi:U\to \mathbb G$ satisfying the corresponding cocycle condition. The morphisms $f_\star\varphi$ and $f_\star\psi$ correspond to the compositions $U\overset{\varphi,\psi}{\longrightarrow} \mathbb G\overset{f}{\longrightarrow}\mathbb G'$ and the equality $f_\star\varphi=f_\star\psi$ means $f\circ \varphi=f\circ \psi$. Since $f$ is a closed immersion this implies $\varphi=\psi$, which proves that the first diagonal morphism is a monomorphism. The proof for the second diagonal morphism $\Delta:\nabla_n\mathscr H^1(C,\mathbb G)\to \nabla_n\mathscr H^1(C,\mathbb G)\times_{\nabla_n\mathscr H^1(C,\mathbb G')} \nabla_n\mathscr H^1(C,\mathbb G)$ works literally in the same way.
}

\ko{}{\label{psalgspace}
The morphism $f_\star:\mathscr H^1(C,\mathbb G)\to \mathscr H^1(C,\mathbb G')$ is representable by an algebraic space. 
In particular for every $\F_q$-morphism $\gtor':S\to \mathscr H^1(C,\mathbb G')$ and the natural projection $p_S:C_S\to S$, the Weil restriction $p_{S\star}(\gtor'/\mathbb G_S)$  is an algebraic space, that equals the fiber product $S\times_{\mathscr H^1(C,\mathbb G')}\mathscr H^1(C,\mathbb G)$.
}
\prof{
Since the diagonal morphism in lemma \ref{diagmono} is a monomorphism it follows by \cite[Corollary 8.1.2]{LMB00} that $f_\star:\mathscr H^1(C,\mathbb G)\to \mathscr H^1(C,\mathbb G')$ is representable by an algebraic space. 
By definition this means that the fiber product $S\times_{\mathscr H^1(C,\mathbb G')}\mathscr H^1(C,\mathbb G)$ is an algebraic space and in particular given by a functor $(\mathbf{Sch}/S)^{op}\to\mathbf{Set}$. We show that this functor coincides with the Weil restriction functor $p_{S\star}(\gtor'/\mathbb G_S)$.
By definition a $T$-valued point of this fiber product $S\times_{\mathscr H^1(C,\mathbb G')}\mathscr H^1(C,\mathbb G)$ is given by a tuple $(g, \gtor, \alpha)$ where $g:T\to S$ is a morphism of schemes, $\gtor$ is a $\mathbb G$-torsor over $C_T$ and $\alpha$ is an isomorphism of $\mathbb G'$-torsors $f_\star\gtor \isom g^\star\gtor '$. Since isomorphisms $f_\star\gtor \isom g^\star\gtor '$ are in bijection with $\mathbb G$-equivariant morphisms $\gtor\to \gtor '$ the category of the tuples above is equivalent to the set of morphisms from $C_T$ to the quotient $\gtor '/\mathbb G_S$.\\
Since $Hom_{C_S}(C_T,[\gtor '/\mathbb G_S])=Hom_S(T,p_{S,\star}([\gtor'/\mathbb G_S]))$ by definition of the Weil restriction, the fiber product $S\times_{\mathscr H^1(C,\mathbb G')}\mathscr H^1(C,\mathbb G)$ is given by $p_{S\star}(\gtor'/\mathbb G_S)$. 
}\mbox{}\\

\theo{}{\label{theounram}
Let $f:\mathbb G\to \mathbb G'$ be a closed immersion of smooth, affine group schemes over $C$. Then the induced morphism $f_\star: \nabla_n\mathscr H^1(C, \mathbb G)\to \nabla_n\mathscr H^1(C,\mathbb G')$ is unramified and schematic.
}

\prof{We first show that $f_\star$ is unramified and then conclude that it is representable by a scheme. Let $B$ be any ring and $I\subset B$ an ideal with $I^2=0$ and $p:Spec\ \overline B:=Spec\ B/I\to Spec\ B$ the natural projection arising in a diagram of the form
\begin{equation*}
\xymatrix{
Spec\ B/I \ar[d]_{p} \ar[rr]^g && \nabla_n\mathscr H^1(C,\mathbb G) \ar[d]^{f_\star} \\
Spec\ B \ar[rr]_{g'} 
\ar@{-->}@<0.6ex>[urr]^{g_1}
\ar@<-0.6ex>@{-->}[urr]_{g_2}
&& \nabla_n\mathscr H^1(C,\mathbb G')\quad .
}
\end{equation*}

To prove that $f_\star:\nabla_n\mathscr H^1(C,\mathbb G)\to \nabla_n\mathscr H^1(C,\mathbb G')$ is unramified, we need to show that for any diagram of this kind there exists at most one dashed arrow making the diagram commutative, that means $g_1=g_2$.
This suffices since $\nabla_n\mathscr H^1(C, \mathbb G)$ and $\nabla_n\mathscr H^1(C,\mathbb G')$ are locally of ind-finite type over the notherian scheme $C^n$.
The morphism $g:Spec\ \overline B\to \nabla_n\mathscr H^1(C,\mathbb G)$ corresponds to a global $\mathbb G$-shtuka $\underline{\overline \gtor}:=(\overline\gtor, \overline{s_1},\dots, \overline{s_n},\tau_{\overline \gtor})$ over $Spec\ \overline B$, where $g_1$ and $g_2$ correspond to global $\mathbb G$-shtukas $\underline\gtor_1=(\gtor_1,s_1',\dots,s_n',\tau_{\gtor_1})$ and $\underline\gtor_2=(\gtor_2,s_1'',\dots,s_n'',\tau_{\gtor_2})$ over $Spec\ B$. The commutativity of the upper triangle means
that there are isomorphisms $\beta_1:p^\star\underline\gtor_1\isom \overline{\underline{\gtor}}$ and $\beta_2:p^\star\underline\gtor_2\isom \overline{\underline{\gtor}}$ of global $\mathbb G$-shtukas over $Spec\ \overline B$. Therefore we have to prove, that the isomorphism $\beta_2^{-1}\circ \beta_1$ arises already from an isomorphism $\underline\gtor_1\to \underline\gtor_2$ of global $\mathbb G$-shtukas over $Spec\ B$.
Furthermore we denote by $\underline \gtor':=(\gtor',s_1,\dots, s_n,\tau_{\gtor'})$ the global $\mathbb G'$-shtuka over $Spec\ B$ corresponding to $g':Spec\ B\to \nabla_n\mathscr H^1(C,\mathbb G')$. The commutativity of the lower triangle gives us isomorphisms $\alpha_1:f_\star\underline\gtor_1\isom \underline \gtor'$ and $\alpha_2:f_\star\underline\gtor_2\isom \underline \gtor'$ of global $\mathbb G'$-shtukas over $Spec\ B$
satisfying $\gamma=p^\star\alpha_2\circ f_\star \beta^{-1}_2=p^\star\alpha_1\circ f_\star\beta _1^{-1}$ where $\gamma:f_\star\underline{\overline\gtor}\isom p_\star\underline\gtor'$ is the isomorphism of global $\mathbb G'$-shtukas over $Spec\ \overline B$ coming from the commutativity of the square.

Now these isomorphisms imply directly that the paws $s_i, s_i'$ and $s_i''$ coincide for all $i$ with $1\klg i\klg n$. Although $f:\mathbb G\to \mathbb G'$ is a closed immersion it is by the following remark \ref{hnotunrami} a priori not so clear that the torsors $\gtor_1$ and $\gtor_2$ are isomorphic, but we now prove this as follows.\\
The $\mathbb G$-torsors $\gtor_1$ and $\gtor_2$ over $C_B$ come with $\mathbb G$-equivariant maps to $\gtor'$ which are induced by $\alpha_1:f_\star\gtor_1\to \gtor'$ and $\alpha_2:f_\star\gtor_2\to \gtor'$. Therefore they define two $C_B$-valued points $h_1,h_2:C_B\to\gtor'/\mathbb G_B$. In other words one can describe them as follows. Since 
\begin{equation*}
\xymatrix{
Spec\ B \ar[r]^{\gtor_1} \ar[rd]_{\gtor'} & \mathscr H^1(C,\mathbb G) \ar[d]^{f_\star} \\
& \mathscr H^1(C,\mathbb G')
}\quad
\mathrm{and}\quad
\xymatrix{
Spec\ B \ar[r]^{\gtor_2} \ar[rd]_{\gtor'} & \mathscr H^1(C,\mathbb G) \ar[d]^{f_\star} \\
& \mathscr H^1(C,\mathbb G')
} 
\end{equation*}
commute, that means $f_\star\gtor_1\simeq \gtor'\simeq f_\star\gtor_2$, the $\mathbb G$-torsors $\gtor_1$ and $\gtor_2$ induce morphisms $h_1,h_2$ from $Spec\ B$ to the fiber product $Spec\ B
\times_{\mathscr H^1(C,\mathbb G')}\mathscr H^1(C,\mathbb G)$. In corollary \ref{psalgspace} this fiber product was seen to be $p_{B\star}(\gtor'/\mathbb G_B)$. Therefore we have $Hom_{C_B}\left(C_B, (\gtor'/\mathbb G_B)\right)=Hom_B(Spec\ B, p_{B\star}(\gtor'/\mathbb G_B))$ by definition of the Weil restriction, so that $h_1$ and $h_2$ correspond consequently to morphisms $C_B\to \gtor'/\mathbb G_B$. First we show that they coincide on $\widetilde C_B:=C_B\backslash\bigcup_i\Gamma_{s_i}$.\\
The $\F_q$-Frobenius induces a morphism $j:B/I\to B,\ \overline b\mapsto b^q$ which is well defined, because $I^2=0$ and in particular $I^q=0$. We get the following commutative diagram:
\begin{equation*}
\xymatrix{
Spec\ B \ar[rrd]_{\sigma_B} \ar[rr]^j && Spec\ B/I \ar[d]_{p} \ar[rr]^{\overline\gtor} && \mathscr H^1(C,\mathbb G)\\
&& Spec\ B 
\ar@<0.7ex>[urr]^{\gtor_1}
\ar@<-0.7ex>[urr]_{\gtor_2} &&
}
\end{equation*}
which implies 
\begin{equation*}
\sigma^\star\gtor_1 = j^\star p^\star \gtor_1\underset{j^\star\beta_1}{\isom} j^\star\overline\gtor\underset{j^\star\beta_2^{-1}}{\isom} j^\star p^\star\gtor_2 = \sigma^\star\gtor_2. 
 \end{equation*}
By restricting this isomorphism to $\widetilde C_B$ and composing with $\tau_{\gtor_1}$ and $\tau_{\gtor_2}$ we get 
\begin{equation*}
\delta_0:\gtor_1\big|_{\widetilde C_B}\underset{\tau^{-1}_{\gtor_1}}{\isom}\sigma^\star \gtor_1\big|_{\widetilde C_B}\underset{j^\star\beta_1}{\isom} j^\star\overline \gtor\big|_{\widetilde C_B}\underset{j^\star\beta_2^{-1}}{\isom} \sigma^\star\gtor_2\big|_{\widetilde C_B}\underset{\tau_{\gtor_2}}{\isom} \gtor_2\big|_{\widetilde C_B}\ ,
\end{equation*} an isomorphism $\delta_0:=\tau_{\gtor_2}\circ j^\star\beta_2^{-1}\circ j^\star\beta_1\circ \tau^{-1}_{\gtor_1}$ of $\mathbb G$-torsors over $\widetilde C_B$. It satisfies 
\begin{align*}
&\alpha_2|_{\widetilde C_B}\circ f_\star\delta_0\circ \alpha_1^{-1}|_{\widetilde C_B} = \underbrace{\alpha_2\circ f_\star\tau_{\gtor_2}}_{=\tau_{\gtor'}\circ\sigma^\star\alpha_2}\circ f_\star j^\star\beta_2^{-1}\circ \underbrace{f_\star j^\star\beta_1\circ \sigma^\star\alpha_1^{-1}}_{=j^\star(f_\star\beta_1\circ p^\star\alpha_1^{-1})} \circ \tau_{\gtor'}^{-1} \\
=&\tau_{\gtor'}\circ \sigma^\star\alpha_2\circ j^\star f_\star\beta_2^{-1} \circ j^\star(f_\star\beta_2\circ p^\star\alpha_2^{-1})\circ \tau_{\gtor'}^{-1}= id_{\gtor'}|_{\widetilde C_B}.
\end{align*}
In other words $\delta_0$ is an isomorphism from $(\gtor_1|_{\widetilde C_B}, \alpha_1|_{\widetilde C_B})$ to $(\gtor_2|_{\widetilde C_B}, \alpha_2|_{\widetilde C_B})$
of $\widetilde C_B$-valued points  in $\gtor'/\mathbb G_B$. Therefore the restriction of $(h_1,h_2):C_B\to \gtor'/\mathbb G_B\times_{C_B} \gtor'/\mathbb G_B$ to the open subscheme $\widetilde C_B$ in $C_B$
 factors through the diagonal in the following diagram
\begin{equation}\label{torsorcb}
\xymatrix{
C_B  \ar[rr]^{(h_1,h_2)} \ar@{-->}[rrd]  &&  \gtor'/\mathbb G_B\times_{C_B} \gtor'/\mathbb G_B  \\
\widetilde C_B \ar[rr] \ar[u] && \gtor'\backslash \mathbb G_B \ar[u]^{\Delta}
  \qquad\quad.} 
\end{equation}
To see that $\gtor_1\simeq \gtor_2$ over $C_B$ we have to show that the morphism $(h_1,h_2)$ factors through the diagonal $\Delta$ as well. Now since $f$ is a closed immersion, the quotient $\gtor'/\mathbb G_B$ exists as a scheme by \cite[Theorem 4.C]{Ana73} and it is smooth and separated by \cite[VI$_{\rm B}$, Proposition~9.2(xii) and (x)]{SGA3}. In particular the diagonal $\Delta$ is a closed immersion. Therefore $C_B$ factors through the diagonal if 
the scheme theoretic image of $\widetilde C_B$ in $C_B$ equals $C_B$. This was proven in lemma \ref{CSdense}.
As a result of this, we conclude that $\delta_0$ extends to an isomorphism $\delta:\gtor_1\isom \gtor_2$ of $\mathbb G$-torsor over $C_B$.
The computation 
\begin{align*}
\delta_0^{-1}\circ \tau_{\gtor_2} \circ \sigma^\star \delta_0 &= \tau_{\gtor_1}\circ j^\star \beta_1^{-1}\circ \underbrace{ j^\star\beta_2 \circ \tau_{\gtor_2}^{-1}\circ \tau_{\gtor_2}\circ \sigma^\star\tau_{\gtor_2}\circ\sigma^\star j^\star \beta_2^{-1}}_{=j^\star\tau_{\overline\gtor}}\circ \sigma^\star j^\star\beta_1\circ \sigma^\star \tau_{\gtor_1}^{-1} \\
&= \tau_{\gtor_1}\circ \underbrace{ j^\star \beta_1^{-1}\circ j^\star\tau_{\overline\gtor}\circ \sigma^\star j^\star\beta_1}_{= \sigma^\star \tau_{\gtor_1}} \circ \sigma^\star \tau_{\gtor_1}^{-1} = \tau_{\gtor_1}
\end{align*}
shows that $\delta:\underline\gtor_1\isom\underline\gtor_2$ is an isomorphism of $\mathbb G$-shtukas over $B$, which finishes the proof that the morphism $f_\star:\nabla_n\mathscr H^1(C, \mathbb G)\to \nabla_n\mathscr H^1(C,\mathbb G')$ is unramified. 
\\
It rests to show that this morphism is schematic.
We have proven in lemma \ref{diagmono} that the diagonal $\Delta_{f_\star}$ of $f_\star$ is a monomorphism, which implies together with \cite[Corollary 8.1.2]{LMB00} that the morphism $f_\star$ is representable by an algebraic space. It is clear that $f_\star$ is a separated morphism, since the moduli spaces of global $\mathbb G$-shtukas are separated. Furthermore we have proven that $f_\star$ is unramified and in particular locally quasi-finite \cite[Corollaire 17.4.3]{Gro67}. All together this allows us to apply \cite[Theorem A.2]{LMB00} which states that a separated, locally quasi-finite morphism of algebraic stacks that is representable by an algebraic space is already schematic. This finishes the proof of the theorem.
}\mbox{}\\

\rem{}{
\label{hnotunrami}

Note that this is a particular property of the morphism of shtukas. Even if $f:\mathbb G\to \mathbb G'$ is a closed immersion, it is not true that $\mathscr H^1(C,\mathbb G)\to \mathscr H^1(C,\mathbb G')$ is an unramified morphism.
}

\ko{}{\label{kounramified}
Let $(id_C,f):(C,\mathbb G,\uv, \hat Z_\uv, H)\to (C,\mathbb G',\uv, \hat Z'_\uv, H')$ be a morphism of shtuka data, where $f:\mathbb G\to \mathbb G'$ is a closed immersion of smooth, affine group schemes over $C$. Then the induced morphism \begin{equation*}
f_\star: \nabla_n^{\hat Z_\uv, H}\mathscr H^1(C, \mathbb G)\to \nabla_n^{\hat Z'_\uv, H'}\mathscr H^1(C,\mathbb G')
\end{equation*}
is unramified and schematic. 
}
\prof{
We first consider the induced morphism $f_\star: \nabla_n^{\hat Z_\uv}\mathscr H^1(C, \mathbb G)\to \nabla_n^{\hat Z'_\uv}\mathscr H^1(C,\mathbb G')$ of the moduli spaces of global $\mathbb G$-shtukas without level structures. We have the following commutative diagram
\begin{equation*}
\xymatrix{
\nabla_n^{\hat Z_\uv}\mathscr H^1(C, \mathbb G)\ar[rr] \ar[d]_{f_\star}
&&
\nabla_n\mathscr H^1(C, \mathbb G)^\uv \ar[d]^{f_\star}
\\
 \nabla_n^{\hat Z'_\uv}\mathscr H^1(C,\mathbb G') \ar[rr]
&&
 \nabla_n\mathscr H^1(C,\mathbb G')^\uv.
}
\end{equation*}
The vertical arrow on the right is an unramified morphism by theorem \ref{theounram}, where the horizontal arrows are closed immersions and in particular also unramified. As a consequence the vertical arrow on the left is unramified as well. 
To prove the statement for the morphism $f_\star$ of moduli spaces of global  $\mathbb G$-shtukas with level $H$-structure, we choose similar to  
\ref{basechangegeneral}
some finite subscheme $D\subset C$ such that 
$H_D:=ker(\mathbb G(\OO^\uv)\to \mathbb G(\oo_D))$ and $H'_D:=ker(\mathbb G'(\OO^\uv)\to \mathbb G'(\oo_D))$ are subgroups of finite index in $H$ (resp. H'). Then we have by \ref{hlevel} the following commutative diagram 

\begin{equation*} \xymatrix{
\nabla_n^{\hat Z_{\underline v}, H}\mathscr H^1(C,\mathbb G) \ar[d] 
&& \nabla_n^{\hat Z_{\underline v},H_D}\mathscr H^1(C,\mathbb G) \ar@<6.6ex>[d] 
\ar[ll]_(0,6){finite}^(0,6){\acute{e}tale}
\isom \nabla_n^{\hat Z_{\underline v}}\mathscr H^1_D(C,\mathbb G) \ar@<-6.6ex>[d] \ar[rr]^(0,6){finite}_(0,6){\acute{e}tale}
&& \nabla_n^{\hat Z_{\underline v}}\mathscr H^1(C,\mathbb G) \ar[d]_{unramified}\\
\nabla_n^{\hat Z_{\uv},H'}\mathscr H^1(C,\mathbb G')
&& \nabla_n^{\hat Z_{\uv},H'_{ D}}\mathscr H^1(C,\mathbb G') 
 \ar[ll]_(0,6){finite}^(0,6){\acute{e}tale}
\isom \nabla_n^{\hat Z_{\uv}}\mathscr H^1_{D}(C,\mathbb G') \ar[rr]^(0,6){finite}_(0,6){\acute{e}tale}
&&\nabla_n^{\hat Z_{\uv}}\mathscr H^1(C,\mathbb G')}
\end{equation*}
All the horizontal arrows are \'{e}tale and in particular unramified. Furthermore we have seen that the vertical arrow on the right is unramified. It follows
that $\nabla_n^{\hat Z_{\underline v}}\mathscr H^1_D(C,\mathbb G) \to \nabla_n^{\hat Z_{\uv}}\mathscr H^1_{D}(C,\mathbb G')$ is unramified \cite[Proposition 17.3.3\ (v)]{Gro67} and finally that $f_\star: \nabla_n^{\hat Z_\uv, H}\mathscr H^1(C, \mathbb G)\to \nabla_n^{\hat Z'_\uv, H'}\mathscr H^1(C,\mathbb G')$ is unramified \cite[Proposition 17.7.7]{Gro67}. It is clear that it is also schematic, which proves the corollary.
}

\theo{}{\label{theoclosedfinite}
Let $\mathbb G$ be a parahoric Bruhat-Tits group scheme and $f:\mathbb G\to \mathbb G'$ be a closed immersion of smooth, affine group schemes and $\uv=(\pv_1,\dots,\pv_n)$ be a set of closed points in $C$. Then the induced morphism
\begin{equation*}
f_\star: \nabla_n\mathscr H^1(C, \mathbb G)^{\uv}\to \nabla_n\mathscr H^1(C,\mathbb G')^{\uv}\mbox{\quad is proper and in particular finite.}
\end{equation*}

}
\prof{
We know by theorem \ref{theounram} that this morphism is unramified and schematic and in particular locally quasi-finite. Moreover the morphism is quasi-compact. Since $\nabla_n\mathscr H^1(C,\mathbb G)\to \mathscr H^1(C,\mathbb G)$ is of ind-finite type, this follows from \cite[Theorem 2.5]{AH13} after choosing a representation $\rho:\mathbb G'\to \mathrm{GL}(\mathcal V_0)$ for some vector bundle $\mathcal V_0$ such that the quotient $\mathrm GL(\mathcal V_0)/\mathbb G$ is quasi-affine (see \cite[Proposition 2.2]{AH13}).
Therefore it suffices to prove that $f_\star$ satisfies the valuative criterion for properness to see that this morphism is proper and consequently also finite, due to the quasi-finiteness. Thus let $R$ be a complete discrete valuation ring with uniformizer $\pi$ such that its residue field ${\kappa_R}=R/\pi$ is algebraically closed and let $K=Frac(R)$ be the fraction field of $R$. Let us further denote by $K^{alg}$ an algebraic closure of $K$ and by $R^{alg}$ the integral closure of $R$ in $K^{alg}$. We need to prove that in every diagram of the form 
\begin{equation}\label{diagrrkk}
\xymatrix{
Spec\ K^{alg} \ar[rr]^{i_K} \ar[d]_{\widetilde j} && Spec\ K \ar[rr]^{g_1} \ar[d]_(0,3){j} && \nabla_n\mathscr H^1(C,\mathbb G)^\uv \ar[d]^{f_\star} \\
Spec\ R^{alg} \ar[rr]_{i_R} \ar@{-->}[rrrru] && Spec\ R \ar[rr]_{g_2}  && \nabla \mathscr H^1(C,\mathbb G')^\uv
}
\end{equation}
there exists a unique dashed arrow making the diagram commutative. 

Here $g_1,g_2,i_K,i_R,j$ and $\widetilde j$ are defined by the diagram. Choosing the closed embedding $\rho:\mathbb G'\hookrightarrow \mathrm {GL}(\mathcal V_0)$ it suffices, due to the separateness of the moduli spaces, to prove the valuative criterion for the composition $\rho_\star\circ f_\star:\nabla_n\mathscr H^1(C, \mathbb G)^{\uv}\to \nabla_n\mathscr H^1(C,\mathrm {GL}(\mathcal V_0))^{\uv}$.  Therefore we may assume that $\mathbb G'$ equals $\mathrm{GL}(\mathcal V_0)$.
We denote by $\underline\gtor=(\gtor,s_1,\dots,s_n, \tau_\gtor)$ the global $\mathbb G$-shtuka over $K$ corresponding to $g_1$ and by $\underline\gtor'=(\gtor',s_1',\dots,s_n',\tau_{\gtor'})$ the global $\mathbb G'$-shtuka over $R$ corresponding to $g_2$. Furthermore the commutativity of the square gives an isomorphism $\alpha:f_\star \underline\gtor\to j^\star\underline \gtor'$ of global $\mathbb G'$-shtukas over $K$.

Let $S:=\{v\in C\ |\ \mathbb G\times_C{\F_v}\mbox{\ is not reductive }\}\bigcup\uv$. Then $\mathbb G'/\mathbb G\times_C{(C\backslash S)}$ is by \cite[Theorem 9.4.1\ and Corollary\ 9.7.7]{Alp14} an affine scheme over $C\backslash S$ and in particular $\gtor'/ \mathbb G_R\times_{C_R}(C\backslash S)_R$ is  an affine scheme over $(C\backslash S)_R$. Now 
the $\mathbb G$-torsor $\gtor|_{(C\backslash S)_R}$ with its $\mathbb G$ equivariant morphism $\gtor\to \gtor'$ induced by $\alpha$ defines an $(C\backslash S)_K$ valued point of the quotient $\gtor'/\mathbb G_R$.
\begin{equation}\label{cohnesdiag}
\xymatrix{
(C\backslash S)_K \ar[rr]^{s} \ar[rrd] && \gtor'/ \mathbb G \times_{C_R} (C\backslash S)_R \ar[d] \\
&& (C\backslash S)_R
}
\end{equation}

Now the proof consists of several steps.
In a first step we want to show that $s$ factors through $(C\backslash S)_R$ which means that it gives a section $s_R:(C\backslash S)_R\to \gtor'/\mathbb G_R\times_{C_R}(C\backslash S)_R$ of the vertical morphism in diagram \eqref{cohnesdiag}. This morphism $\ds s_R$ corresponds to a unique $\mathbb G$-torsor $\widetilde{\mathcal E}$ over $(C\backslash S)_{R}$ together with an isomorphism $\ds\alpha_R:f_\star\widetilde {\mathcal E}\isom \gtor'|_{(C\backslash S)_R}$ satisfying $j^\star \widetilde{\mathcal E}=\gtor|_{(C\backslash S)_K}$.\\
In the second step of the proof we then show that the base change of $\widetilde {\mathcal E}$ to $R^{alg}$ extends uniquely to a $\mathbb G$-torsor over the whole relative curve $C_{R^{alg}}$. More precisely, we show that there is a $\mathbb G$-torsor $\mathcal E$ over $C_{R^{alg}}$ such that firstly the restriction $\mathcal E|_{(C\backslash S)_{R^{alg}}}$ is isomorphic to the $\mathbb G$-torsor $i_R^\star\widetilde {\mathcal E}$ over $\ds (C\backslash S)_{R^{alg}}$ and secondly $f_\star\mathcal E\simeq i_R^\star \gtor'$ and ${\widetilde j}^\star\mathcal E\simeq i_K^\star\gtor$.
Then we show in the third step that this $\mathbb G$-torsor $\mathcal E$ over $C_{R^{alg}}$ gives rise to a unique $\mathbb G$-shtuka $(\mathcal E, r_1,\dots, r_n, \tau_{\mathcal E})$ in $\nabla_n\mathscr H^1(C,\mathbb G)^\uv(R^{alg})$ making the diagram \eqref{diagrrkk} commutative. This will then finish the proof.\\

\textsl{(Step 1)}\quad We can assume that $Spec\ A= C\backslash S$ is affine by enlarging $S$ if necessary. Since we have seen that $\gtor'/\mathbb G_R\times_{C_R}(C\backslash S)_R$ is affine over $(C\backslash S)_R=Spec\ A_R:=Spec\ (A\times_{\F_q}R)$ we can set $\gtor'/\mathbb G_R\times_{C_R}Spec\ A_R=:Spec\ B$ for some ring $B$. Therefore, to prove the assertion of the first step, namely that $s$ in diagram \eqref{cohnesdiag} factors through $(C\backslash S)_R$ it is enough to show that the ring morphism $s^\star:B\to A\otimes_{\F_q}K=:A_K$ factors through $A_R$. We write $L:=Frac(A_R)$ for the function field of $C_R$ and $\mathcal O:=(A_R)_{(\pi)}\subset L$ for the localisation of $A_R$ at the prime ideal $(\pi):=ker(A_R\to A_{{\kappa_R}})$. The fact that $A_R$ is normal due to the smoothness of $C_R$ over $R$ and the fact that the prime ideal $(\pi)\subset A_R$ corresponding to the generic point of $Spec\ A_{\kappa_R}$ is of height $1$ in $A_R$, implies that $\mathcal O$ is a discrete valuation ring with uniformizer $\pi$. The normality of $A_R$ allows us also by \cite[chapter II, 6.3.A]{Har77} to write $A_R=\ds \bigcap_{\mathfrak p\subset A_R\ \mathfrak p\ of\ height\ 1}A_{R,\mathfrak p}$. For all prime ideals $\mathfrak p\subset A_R$ of height $1$ we have either $\mathfrak p=(\pi)$ or $\pi\notin\mathfrak p$. In the second case $\mathfrak p$ comes from a closed point in $Spec\ A_K$ which means $A_{R,\mathfrak p}=A_{K,\mathfrak p}$. Since $A_K=\underset{\mathfrak q\subset A_K max. ideal}{\bigcap}A_{K,\mathfrak q}$ we conclude 
\begin{equation*}
A_R=\mathcal O\cap A_K\subset L.
\end{equation*}
Due to this equation it is enough to show that the composition
$s_L:Spec\ L\xrightarrow{\eta} Spec\ A_K\xrightarrow{s|_{A_K}} Spec\ B$ of $s|_{A_K}$ with $\eta:Spec\ L\to Spec\ A_K$ factors through $Spec\ \mathcal O$. \\
The Frobenius pullback $(\sigma^\star\gtor_L,\sigma^\star\alpha_L)$ with $\sigma^\star\alpha_L:(f_\star\sigma^\star\gtor)\to \sigma^\star\gtor'_L$ gives an $L$-valued point of the quotient $\sigma^\star\gtor'/\mathbb G_R$. As before this quotient is affine over $A_R$ and given exactly by $\sigma^\star\gtor'/\mathbb G_R\times_{C_R}Spec\ A_R=Spec (B\otimes_{A_R,\sigma}A_R)\to Spec\ A_R$, where the $A_R$-algebra structure of $B\otimes_{A_R,\sigma}A_R$ is given by multiplication in the second component. This means that the $L$-valued point $(\sigma^\star\gtor_L, \sigma^\star\alpha_L)$ is given by an $A_R$-morphism $Spec\ L\to Spec(B\otimes_{A_R,\sigma}A_R)$. In other words we can describe this morphism as follows. The Frobenius $\sigma:=id_A\otimes\sigma_R:A_R\to A_R$ induces of course a morphism of the fraction field $L$ which we denote again by $\sigma:L\to L,\ \frac{a}{b}\mapsto \frac{\sigma(a)}{\sigma(b)}$ for $a,b\in A_R$. It is not the absolut $\F_q$-Frobenius. Now the composition $\sigma\circ s_L^\star:B\to L$ is not an $A_R$-linear morphism, but it induces a unique $A_R$-linear morphism $\sigma^\star s_L^\star:B\otimes_{A_R,\sigma}A_R\to L$ making the following diagram commutative.
\begin{equation*}
\xymatrix{
b \ar@{|->}[d] & B \ar[rr]^{s_L^\star} \ar[d] && L \ar[d]^{\sigma} \\
b\otimes 1  & B\otimes_{{A_R},\sigma} A_R \ar[rr]^{\sigma^\star s_L^\star} && L  
}
\end{equation*}
This morphism $\sigma^\star s_L^\star$ is the one coming from the tuple $(\sigma^\star\gtor_L, \sigma^\star\alpha_L)$.

The $\mathbb G'$-shtuka $\underline\gtor'$ is defined over $R$. In particular the restriction of $\tau_{\gtor'}$ to $Spec\ A_R$ is an isomorphism $\sigma^\star\gtor'|_{A_R}\to \gtor'|_{A_R}$ that induces an isomorphism $\overline{\tau_{\gtor'}}$ of $A_R$-algebras 
\begin{equation}\label{armorphism}
\overline{\tau_{\gtor'}}:Spec\ (B\otimes_{A_R,\sigma}A_R)\to Spec\ B.
\end{equation}
 It sends a $T$-valued point $(\mathcal E_0, \delta)$ with $\delta:f_\star\mathcal E_0\to \sigma^\star\gtor'$ to $(\mathcal E_0, \tau_{\gtor'}\circ \delta)$. We then would like to know, that the following diagram 
\begin{equation}\label{diagblb}
\xymatrix{
B \ar[rd]_{s_L^\star} \ar[rr]^{\overline{\tau^\star_{\gtor'}}} & & B\otimes_{A_R,\sigma}A_R \ar[ld]_{\sigma^\star s_L^\star} \\
& L &
}
\end{equation}
of $A_R$-morphisms commutes, which can be seen as follows.
By assumption (see diagram \eqref{diagrrkk}) the diagram 
\begin{equation}\label{diagffgg}
\xymatrix{
f_\star\sigma^\star\gtor_L \ar[rr]^{\sigma^\star\alpha_L} \ar[d]_{f_\star\tau_\gtor
} && \sigma^\star\gtor'_L \ar[d]_{\tau_{\gtor'}
} \\
f_\star\gtor_L \ar[rr]^{\alpha_L} && \gtor'_L
}
\end{equation}
is a commutative diagram of isomorphisms of $\mathbb G'$-torsors over $L$. (Actually the whole diagram is already defined over $A_K$ and the vertical arrow on the right is even defined over $A_R$.) 
Now $\sigma^\star s_L$ was corresponding to $(\sigma^\star\gtor_L, \sigma^\star\alpha_L)$, so that by the description of the morphism \eqref{armorphism} the composition $\overline{\tau_{\gtor'}}\circ \sigma^\star s_L$ corresponds to the $L$-valued point $(\sigma^\star\gtor_L, \tau_{\gtor'}\circ \sigma^\star\alpha_L)$ of $Spec\ B$. This point in the fiber category $(Spec\ B)(L)$ is by $\tau_\gtor^{-1}$ isomorphic to $(\gtor_L, \tau_{\gtor'}\circ\sigma^\star\alpha_L\circ  f_\star\tau_\gtor^{-1})$, which is by diagram \eqref{diagffgg} equal to $(\gtor_L, \alpha_L)$. Since $(\gtor_L, \alpha_L)$ is exactly the $L$-valued point $s_L$ the commutativity of diagram \eqref{diagblb} follows.\\

Now we choose a closed point $\pv\in C\backslash S$. Then we can consider the associated \'etale local $\widetilde {\mathbb G'_\pv}$-shtuka $L_\pv(\underline\gtor')=(L_\pv^+(\gtor'), \tau'_\pv:=L_\pv(\tau_{\gtor'}))$ over $R$, which arises from the formal $\hat {\mathbb G}'_\pv$-torsor $\gtor'\times_{C_R}Spf\ A_{\pv,R}$ as described in \ref{parglobloc}. Since $R$ is strictly henselian the $L^+\widetilde{\mathbb G'_\pv}$-torsor $L^+_\pv(\gtor')$ is trivial so that we choose a trivialization $\beta:L^+_\pv(\gtor')\isom L^+\widetilde{\mathbb G'_\pv}$. In particular the composition $\beta\circ\tau'_\pv\circ \sigma^\star\beta^{-1}:L^+\widetilde{\mathbb G'_\pv}\isom L^+\widetilde{\mathbb G'_\pv}$ is given by an element $b\in L^+\widetilde{\mathbb G'_\pv}(R)$ so that $\beta:L_\pv(\underline\gtor')\isom (L^+\widetilde{\mathbb G'_\pv},b)$.\\
We define $R_i:=R/\pi^{q^i}$ and $b_i\in L^+\widetilde{\mathbb G'_\pv}(R_i)$ as the image of the projection of $b$ under the map $L^+\widetilde{\mathbb G'_\pv}(R)\to L^+\widetilde{\mathbb G'_\pv}(R_i),\ b\mapsto b_i$.
Since $R_0={\kappa_R}$ is algebraically closed, there exists by \cite[Corollary 2.9]{AH14} a $c_0\in L^+\widetilde{\mathbb G'_\pv}$ with $c_0=b_0\cdot\sigma^\star c_0$. Note that $\sigma^\star c_0\in L_\pv^+(\widetilde{\mathbb G'_\pv})(R_1)$. We set inductively $c_i:=b_i\cdot \sigma^\star c_{i-1}$ for $i\grg 1$ and $c:=\underset{i\to\infty}{\lim} c_i=\underset{k\to\infty}{\lim} b\cdot \sigma^\star b\cdot \dots \cdot \sigma^{(k-1)\star}b\sigma^{k\star}c_0\in L^+_\pv(\widetilde{\mathbb G'_\pv})(R)$ which satisfies $c=b\cdot\sigma^\star c$. Replacing the trivialization $\beta$ by $c^{-1}\circ \beta$ gives therefore an isomorphism of local $\widetilde {\mathbb G'_\pv}$-shtuka $c^{-1}\cdot\beta:L_\pv(\underline\gtor')\isom (L^+_\pv{\widetilde{\mathbb G'_\pv}}, id)$ as becomes clear from the diagram
\begin{equation*}
\xymatrix{
\sigma^\star L^+_\pv(\underline\gtor') \ar[r]^{\sigma^\star\beta} \ar[d]^{\tau_\pv} & L^+_\pv{\widetilde{\mathbb G'_\pv}} \ar[r]^{\sigma^\star c} \ar[d]^{b} & L^+_\pv{\widetilde{\mathbb G'_\pv}} 
\ar[d]^{id}
\\
 L^+_\pv(\underline\gtor') \ar[r]^{\beta} & L^+_\pv{\widetilde{\mathbb G'_\pv}} \ar[r]^{c^{-1}} & L^+_\pv{\widetilde{\mathbb G'_\pv}}
}
\end{equation*}
Let $A_{\pv,R}:=A_\pv\hat\otimes_{\F_\pv}R$ and $\Gamma(A,\mathbb G'/\mathbb G)$ the ring of sections of $\mathbb G'/\mathbb G$ over $Spec\ A$. The trivializations $c^{-1}\circ \beta: L^+_\pv(\underline\gtor') \to L^+_\pv{\widetilde{\mathbb G'_\pv}}$ and $\sigma^\star (c^{-1} \circ \beta) : \sigma^\star L^+_\pv(\underline\gtor') \to L^+_\pv{\widetilde{\mathbb G'_\pv}}$ and the isomorphism $\tau_\pv$ induces after passing to the $v$-adic completion morphisms $\overline{c^{-1}\beta}$, $\overline{\sigma^\star (c^{-1}\beta)}$ and $\overline{\tau_\pv}$ as in the following diagram
\begin{equation}\label{diagggbbl}
\xymatrix{
\Gamma(A,\mathbb G'/\mathbb G)\otimes_A A_{\pv,R} 
\ar[dd]^{id} \ar[rr]^{\overline{c^{-1}\beta}}
&& B \otimes_{A_R} A_{\pv,R} \ar[dd]^{\overline{\tau_\pv^\star}} \ar[rrd]^{s_L^{\star\pv}} &&\\
&& &&  L^\pv:=L\otimes_{A_R} A_{\pv,R} \\
\Gamma(A,\mathbb G'/\mathbb G)\otimes_A A_{\pv,R}
\ar[rr]^{\overline{\sigma^\star(c^{-1}\beta)}}
 && B\otimes_{A_R,\sigma} A_{\pv,R} \ar[rru]_{\sigma^\star s_L^{\star\pv}} &&
}
\end{equation}

The right hand side of the diagram arises as the $v$-adic completion of the diagram \eqref{diagblb}, where $s_L^{\star\pv}$ and $\sigma^\star s_L^{\star\pv}$ denote the induced morphism of the completion.
Since $ord_\pi(\sigma(x))=q\cdot\mathrm{ord}_\pi(x)$ for all $x\in L^\pv$ the diagram \eqref{diagggbbl} implies 
\begin{equation*}
ord_\pi(s_L^{\star\pv}\circ \overline{c^{-1}\beta}(y))=ord_\pi(\sigma^{\star}(s_L^{\star\pv}\circ \overline{c^{-1}\beta})(y))=q\cdot ord_\pi(s_L^{\star\pv}\circ \overline{c^{-1}\beta}(y))
\end{equation*}
for $y\in \Gamma(A,\mathbb G'/\mathbb G)$. This means that $ord_\pi(s_L^{\star\pv}\circ \overline{c^{-1}\beta}(y))$ equals $0$ or $\infty$. In particular we have
\begin{equation*}
s_L^{\star\pv}\circ \overline{c^{-1}\beta}:\Gamma(A,\mathbb G'/\mathbb G)\otimes_A A_{\pv,R}\to \{x\in L^\pv\ |\ ord_\pi(x)\grg 0 \}
\end{equation*}
which implies 
\begin{equation*}
\xymatrix{
B\otimes_{A_R}A_{\pv,R} \ar[rr]^{s_L^{\star\pv}\qquad} && \{x\in L^\pv\ |\ ord_\pi(x)\grg 0 \}\qquad\qquad\mbox{} \\
B \ar[rr]^{s_L^{\star}\qquad} \ar[u] && \{x\in L^\pv\ |\ ord_\pi(x)\grg 0 \}\cap L= \mathcal O
}
\end{equation*}
This finishes the first step.\\

\textsl{(Step 2)}\ 
As we have described above, the proof of the first step gives us a $\mathbb G$-torsor $\widetilde {\mathcal E}$ over $(C\backslash S)_R$ with $j^\star\widetilde{\mathcal E}=\gtor|_{(C\backslash S)_K}$ and an isomorphism $\widetilde {\alpha_R}: f_\star\widetilde{\mathcal E}\isom \gtor'|_{(C\backslash S)_R}$. We now show that $\widetilde{\mathcal E}\times_{(C\backslash S)_R}(C\backslash S)_{R^{alg}}$ extends to a $\mathbb G$-torsor $\mathcal E$ over $C_{R^{alg}}$ with $\mathcal E\times_{C_{R^{alg}}}C_{K^{alg}}=\mathcal G_{K^{alg}}$ and $\alpha_R:f_\star\mathcal E\isom \gtor'_{R^{alg}}$. 

Now the field $\widetilde L:=Quot(A_{R^{alg}})$ has transcendence degree one over $K^{alg}$ so that its cohomological dimension equals one by \cite[\S 2.3 Th\'{e}oreme 1 and remark page 140]{Ser94}. Since $\mathbb G_L$ is reductive this implies by \cite[subsection 8.6]{BS68} that $\widetilde{\mathcal E}$ is trivial over $\widetilde L$. Therefore we can choose a finite extension $K'/K$ and a trivialization $\gamma_{L'}:\widetilde{\mathcal E}\to \mathbb G_{L'}$, where $L':=Quot(A_{R'})$ and where $R'$ is the integral closure in $K'$. We recall that we denoted by $z_\pv$ a uniformizer of $C$ at $v$, so that $A_{\pv}\hat\otimes R'=R'\sem{z_\pv}$ (do not confuse $A_\pv$ with $A$) and $L'$ is contained in $Quot(R'\sem{z_\pv})$. In particular the trivialization $\gamma_L$ implies that the $\mathbb G$-torsor $\widetilde{\mathcal E}$ is trivial over $Quot(R'\sem{z_\pv})$. This fact allows us to apply \cite[1) in Theorem 1.2]{Ans18}
to see that $\widetilde{\mathcal E}_{Quot(R'\sem{z_pv})}$ extends to a $\mathbb G$-torsor $\widetilde{\mathcal E_\pv}$ over $R'\sem{z_\pv}$. (Note that \cite{Ans18} use the notation $\mathcal O_E$ for our ring ${\kappa_R}\sem{z_\pv}$ and $z$ for a uniformizer $\pi'$ in our ring $R'$.) This corresponds by \cite[Proposition 2.4]{AH14} to a $L^+\widetilde{\mathbb G}_\pv$-torsor $L^+(\widetilde{\mathcal E_\pv})$ over $R'$ which becomes trivial after base change to the strictly henselian ring $R^{alg}$. We fix such a trivialization $\beta_v:L^+_\pv(\widetilde{\mathcal E_\pv}_{R^{alg}})\isom L^+\widetilde{\mathbb G_\pv}_{R^{alg}}$ for all $\pv\in S$. 
They induce trivializations
$L(\beta_v):L_\pv(\widetilde{\mathcal E}_{R^{alg}})\isom L\widetilde{\mathbf G_\pv}_{R^{alg}}$
and therefore isomorphisms 
\begin{equation*}
L_\pv{\widetilde{\mathcal E}_{R^{alg}}}\big/L^+\widetilde{\mathbb G_{\pv}}_{,R^{alg}}\isom \left(\fl_{\widetilde{\mathbb G_\pv}}\right)\times_{\F_q}R^{alg}=:{\fl_{\pv}}_{R^{alg}}.
\end{equation*}
An $T$-valued point $(\mathcal L^+,\delta)$ with $\delta:\mathcal L\to L_\pv\widetilde{\mathcal E}_{R^{alg}}$ is send to $(\mathcal L^+,\beta_\pv\circ\delta)$. By the theorem of Beauville-Laszlo in \ref{beauvillelaszlo} we have the following cartesian diagram
\begin{equation}\label{diagbvhhhh}
\xymatrix{
\mathscr H^1(C,\mathbb G) \ar[rr] \ar[d]^{\underset{\pv\in S}{\prod}L^+_\pv} && \mathscr H^1_e(C\backslash S,\mathbb G) \ar[d]^{\underset{\pv\in S}{\prod}L_\pv} \\
\underset{\pv\in S}{\prod}\mathscr H^1(\F_q,L^+\widetilde{\mathbb G_\pv}) 
\ar[rr]_L
&& \underset{\pv\in S}{\prod} \mathscr H^1(\F_q, L\widetilde{\mathbf G_\pv})
}.
\end{equation}
Due to this diagram the torsor $\gtor_{K^{alg}}$ corresponds to the tuple 
$\big(\gtor|_{(C\backslash S)_{K^{alg}}}, \underset{\pv\in S}{\prod}L^+_\pv(\gtor_{K^{alg}}) , (\epsilon_\pv)_{\pv\in S}\big)$  with $\epsilon_\pv=id:L(L_\pv^+(\gtor))\to L_\pv(\gtor|_{(C\backslash S)_{K^{alg}}})$. Now for all $\pv\in S$ the tuple $(L_\pv^+(\gtor_{K^{alg}}),\ (\beta_\pv\times id_{K^{alg}})\circ \epsilon_\pv)$ gives an $K^{alg}$-valued point of the affine flag variety ${\fl_{\pv}}_{R^{alg}}$. By assumption $\widetilde{\mathbb G_\pv}$ is parahoric so that ${\fl_{\pv}}_{R^{alg}}$ is ind-projective over $R^{alg}$ by \cite[Theorem A]{2Ric16}. As a consequence we can lift $(L_\pv^+(\gtor_{K^{alg}}),(\beta_\pv\times id_{K^{alg}})\circ \epsilon_\pv)$ to a unique $R^{alg}$-valued point $(\mathcal E_\pv, \beta_\pv\circ\delta_\pv)\in \fl_{\pv}(R^{alg})$ with $(\mathcal E_\pv\times_{R^{alg}}K^{alg}, (\beta_\pv\circ\delta_\pv)\times id_{K^{alg}})\simeq(L_\pv^+(\gtor_{K^{alg}}),(\beta_\pv\times id_{K^{alg}})\circ \epsilon_\pv)$. \\
In particular the tuple $(\widetilde{\mathcal E}_{(C\backslash S)_{R^{alg}}}, \prod_{\pv\in S}{\mathcal E_\pv},\delta_\pv)$ defines a unique $R^{alg}$-valued point in $\mathscr H^1(C,\mathbb G)$ given by a $\mathbb G$-torsor $\mathcal E$ over $C_{R^{alg}}$ with $\mathcal E\times_{C_{R^{alg}}}C_{K^{alg}}=\gtor_{K^{alg}}$. By diagram \eqref{diagbvhhhh} with $\mathbb G$ replaced by $\mathbb G'$ we get an isomorphism $\alpha_{R^{alg}}:f_\star\mathcal E\isom \gtor'_{R^{alg}}$. This finishes the second step.\\

\textsl{(Step 3)}\ We need to show that the $\mathbb G$-torsor $\mathcal E$ over $C_{R^{alg}}$ is part of a global $\mathbb G$-shtuka $\underline{\mathcal E}=(\mathcal E,r_1,\dots,r_n,\tau_{\mathcal E})$ defining the dashed arrow in diagram \eqref{diagrrkk}. The condition that $\alpha_{R^{alg}}:f_\star\underline {\mathcal E}\isom i_R^\star\underline\gtor'$ needs to be an isomorphism of global $\mathbb G'$-shtukas defines $r_i$ by $r_i=s_i'\circ i_R$ for all $1\klg i\klg n$. So we have to construct $\tau_{\mathcal E}$. \\
From the proof of the first step we get the commutative diagram
\begin{equation}\label{diagccc}
\xymatrix{
& (C\backslash S)_R \ar[dl]_{\sigma^\star s_R} \ar[dr]^{s_R} & \\
\sigma^\star\gtor'/\mathbb G_R\times_{C_R}(C\backslash S)_R \ar[rr]^{\overline{\tau_{\gtor'}}} && \gtor'/\mathbb G_R\times_{C_R}(C\backslash S)_R 
}
\end{equation}
We defined $(\widetilde {\mathcal E}, \widetilde{\alpha_R})$ with $\widetilde{\alpha_R}:f_\star\widetilde{\mathcal E}\to \gtor'$ to be the $(C\backslash S)_R$-valued point in $\gtor'/\mathbb G_R$ corresponding to $s_R$. Hence $(\sigma^\star\widetilde{\mathcal E}, \sigma^\star \widetilde{\alpha_R})$ corresponds to $\sigma^\star s_R$ and the composition $\overline{\tau_\gtor'}\circ \sigma^\star s_R$ corresponds to the tuple $(\sigma^\star\widetilde{\mathcal E}, \tau_{\gtor'}\circ \sigma^\star \widetilde{\alpha_R})$. The commutativity of \eqref{diagccc} means that $(\sigma^\star\widetilde{\mathcal E}, \tau_{\gtor'}\circ\sigma^\star\alpha_R)$ and $(\widetilde{\mathcal E}, \widetilde{\alpha_R})$ are isomorphic as $(C\backslash S)_R$-valued points in $\gtor'/\mathbb G_R$. This gives us therefore an isomorphism $\tau_{\widetilde{\mathcal E}}:\sigma^\star\widetilde{\mathcal E}\to \widetilde{\mathcal E}$ of $\mathbb G$-torsors over $(C\backslash S)_R$ satisfying $\tau_{\gtor'}\circ \sigma^\star \widetilde{\alpha_R}=\widetilde{\alpha_R}\circ f_\star \tau_{\widetilde {\mathcal E}}$. This defines the isomorphism $\tau_{\mathcal E}$ restricted to $(C\backslash S)_{R^{alg}}$ by $\tau_{\mathcal E}|_{(C\backslash S)_{R^{alg}}}= \tau_{\widetilde{\mathcal E}}\times_{R}id_{R^{alg}}$ and we have to 
extend it to $C_{R^{alg}}\backslash \underset{i}{\bigcup}\Gamma_{r_i}$. We know additionally by ${\mathcal E}|_{C^{K^{alg}}}=\gtor$ and $\alpha:f_\star \underline\gtor\isom j^\star\underline\gtor'_{K^{alg}}$ that $\tau_{\mathcal E}$ extends to $C_{K^{alg}}\backslash \underset{i}{\bigcup}\Gamma_{r_i}$. Therefore we only have to extend $\tau_{\mathcal E}$ at finitely many closed points of $C_{R^{alg}}\backslash \underset{i}{\bigcup}\Gamma_{r_i}$. This works similar as at the end of the proof of proposition \ref{finitemorbasechange}. 
So for $p\in C_{R^{alg}}\backslash (\underset{i}{\bigcup}\Gamma_{r_i}
\bigcup C_{K^{alg}})$ we choose an open neighboorhood  $V\subset C_{{\kappa_R}}$ with $(V\times_{\kappa_R} R^{alg})\bigcap((\underset{i}{\bigcup}\Gamma_{r_i})\bigcup (C\backslash S)_{R^{alg}})=p$. We write $\widetilde V:=V\backslash p$ so that 
$\tau_{\mathcal E}$ is defined on $\widetilde V_{R^{alg}}$ and need to be extended to $V_{R^{alg}}$. Moreover the $\mathbb G$-torsors $\sigma^\star\mathcal E|_{V_{R^{alg}}}$ and ${\mathcal E}|_{V_{R^{alg}}}$ are two $R^{alg}$-valued points in $\mathscr H^1(V,\mathbb G)(R^{alg})$ so that $\tau_{\mathcal E}|_{\widetilde{V}_{R^{alg}}}$ is an isomorphism in $\mathscr H_e^1(\widetilde V,\mathbb G)(R^{alg})$. Thanks again to Beauville and Laszlo (\ref{beauvillelaszlo}) the cartesian diagram
\begin{equation*}
\xymatrix{
\mathscr H^1(V,\mathbb G) \ar[rr] \ar[d]^{L^+_p} &&  \mathscr H^1(\widetilde V, \mathbb G) \ar[d]^{L_p} \\
\mathscr H^1({\kappa_R}, L^+\widetilde {\mathbb G_p}) \ar[rr] && \mathscr H^1({\kappa_R},L\widetilde {\mathbf G_p})
}
\end{equation*}
makes it sufficient to show that the isomorphism $L_p(\tau_{\mathcal E}):L_p(\sigma^\star{\mathcal E})\to L_p(\mathcal E)$ in $\mathscr H^1({\kappa_R}, L\widetilde{\mathbf G_p})$ comes from an isomorphism in $\mathscr H^1({\kappa_R}, L^+\widetilde{\mathbb G_p})$. After trivializing $L_p(\mathcal E)$ the morphism $L_p(\tau_{\mathcal E})$ is given by an element $h\in L\widetilde {\mathbf G_p}(R^{alg})$. Since $\tau_{\mathcal E}$ is already defined on $V_{K^{alg}}$ the pullback of $h$ to $h_K\in L\widetilde{\mathbf G_p}(K^{alg})$ is already given by an element in $L^+\widetilde{\mathbb G_p}(K^{alg})$. Since $L^+\widetilde{\mathbb G_p}\subset L\widetilde{\mathbf G_p}$ is a closed subgroup we conclude that $h$ is already an element in $L^+\widetilde{\mathbb G_p}(R^{alg})$. This shows that $\tau_{\mathcal E}$ extends uniquely to $V_{R^{alg}}$ and hence to $C_{R^{alg}}\backslash \underset{i}{\bigcup}\Gamma_{r_i}$. Hereby we found the $\mathbb G$-shtuka $\underline{\mathcal E}$ over $R^{alg}$ defining a unique dashed arrow in the diagram \eqref{diagrrkk}, which ends the proof of the theorem.
}

\ko{}{
Let $\mathbb G$ be a parahoric Bruhat-Tits group scheme and $(id_C,f):(C,\mathbb G,\uv, \hat Z_\uv, H)\to (C,\mathbb G',\uv, \hat Z'_\uv, H')$ be a morphism of shtuka data, where $f:\mathbb G\to \mathbb G'$ is a closed immersion of smooth, affine group schemes over $C$. Then the induced morphism \begin{equation*}
f_\star: \nabla_n^{\hat Z_\uv, H}\mathscr H^1(C, \mathbb G)\to \nabla_n^{\hat Z'_\uv, H'}\mathscr H^1(C,\mathbb G')
\end{equation*}
is finite.
}
\prof{
The proof of this corollary works literally in the same way as the proof of corollary \ref{kounramified} with replacing unramified by finite.
}

\rem{}{\label{remAO}
The results of this section can potentially be used 
to formulate and prove an analog of the Andr\'{e}-Oort conjecture for global $\mathbb G$-shtukas. To formulate such a conjecture one needs the notion of special points and special subvarieties. In the case of Drinfeld modular curves an analog of the Andr\'{e}-Oort conjecture has been formulated and proved in \cite{Bre05}. Later the notion of special subvarieties and the formulation of the Andr\'{e} Oort conjecture was generalized in \cite{Bre12} to the higher dimensional Drinfeld modular varieties. In the same paper this Andr\'{e}-Oort conjecture was proven in some special cases. These results were extended in \cite{Hub13}. To define Drinfeld modular varieties, one fixes a point $\infty\in C$ so that $C\setminus \infty=:Spec\ A$ is affine and $\mathcal M_A^r$ is the moduli space for Drinfeld $A$-modules of rank $r$. Now for certain finite extensions $A'\subset A$ coming from a morphism $C\to C'$ of curves, Breuer shows that there is a proper morphism $\mathcal M_A^r\to \mathcal M_{A'}^{r\cdot[A:A']}$ of moduli spaces which is also finite by \cite{HH06}. 
Then Breuer uses the image of this morphism to define special subvarieties.\\
Drinfelds modular variety $\mathcal M_A^r$ can be embedded into $\nabla_2^{\hat Z_\uv}\mathscr H^1(C,\mathrm {GL_r})$ for $n=2$ and some specific choosen bound $\hat Z_\uv$. The morphism $\mathcal M_A^r\to \mathcal M_{A'}^{r\cdot[A:A']}$ then corresponds to a morphism $\nabla_2^{\hat Z_\uv}\mathscr H^1(C,\mathrm {GL_r})\to \nabla_2^{\hat Z'_{\underline w}}\mathscr H^1(C',\mathrm {GL}_{r\cdot[C:C']})$ coming from a morphism of 
shtuka data $(C,\mathrm{GL}_r, \uv, {\hat Z}_\uv)\to (C',\mathrm {GL}_{r\cdot[C:C']},\underline w,\hat Z_{\underline w})$. So extending the coefficients for Drinfeld modules generalizes to changing the curve for global $\mathbb G$-shtukas as in subsection \ref{changeofcurve}, since we are not restricted to choose $n=2$, $\mathbb G=\mathrm {GL}_r$ or some specific bound. Moreover we have seen that additionally to changing the curve, we can also change the group scheme $\mathbb G$ as in subsection \ref{secchangegroup}. 
Although we do not know if this is precisely the correct definition it is conceivable to define a special subvariety of $\nabla_n^{\hat Z'_{\underline w}, H'}\mathscr H^1(C',\mathbb G')$ to be the image of the morphism 
\begin{equation*}
f_\star\circ\pi_\star:\nabla_n^{\hat Z_\uv,H}\mathscr H^1(C,\mathbb G)\to\nabla_n^{\hat Z'_{\underline w}, H'}\mathscr H^1(C',\mathbb G')
\end{equation*}
arising from a morphism $(\pi,f)$ of shtuka data, where $f:\pi_\star\mathbb G\hookrightarrow \mathbb G'$  is a closed immersion of (Bruhat-Tits) group schemes.  Special points in $\nabla_n^{\hat Z'_{\underline w}, H'}\mathscr H^1(C',\mathbb G')$ would then be defined to be those points which arise in the image of a morphism $\widetilde f_\star: \nabla_n^{\hat Z_{\underline w}, H}\mathscr H^1(C',\mathbb T)\to \nabla_n^{\hat Z'_{\underline w}, H'}\mathscr H^1(C',\mathbb G')^{\uv}$, where $\widetilde f:\mathbb T\to \mathbb G'$ is a closed (Bruhat-Tits) group scheme that is generically a torus in $\mathbb G'$.\\
Following this, an Andr\'{e}-Oort conjecture for global $\mathbb G$-shtukas would then say that given a set $S$ of special points, the Zariski closure of these points is a finite union of special subvarieties.
Again, this is not a precise formulation but should give an impression of the flavor of a possible statement.

}

\addcontentsline{toc}{section}{References}
\printbibliography

Paul Breutmann, CNRS,\\
Institut de Math\'ematiques de Jussieu - Paris Rive Gauche (IMJ-PRG)\\
4 place Jussieu,\ 75252 Paris Cedex 5,\ France

\end{document}